\documentclass[sn-basic]{sn-jnl}% Math and Physical Sciences Reference Style

\jyear{2024}%

\theoremstyle{thmstyleone}%
\theoremstyle{thmstyletwo}%

\theoremstyle{thmstylethree}%

\usepackage{lineno}
\usepackage{natbib}
\usepackage{graphicx}
\usepackage{amsmath}
\usepackage{amsfonts}
\usepackage{xcolor}
\usepackage{todonotes}
\usepackage{rotating}
\usepackage{cleveref}
\usepackage{multirow}
\usepackage{svg}
\usepackage{comment}
\svgsetup{inkscapelatex=false}
\usepackage{subfigure}
\usepackage{stfloats}
\usepackage{hyperref}
\modulolinenumbers[5]
\usepackage[margins]{trackchanges}

\newcommand{\unoref}{\textcolor{black}}
\newcommand{\dueref}{\textcolor{black}}

\definecolor{darkgreen}{rgb}{0.0, 0.85, 0.0}

\begin{document}
\fontsize{9.5}{12}\selectfont

\title[Optimization issues in DN]{Computational issues in Optimization for Deep networks}

\author*[1]{Corrado Coppola}\email{corrado.coppola@uniroma1.it}

\author[1]{Lorenzo Papa}\email{lorenzo.papa@uniroma1.it}

\author[2]{Marco Boresta}\email{marco.boresta@iasi.cnr.it}

\author[1]{Irene Amerini}\email{irene.amerini@uniroma1.it}

\author[1,2]{Laura Palagi}\email{laura.palagi@uniroma1.it}

\affil[1]{Sapienza University of Rome, Department of Computer, Control, and Management Engineering Antonio Ruberti - Via Ariosto 25 - Roma}

\affil[2]{Istituto di Analisi dei Sistemi ed Informatica Antonio Ruberti - CNR - Via dei Taurini - Roma}

\abstract{The paper aims to investigate relevant computational issues of deep neural network {architectures} with an eye to the interaction between the optimization algorithm and the classification performance. In particular, we aim to analyze the behaviour of state-of-the-art optimization algorithms   in relationship to  {their} hyperparameters setting  {in order to} detect robustness with respect to the choice of a certain starting point in ending on different local solutions.
We conduct extensive computational experiments using nine open-source   optimization algorithms to train {deep} Convolutional Neural Network { architectures on  {an} image multi-class classification task. Precisely, we}
consider several architectures by changing the number of layers and neurons per layer,  {in order to} evaluate the impact of  {different width} and depth {structures} on the computational optimization performance.}


\maketitle

\section{Introduction}
\label{sec:Intro}
One of the key areas of artificial intelligence is \unoref{supervised} machine learning (ML), which   {involves the development of algorithms and models capable of learning \unoref{a model} on a given set of samples and making predictions or decisions based on previously unseen data.
In ML community this property is commonly known as generalization capability}. 
In most real-world applications, as well as in the case analyzed in this paper, the system is a neural network trained by minimizing a differentiable loss function measuring the dissimilarity between some target values and the values returned by the network itself.
The training process consists in the minimization of a loss function with respect to the network weights {and biases}, which can result in a complex and large-scale optimization problem. 

The crucial role played by optimization algorithms  in   {machine learning}, acknowledged since the birth of this research field, has been widely and deeply discussed in the literature, both from an {operations research} and a  {computer science} perspective. 
Training a supervised ML model involves, indeed, addressing the optimization problem (\cite{gambella2021optimization}) of minimizing a function, which measures the dissimilarity between predicted and correct values. 
Simple examples are the different type of regression (\cite{lewis2015applied,lavalley2008logistic,ranstam2018lasso}), the neural network optimization (\cite{goodfellow2014qualitatively}), the decision trees (\cite{carrizosa2021mathematical,rokach2010classification,bertsimas2017optimal,buntine2020learning}), support vector machines (  \cite{steinwart2008support,tatsumi2014support,suthaharan2016support,pisner2020support}).
While the underlying idea of stochastic gradient-like methods, proposed by \cite{robbins1951stochastic}, dates back to the 1950s, the research community has deeply investigated its theoretical and computational properties and a vast amount of new algorithms have been developed in the past decades.

 Despite all the studies that have already been carried out {in this field}, {which will discuss further,} to the best of our knowledge, only a few have tried to answer some relevant computational questions encountered when using optimization methods to train deep \unoref{neural} networks (\unoref{DNNs}). 
In this paper, we point out and address the following issues in solving training optimization problems  {to assess} the influence of optimization algorithm settings and architectural choices on generalization performances. 

\begin{itemize}
\setlength \itemsep{0.5em}
\item Convergence to 
local versus global minimizers and the effect of %whether 
the quality \unoref{(in terms of training loss) of the solution }on  the generalization performances;
\item effect of the non-monotone behaviour of mini-batch methods with respect to traditional batch methods (L-BFGS) on the computational performances;
\item different role of the starting point and regions of attraction on L-BFGS than on mini-batch algorithms;
\item importance of the optimization algorithm's hyperparameters tuning on the optimization and the generalization performances;
% \item how the overall performances depend on the network architecture and whether hyperparameters tuned on a specific problem are robust to the increase of the network in depth and in width
\item \unoref{the 
robustness of hyperparameters setting tuned on a specific architecture and dataset by modifying the number of layers and neurons and the datasets}
%robustness of selected hyperparameters settings {both} on a specific problem increasing the depth and the width of the network {and on different problems}.
\end{itemize}

% \subsection{Our contributions}
% \label{subsec:our_contr}

In this paper, we discuss and try to answer some questions \unoref{regarding the aforementioned issues}.
\unoref{We conduct extensive computational experiments} to enforce our main claims:

\begin{itemize}
\setlength \itemsep{0.5em}
\item[i)] generalization performances can be influenced by the  solution found in the training process. Local minima can be very different from each other and result in very different test performances;
\item[ii)] traditional batch methods, like L-BFGS, are less efficient and also {more sensitive} to the starting point
than mini-batch online algorithms;
%\item[iii)] the optimizers' hyperparameters (such as learning rate, momentum coefficient, etc.) play an important role in the final generalization performances and in the computational efficiency;
%\item[iii)] \unoref{some of the algorithms widely used in ML (FTLR, Adadelta and Adagrad)  cannot find good solutions on our experimental testbed, regardless the initialization seed and the hyperparameters setting; }
\item[iii)] \unoref{hyperparameters tuned on a specific baseline problem, namely a given baseline architecture trained on an instance of a class of problems, can achieve better generalization performance than the default ones even on different problems, changing either the architecture and/or the instance in the given class.}
\end{itemize}

\dueref{To the aim above, we consider the task of training convolutional neural networks (CNNs) for an image classification task.}
We use {three} open-source datasets to \unoref{carry out our} experiments. % \unoref{on training} convolutional neural network{s} (CNN) for an image classification task.
We train the networks using \dueref{nine} optimization algorithms implemented in open-source state-of-the-art libraries for optimization and ML.
We show that not all the algorithms reach a neighbourhood of a global optimum, getting stuck in local minima. \unoref{In particular, FTLR, Adadelta and Adagrad  cannot find good solutions on our experimental testbed, regardless the initialization seed and the hyperparameters setting.}
We also notice that test performance,  {i.e.,} the classification test accuracy, is remarkably higher when a good approximation of the global solution is reached and that better solutions can be achieved by carefully choosing the optimization hyperparameters setting. 
We carry out a thorough computational analysis to assess the robustness of the tuned hyper-parameters configuration on a baseline problem \unoref{(the image classification task on the open-source dataset UC Merced (\cite{uc_merced}) using a customized deep convolutional neural network)} with respect to architectural changes of the network and to new datasets for image classifications \unoref{and we find that the hyperparameters tuned on the baseline problem give often better out-of-sample performance than the default settings even on different image classification datasets.}
\unoref{Notice that we define a problem as a couple dataset-network, e.g., UC Merced-Baseline architecture.}
%\todo[inline]{dare la definizione di problema come coppia architettura - dataset}
The paper is organized as follows.
In \Cref{sec:RelWor} we discuss some relevant literature \unoref{highlighting both the importance of what has already been produced by the ML research community and the novelty of our contributions.}
In \Cref{sec:NetArc} we describe the network architecture, while in \Cref{sec:OptProb}, we formalize the optimization problem behind the image classification task, mathematically describing the convolution operation performed by the network layers.
In \Cref{sec:OptAlg} we briefly describe each of the nine different open-source algorithms we have tested on our task. 
In \Cref{sec:dataset} we describe the composition of  {the three open-source datasets}, and in \Cref{sec:impdet} implementation details are reported.  
{In} \Cref{sec:Res}, we describe in detail the computational tests we have carried out on  {different} networks {and datasets}.  
We present our conclusions in \Cref{sec:Conc}.

\section{Related literature}
\label{sec:RelWor}
Several attempts have been made in the scientific literature to {address the main issues discussed in this paper, both in the form of  {a} survey and in the form of a comparative analysis and computational study.} 
For instance, in order to understand how different types of data and tuning of algorithm parameters affected performances, \cite{lim2000comparison} carried out a thorough comparative analysis of nearly all the algorithms available at the time for classification tasks.
As {machine learning} and, in particular, deep learning, gained steadily growing interest in the community, this comparative analysis methodology became a standard framework applied to specific methods and neural architectures.
More recently, some other specific surveys have been produced, in particular comparing the behaviour of different optimization algorithms on image classification tasks \citep{dogo2018comparative,kandel2020comparative,haji2021comparison}, but they are mostly focused on mere computational aspects \dueref{rather than to performance in respect of the ML task}.
\cite{pouyanfar2018survey} and  \cite{BRAIEK2020110542} provided methodological surveys on different approaches to ML problems. In the same years, \dueref{ algorithms used in   {machine learning} have been widely studied also from an   optimization perspective. } 
%optimization algorithms used in   {machine learning} have been widely studied also from an   {operations research} perspective. 
\cite{bottou2018optimization} studied different first-order optimization algorithms applied to large-scale   {machine learning} problems, while \cite{baumann2019comparative} produced a thorough comparative analysis of first-order methods in a  {machine learning} framework and traditional combinatorial methods on the same classification tasks. 
Following the increased need for a high-level overview, some other papers on first-order methods have been published by \cite{lan2020first}, which provides a detailed survey on stochastic optimization algorithms, and by \cite{sun2019survey}, which  {compares} from a theoretical perspective the main advantages and drawbacks of some of the most used methods in   {machine learning.} 

{Some recent literature (see e.g. \citep{palagi2019global}) also discusses the role of global optimization in the training of neural networks, as well as the problem of hyper-parameters optimization.}
The role of global optimization in ML is also strictly linked to the emerging practice in ML of perfect interpolation, i.e., of training a model to fit the dataset perfectly.
Advanced studies in this direction have been carried out over the past years, starting from \cite{zhang2016understanding} and \cite{zhang2021understanding}, who reconsidered the classical bias-variance trade-off, remarking that most of the state-of-the-art neural models, especially in the field of image classification, are trained to reach close-to-zero training error, i.e., a global minimizer of the loss function.
\cite{sun2019optimization} investigates the problem of choosing the best initialization of parameters and the  {best-performing} algorithms for a given dataset, namely Global Optimization of the Network {framework}.
% {Sun} is one of the {first} to have discussed the importance of reaching a globally optimal solution when training a neural network, emphasizing how finding a global optimum leads to better test performances.
other computational studies (\cite{advani2020high,spigler2019jamming,geiger2019jamming}) enforce the idea that larger or more trained \unoref{(i.e., trained for a larger number of epochs)} models also generalize better.

%\cite{belkin2019reconciling} refers to the improved performance of over-parametrized models as the double descent phenomenon. 

Another particularly valuable work for our research is \cite{im2016empirical}, where a loss function projection mechanism is used to discuss how different algorithms can have remarkably different  {performances} on the same problem. 
{Eventually, the issue of having plenty of local minimizers, some of which are better than  {others} in the sense that they lead to better test performances \unoref{(which is, indeed, the main point of our claim i))}, has been actively addressed both from a theoretical and from a computational perspective.
\cite{ding2022suboptimal} provided detailed mathematical proof of the existence of sub-optimal local minima for deep neural networks with smooth activation. 
The authors show how it is not possible to create general mathematical rules to guarantee convergence to good local minima.
Recent research shows that local minima can, in practice, be distinguished  {by} visualizing the loss function (\cite{sun2020global}) and, in particular, the occurrence of bad local minima can be empirically reduced with some architectural choices (\cite{li2018visualizing}).
}\dueref{We show the role of selecting different stationary points in \Cref{sec:opt_default}, where a multistart approach is also used on a subset of algorithms that seem particularly affected by the starting point.}

The issues surrounding hyper-parameters of optimization algorithms {are} also an important field for the \dueref{ML} research community.
\dueref{These hyper-parameters are often treated in the same way as the hyper-parameters defining the architecture (layers, neurons, activation functions, etc.), thus causing possible confusion about the reason for the good/bad performance of the obtained classification model.
More specifically, despite the problem of local minimizers being certainly well-known, this has been studied more in relation to the loss landscape,  namely in relation to architecture hyper-parameters, which can have an impact on shaping the loss landscape.}
However, recent research highlights that hyper-parameters setting can have a strong influence on algorithms' behavior, if they are specifically tuned on a given {task.} 
\cite{xu2020second} are amongst the firsts to point out the problem of the robustness of hyperparameters; they discuss how traditional first-order methods can get stuck in bad local minima or saddle points when tackling non-convex ML problems and how computational results can depend on the hyperparameters setting. 
\cite{jais2019adam} carry out a thorough analysis of Adam algorithm performance on a classification problem, focusing on optimizing the network structure as well as Adam parameters.
Nonetheless, Hyper-parameters are often set to a default value, which is obtained by maximizing the aggregated (in most cases, the average) performance across a variety of tasks, balancing a trade-off between efficiency and adaptability to different datasets (\cite{probst2019tunability,yang2020hyperparameter,bischl2023hyperparameter}).
\unoref{To our knowledge, no one has systematically addressed the question of whether it could be convenient to tune the hyper-parameters on a baseline problem (small network and small dataset) and use the tuned configuration on other problems (network-dataset) rather than using the default setting, which is our claim iii).}
\dueref{Indeed, performing a grid search is computationally expensive for the considered task due to the high amount of training time needed for each possible combination of hyperparameters.
Thus,
we aim to show that performing a single grid search for hyperparameters and tuning them for the baseline network on a simple  dataset can also have advantages on more complex problems (network dataset). The grid search on the baseline problems is reported in \Cref{sec:grid_search}. We then reuse the best-identified hyperparameter setting 
to investigate the effect when the architecture changes (in \Cref{sec:deep_wide}), and as the dataset varies, (\Cref{subsec:new_datasets}).
Indeed, we aim to analyze if the high-demanding operation of the grid search is more dataset-oriented or architecture-oriented, i.e., if the hyperparameters are more sensitive when the architecture or dataset change, given the same (classification) task.}

\section{The task and the Network Architectures}
\label{sec:NetArc}
The chosen task is multi-class {image} classification, which is a predictive {modeling} problem where a class, among the set of classes, is predicted for a given input data. 
More formally, we are given a training set made up of $P$ pairs $(x^j, y^j)$, $j=1,\dots,P$, of two-dimensional  input colourful images $x $ represented by \unoref{$H \times W$} pixels for each of the three colour channels (red, green, blue) thus as a tensor $(3 \times H \times W)$, and the corresponding class label $y$.
% \todo[inline]{Poiché i datasets hanno risoluzione diversa: UC Merced 256 mentre CONV 10 e 100 solo 32, è necessario specificare qui quanti sono i pixel ? \\
% \newline
% Lorenzo: I dataset sono stati forse catturati con sensori diversi o magari normalizzati dagli autori in determinate risoluzioni, è normale in computer vision - Volendo possiamo non specificarlo comunque, non è strattamente necessario.}
We denote with $N$ the number of possible classes, so that $y^j\in \{0,1\}^N$  and the target class value of sample $j$ is $y^ j_{i}= 1$ if the sample image $j$ belongs to class $i$ and $y^ j_{i}= 0$ otherwise.
% $$y^ j_{i}=\begin{cases}
% 1 & \text{if the sample image } j \text{ belongs to class } i \cr
% 0 &\text{otherwise}
% \end{cases}
% $$ 
\dueref{In this paper, we consider three well-known 2D input images datasets, described in Section \ref{sec:dataset}. }

Developed and formalized by (\cite{lecun1995convolutional}), {deep Convolutional Neural Networks (CNN)}, a special type of deep neural network (DNN) architecture, are one of the most widespread types of neural network for image processing, e.g. image recognition and classification \cite{hijazi2015using}, monocular depth estimation (\cite{9762984}), semantic segmentation (\cite{guo2018review}), video recognition (\cite{ding2017trunk}), and vision, speech, and image processing tasks (\cite{Abbaschian2021DeepLT, Kuutti2021ASO, Shorten2021DeepLA}).
In the literature, well-known Deep CNN models have been developed to face multi-class classification.
Among them, we cite the DenseNet, (\cite{Huang2017DenselyCC}), the ResNet, (\cite{He2016DeepRL}), and the MobileNet, (\cite{Howard2017MobileNetsEC}).
These architectures are distinguished by specific and complex designs composed of stacked operational blocks.  

In this paper,  different optimizers' are tested over the three datasets and different CNN architectures. 
In particular, we \dueref{specifically} design a lightweight low-complexity \textsc{Baseline} CNN model composed of elementary operations such as Convolution (Conv2D), Pooling, and Fully Connected (FC) layers, briefly described below.
{Moreover, starting from the \textsc{Baseline} model, we designed three architectural variants based on the same elementary blocks, varying the number of units per layer \dueref{(\textsc{Wide})}, the number of layers \dueref{(\textsc{Deep})},  and both of them  \dueref{(\textsc{Deep}\&\textsc{Wide})}.
\dueref{We refer to these architectures as \emph{Synthetic Networks.}}

\dueref{The Synthetic Networks}
%Those models, namely \textsc{Wide}, \textsc{Deep}, and \textsc{Deep}\&\textsc{Wide},
have been designed to analyze if the \unoref{tuned} set-up 
of the hyperparameters for the \textsc{Baseline} architecture \dueref{on a baseline dataset  shows similar improvements}  across  architectural  \dueref{and dataset changes}. 
\dueref{To check whether this behaviour can be generalized to other architectures,} %To the same aim, 
we also used two traditional CNN architectures, namely Resnet50}  (\cite{He2016DeepRL}) {and Mobilenetv2} (\cite{sandler2018mobilenetv2}).

In this section, we present the architectural aspects of the \textsc{Baseline} CNN architecture and its  \textsc{Wide}, \textsc{Deep}, and \textsc{Deep}\&\textsc{Wide} variants. Details and the mathematical formalization of the operations performed by the different layers are presented in Section \ref{sec:OptProb}.
%describe the \textsc{Baseline} CNN model and its  \textsc{Wide}, \textsc{Deep}, and \textsc{Deep}\&\textsc{Wide} variants that we have designed specifically for our experiments. 
%In this section, we highlight the architectural aspects only. Details and the mathematical formalization of the operations performed by the different layers are presented in Section \ref{sec:OptProb}.

\begin{figure*}[h!]
    \centering
    \subfigure[]{\includegraphics[scale=0.6]{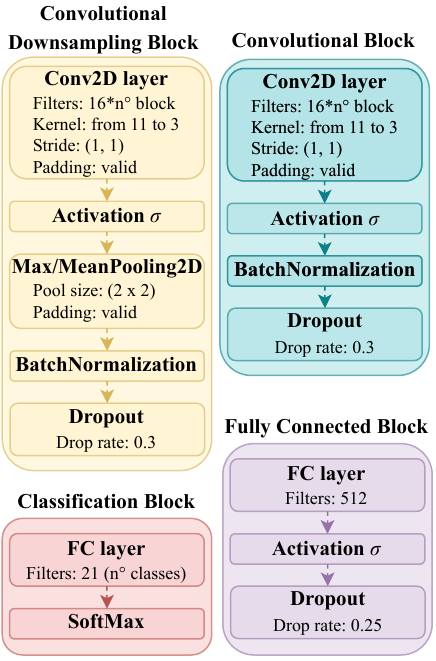}
    \label{img:layers}
    } \hfill
    \subfigure[]{\includegraphics[scale=0.33]{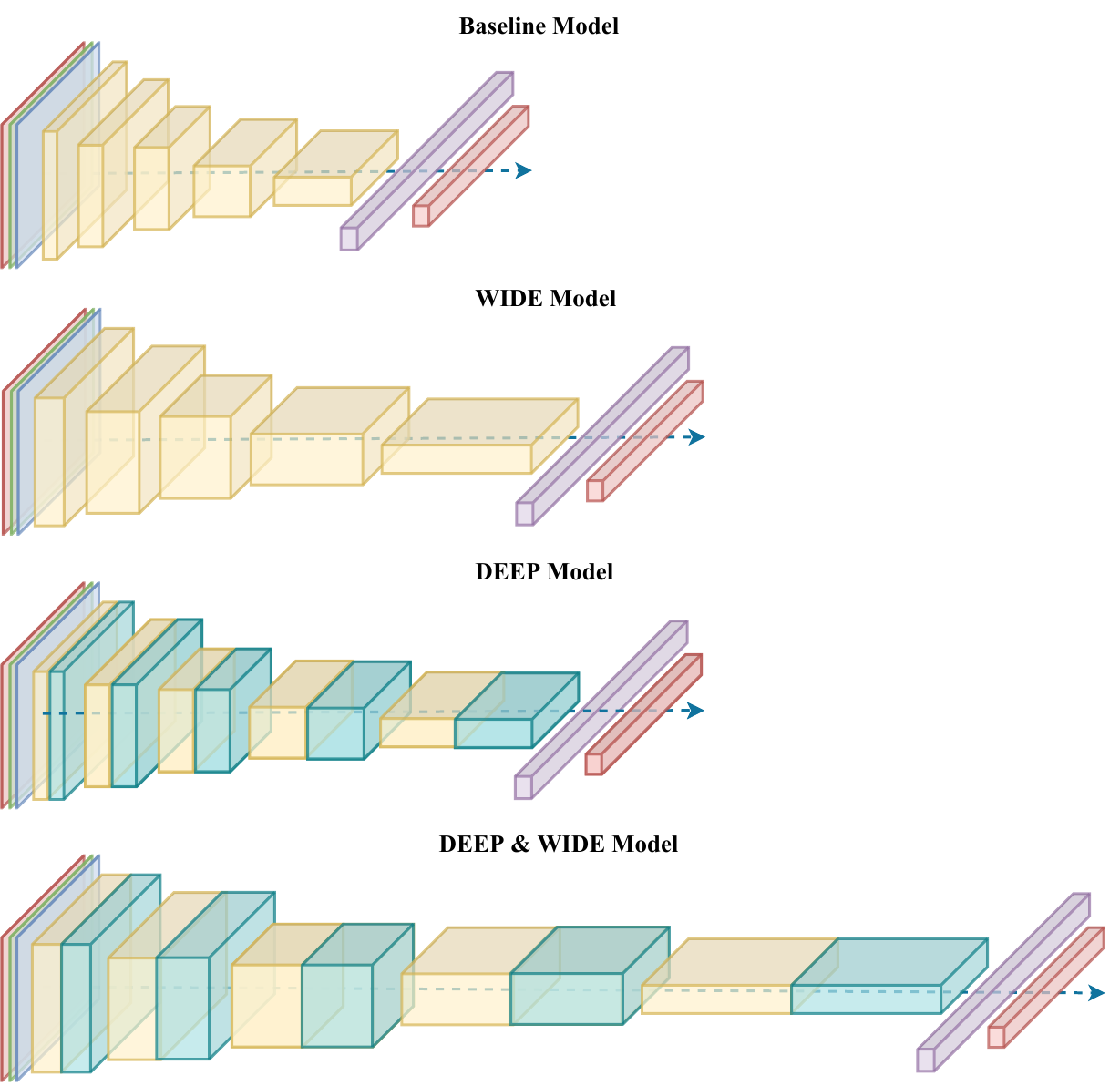}
    \label{img:models}
    }
    \caption{Overview of the \textsc{Baseline} model and corresponding blocks used for the study and its variants \textsc{Wide}, \textsc{Deep}, and \textsc{Deep}\&\textsc{Wide}. {A generic activation function $\sigma$ is used.}}
    \label{fig:structure_overview}
\end{figure*}

{ The \textsc{Baseline}
  CNN is composed  of a cascade of\begin{itemize}
      \item[-] five Convolutional Downsampling Blocks (CDBs) 
      \item[-] one Fully Connected Block (FCB) 
      \item[-] one final  Classification Block (CB).
  \end{itemize}}
{
A graphical representation of the models and a detailed block diagram representation, with layers operations and respective parameters, are reported in \Cref{fig:structure_overview}.
Each CDB block, {represented with the yellow blocks in Figure} \ref{fig:structure_overview},
performs a sequence of operations}

{
\begin{itemize}
    \setlength \itemsep{0.25em}
    \item[-] a  standard 2D-convolution (Conv2D layer) which takes in input a tensor of dimension $C \times H \times W$ 
and produces in output a new tensor with the spatial feature dimensions $H$ and $W$ decreased and the channels $C$ increased;
\item[-] application  of an activation function $\sigma$;
\item[-]   2D-max-pooling or 2D-mean-pooling  operation allows downsampling the extracted features along  {their} spatial dimensions by taking the maximum value or the mean over a fixed-dimension, %\change{input window}
{known as pool size};
\item[-] batch normalization.
\end{itemize}}

{Both CDB and FCB allow dropout with a given drop rate. Dropout consists of removing randomly parameters during optimization. Thus, it affects the structure of the objective function during the iterations by fixing some variables, and it can be seen as a sort of decomposition over the variables. }

{
The  {FC}  layer is a shallow Feed-forward Neural  Network (FFN) where all the possible layer-by-layer connections are established. The FCB block uses the Dropout operation, which is considered a trick to prevent overfitting. The last Classification layer is made up of a FC layer too, followed by the SoftMax activation function in order to extract the probability of each class.
}

 An overview of the input-output shapes of the \textsc{Baseline} model and the respective number of trainable parameters is reported in \Cref{tab:architecture_structure}.

%\remove{The first three operations are used for the feature extraction, whereas the other two operations have been applied to regularize and improve both the training process and the model generalization capability while avoiding overfitting.  Indeed, the standard 2D-convolution (Conv2D layer) takes in input a tensor of dimension $C \times H \times W$ and produces in output a new tensor with the spatial feature dimensions $H$ and $W$ decreased and the channels increased. The 2D-max-pooling or 2D-mean-pooling  operation allows one to downsample the extracted features along their spatial dimensions by taking the maximum value or the mean over a fixed-dimension input window. The proposed architecture, with its simple design and a limited number of trainable parameters, is ideal for the proposed analysis.}

%\remove{A graphical representation of the developed reference model and a detailed block diagram representation, with layers operations and respective parameters, are reported in , and}

\begin{table}[h]
    \footnotesize
    \centering
    \caption{\textsc{Baseline} architecture: input-output spatial dimensions, reported in the {[Channels, Height, Width]} format, and number of trainable parameters (N° Param.) of each layer \unoref{in the case of a $256 \times 256$} input.}
    \begin{tabular}{ l | c c }
        \multirow{2}*{Operations sequence} & Input Shape & Output Shape  \\
         & {[C,H,W]} & {[C,H,W]}  \\
        \hline
        Convolutional Downsamplig Block$_1$ & (3, 256, 256) & (16, 123, 123) \\
        Convolutional Downsamplig Block$_2$ & (16, 123, 123) & (32, 57, 57)  \\
        Convolutional Downsamplig Block$_3$ & (32, 57, 57) & (64, 25, 25 \\
        Convolutional Downsamplig Block$_4$ & (64, 25, 25) & (128, 10, 10))  \\
        Convolutional Downsamplig Block$_5$ & (128, 10, 10) & (256, 4, 4)  \\
%        Flatten & (4, 4, 256) & (1, 1, 4096) & 0 \\
   Fully Connected Block & (256, 4, 4) & (512, 1, 1) \\
        Classification Block & (512, 1, 1) & ($N$, 1, 1) \\
    \end{tabular}
    \label{tab:architecture_structure}
\end{table}
% \todo[inline]{La prima riga è 123 o 128 ? 
% Se i pixel non sono 256 ma 32, potrei scrivere tutto parametrico ocn il numero di pixel iniziali \#p ? oppre forse questa dovremo spostarla quando descriviamo i datasets.
% Mi semmbra che anche output finale del calssification block dipende dal numero di classi del dataset e in questo caso UC Merced ! QUIndi non è legato solo all'architettura ma anche al data set usato. SOlo che ancora il dataset NON lo abbiamo introdotto. 
% }
%\remove{Moreover, starting from the baseline model, we designed three architectural variants based on the same elementary blocks. Those models, namely \textsc{Wide}, \textsc{Deep}, and \textsc{Deep}\&\textsc{Wide}, are designed to study the behaviour of the optimized hyperparameters, i.e. understand if their optimal set-up for the baseline model can be generalized across its architectural changes.The quantitative analysis is then provided in .}
The %\change{applied variations}
{\textsc{Wide}, \textsc{Deep}, and \textsc{Deep}\&\textsc{Wide} architectures} are detailed below, while a block diagram for each designed architecture is shown in \Cref{fig:structure_overview}.

\begin{itemize}
    \setlength \itemsep{0.5em}
    \item The \textsc{Wide} model is designed by doubling the dimension of the output of each Conv2D layer, i.e. the number of output filters of the \textsc{Baseline} model in the convolution. 
    \item The \textsc{Deep} model is designed by doubling the number of convolutional operations, i.e. stacking to each CDB a further Convolutional Block (CB), as reported in \Cref{fig:structure_overview} (blue blocks), performing the same operations as the CDB except for the downsampling step of the 2D-max/mean pooling.
    \item The \textsc{Deep}\&\textsc{Wide} model is designed by combining the previous \textsc{Wide} and \textsc{Deep} structures.%\remove{, i.e. doubling the output dimensionality and the number of operations of the convolutional layers}.
\end{itemize}

{Finally, in order to assess the \dueref{generality} of the computational results, tests have been carried out also using two well-known neural architectures: Resnet50   (\cite{He2016DeepRL}) and Mobilenetv2   (\cite{sandler2018mobilenetv2}).
Resnet50 is a deep convolutional neural network 
%and has a similar structure  {to} the \textsc{Baseline} model
%\todo[inline]{lo toglierei perché dopo diciamo che è quella più diversa e infatti i risultati sono molto male su resnet. In effetti se c'è una connessine reisduale con output recedente anche il modello di ttimizzazine cambia perché gli output dei layer più esterni passano più volte nei layer successsivi}
with 50 hidden layers \unoref{and with residual connections at each layer, meaning that the output of each layer is added to the output of the subsequent layer in order to prevent the well-known vanishing gradient issue (\cite{borawar2023resnet}).}
Mobilenetv2 is a lightweight convolutional neural network, \dueref{which uses a ligher convolutional operator.}
%, designed to be more efficient and smaller than Resnet50.
Both Resnet50 and Mobilenetv2 are  {among} the most used neural architectures for image classification.}

\section{The optimization problem \dueref{for the Synthetic networks}}\label{sec:OptProb}
The optimization problem related to our task consists in the unconstrained minimization of the Categorical Cross Entropy ({CCE})
%as implemented in TensorFlow\footnote{ \url{https://www.tensorflow.org/api_docs/python/tf/keras/losses/CategoricalCrossentropy}} 
between the predicted  output of the neural model $\hat{y}^ j_{i}(\omega)$ and the correct classes $y^j\in \{0,1\}^N$.
%\remove{output of the neural model and the correct labels, i.e., between the predicted $\hat{y}^ j_{i}(\omega)$, which depends on the network weights $\omega$, and the correct classes $y^j\in \{0,1\}^N$.}

The predicted value $\hat{y}^ j_{i}(\omega)$ is the output of the last Classification block, and it represents the probability, estimated by the neural architecture, that the sample $j$ belongs to class $i$. 
Thus, we have $\sum_{i=1}^{N} 
\hat{y}^ j_{i}(\omega) = 1 \ \forall j = 1,\dots,P \ \forall \omega \in \mathbb{R}^h$. Thus, the unconstrained optimization problem can be written as: 

\begin{equation}
    \min_{\omega \in \mathbb{R}^h} f(\omega) := - \frac 1 P\sum^P_{j=1} \sum^{N}_{i=1} y^j_{i} \log(\hat{y}^j_{i}(\omega))
    \label{eq:loss}
\end{equation}

{We aim to derive the probability output of the \dueref{\textsc{Baseline} synthetic network} %neural model 
as a function of the network parameters $\omega$, i.e., $\hat{y}^ j_{i}(\omega)$, and, in particular, to write the dependency of $\hat{y}^ j_{i}$  from $\omega$ in a closed form}

 As reported in \Cref{fig:structure_overview}, in the \textsc{Baseline} model the input images propagate along 
convolutional downsampling layers (CDB), a fully connected (FCB) and a classification (CB) layers. We formalize this process to get an analytical expression for $\hat{y}^ j_{i}(\omega)$.

{Each \unoref{CDB} performs five operations: a standard 2D-convolution, denoted as Conv2D layer in \Cref{fig:structure_overview}, followed by an activation function, a Max-pooling or a Mean-pooling operation, the batch normalization.
In the FCB, the Max/Mean-pooling is removed.}

%\remove{A 2D-convolutional layer $\ell$, is made up of $d^{\ell-1} \times d^{\ell}$ filters (also called kernels), where $d^{\ell-1}$ is the number of the input channels and $d^{\ell}$ the number of the output channels. The Conv2D layer applies the discrete convolution between the $k$-th filter of the $c$-th input channel, denoted as $w^k_c \in \mathbb{R}^{n \times n}$, and its input $X^{\ell-1}$, $\ell=0,1,\dots,L$ and then, for each output channel, performs the sum of the convoluted output along the input channels. More formally, the input $X^{\ell-1}$ to a convolutional layer $\ell$ is a tensor of dimension $[d^{\ell-1}, m^{\ell-1}, m^{\ell-1}]$ where $m^\ell$ is the height/width for each channel $c=1,\dots, d^{\ell-1}$ and produces  an output  tensor $X^\ell$ of dimension $[d^{\ell},m^{\ell}, m^{\ell}]$. The input $X^0$ at layer $\ell=0$  is the colorful sample image  represented with $m^0=256$ pixels and $d^0=3$ color channels (red, green, blue). We denote by $x_c\in \mathbb{R}^{ m \times m}$ the input $X^0$ of each channel $c=1,\dots, d^0$ and, more in general, by  $X^\ell_c$ the matrix $m^{\ell} \times m^{\ell}$ for each channel $c$.}

{The input $X^{\ell-1}$ to a Conv2D layer $\ell$ is a tensor of dimension $[d^{\ell-1}, m^{\ell-1}, m^{\ell-1}]$ where $d^{\ell-1}$ are the channels and $m^\ell$ is the height/width of each channel $c=1,\dots, d^{\ell-1}$, and the output $X^\ell$ is a tensor of dimension $[d^{\ell}, m^{\ell}, m^{\ell}]$.
The input $X^0$ at layer $\ell=0$ of the CDB layer is the colourful sample image  represented with $m^0=256$ pixels and
$d^0=3$ colour channels (red, green, blue). 
 We denote by $X^\ell_c$ the matrix $m^{\ell} \times m^{\ell}$ for each channel $c=1,\dots, d^{\ell}$ and by $x_c=X^0_c\in \mathbb{R}^{ m \times m}$ for $c=1,\dots, d^0$.
A Conv2D layer $\ell$ applies a discrete convolution on the input $X^{\ell-1}$.
This operation consists in applying filters (also called kernels) $w^k_c \in \mathbb{R}^{n \times n}$ for $k=1,\dots, d^{\ell}$ to the $c$-th input channel $X^{\ell-1}_c$ with $c=1,\dots, d^{\ell-1}$. 
 }

The convolution operation depends on the integer stride $s\ge 1$, representing the amount by which the filter $w^k_c$ shifts around the input $X^{\ell-1}_c$.
\unoref{The stride is commonly set to $s=1$, as we did for all the experiments except for L-BFGS for which the stride has been fixed to $s=2$.}
The dimension of the filter $n$, the number of channels $d^\ell$,  {and} the stride $s\ge 1$, for each layer $\ell$, are network hyperparameters. 
Let us denote as  $\otimes_s$ the discrete convolution operation with stride $s\in \mathbb{N}$.
 The expression componentwise of the convolutional operation between the filter $w^k_c \in  \mathbb{R}^{n \times n}$  and the input feature $X^{\ell-1}_c$ is 
 
%\begin{equation}
%    \label{eq:convolution}
$$
\left[w^k_c \otimes_s X_c\right]_{ij} = \sum_{a=1}^n \sum_{b=1}^n [w^k_c]_{a,b}[X_c]_{i+sa, j+sb} \quad i, j=0,\dots, m-n
$$
%\end{equation}
 where for the sake of simplicity, we avoid the use of the \unoref{superscript} $\ell$.
As reported in \cite{bengio2017deep} Chapter 9 (equation (9.4) with $s=1$), 
the $k$-th convoluted output is the matrix $F_k^\ell \in  \mathbb{R}^{\frac{(m^{\ell-1}-n+1)}{s} \times \frac{(m^{\ell-1}-n+1)}{s}}$ defined as the sum over the channels, namely

\begin{equation}
    \label{feature_map_equation}
    F^\ell_k = \sum_{c=1}^{d^{\ell -1}} w^k_c \otimes_s X^{\ell -1}_c   \qquad k=1,\dots,  {d^{\ell}}
\end{equation}

The convoluted output $F^\ell$ of the Conv2D layer $\ell$ is thus
$$F^\ell=[F^\ell_k]_{k=1,\dots, d^\ell}= (F^1,\dots,F^{d^{\ell}}) \in \mathbb{R}^{d^{\ell} \times\frac{(m^{\ell-1}-n+1)}{s}\times \frac{(m^{\ell-1}-n+1)}{s}},$$

The parameters of the convolutional layer $\ell$ are denoted as 
\begin{equation}\label{eq:pesiconv}
W_{Conv}^{\ell} = \left[w^k_c\right]^{k=1,\dots,d^{\ell}}_{c=1,\dots,d^{\ell - 1}} \in \mathbb{R}^{n \times n \times d^{\ell} \times d^{\ell - 1}}    
\end{equation}

The next block applies a nonlinear activation function  $\sigma$ to the output  of the Conv2D layer $\ell$. The activation function $\sigma$ that has been used in the computational experiments in Section \ref{sec:Res} can be either the  ReLU 
 {$\sigma(z) = \max \{ 0,z \}$}
or the SiLU {$\sigma(z) = \frac{z}{1 + e^{-z}}$} (\cite{mercioni2020p,ramachandran2017searching}). In particular, the SiLU is used when applying  L-BFGS to ensure  {the} smoothness of the objective function \dueref{and avoid failures of the optimization procedure.}

{In  CDB blocks, the Conv2D output is submitted to the pooling operation aimed at further reducing the image dimensionality. 
The pooling operation \texttt{Pool}  involves sliding a two-dimensional 
 non-overlapping $p \times p$ matrix, where $p$ is the pool size,
over each convoluted output $F^\ell_k$ and contracting the features lying within the $p\times p$ region covered by the filter by using the \texttt{max} or the  \texttt{mean} operation.
More formally, let us introduce the set ${\cal P}_{(i,j)}$  of positions i.e., for each $(i,j)$ the rows and columns   of the matrix $F^\ell_k$, that are located in the $p\times p$ region; then \texttt{Pool}: $F^\ell_k\to G^\ell_k$ where the  $k$-th pooled output is $
G^\ell_k \in \mathbb{R}^{\frac{(m^{\ell -1}-n+1)}{sp} \times \frac{(m^{\ell -1}-n+1)}{sp}}$ .
  The output of the Max-Pool or Mean-Pool  layer is computed respectively as: 
\begin{equation}
    (G_k)^{Max}_{ij}  = \max_{(r,c) \in {\cal P}_{(i,j)}} (F_k)_{r,c}
\qquad
\qquad
    (G_k)^{Mean}_{ij} = \frac 1 {\left|{\cal P}_{(i,j)}\right|}\sum_{(r,c) \in {\cal P}_{(i,j)}} (F_k)_{r,c}
\label{eq:CNN}
\end{equation}
where for the sake of simplicity we removed the superscript $\ell$.}
The set ${\cal P}_{(i,j)}$ depends on how drastically we want to reduce the dimensionality.
In our experiments, we have set $p=2$. In this case, the moving region is just a $2 \times 2$  matrix, thus we halve the dimension  of $F$.
%\todo[inline]{ho tolto comment su $p=!$ per non chiamarci ulterior domande}
%\unoref{The only exception is L-BFGS, where, to preserve differentiability, we also performed experiments with $p=1$ and stride $s=2$.}
For instance, ${\cal P}_{(1,1)} = \{(1,1),(1,2),(2,1),(2,2)\}$, meaning that {$G_{11}^k$ is the maximum/mean value between four different values $\{F^k_{1,1},F^k_{1,2},F^k_{2,1},F^k_{2,2}\}$.} 
The Max-pooling operation is widely used in image classification, but it introduces a non-differentiability issue. For this reason, when testing  L-BFGS, where  differentiability is crucial, we use the Mean-pooling. 

The output of a CDB block obtained by 
\Cref{feature_map_equation} and \Cref{eq:CNN} is { $Z^\ell=\left\{\texttt{Pool} \left[ \sigma\left(Z^\ell_k\right)\right]\right\}_{k = 1,\dots, d_{\ell}}$  and finally given as
\begin{equation}
    \label{eq:convLAYER}
  Z^\ell= \left\{\texttt{Pool} \left[ \sigma \left(F^\ell_k\right) \right]\right\}_{k = 1,\dots, d_{\ell}}
\end{equation}
whereas the output of a CB block does not use the Pool operations and thus is given simply by 
$ Z^\ell= F^\ell$.
In both cases,
 $Z^\ell $ is then normalized to stabilize and speed up the training process. 
The normalization is performed following the standard batch-normalization procedure described in \cite{ioffe2015batch}, i.e., subtracting the mean and dividing by the standard deviation.
% \begin{equation}
%     \label{eq:batch_norm}
%     Z^\ell_{norm} = \frac{Z^{\ell}-\mathbb{E}\left[Z^{\ell}\right]}{\sqrt{\texttt{Var}\left[Z^{\ell}\right] + \epsilon}}
% \end{equation}
% where $\epsilon>0$. $Z^\ell_{norm}$ is the output $X^\ell$ of the $\ell$ CDB layer.}

The output of the last layer $X^L$, \unoref{being $L$ the total number of layers}, of either CDB or CB is finally sent into the Fully Connected Block (FCB) and then into the Classification Block.
The FCB is a shallow  Feed-forward Neural (FFN) Network with $M_{FCB}$ neural units and the activation function $\sigma$ (ReLU or SiLU, as  before).
The output of the FCB is given by:
\begin{equation}
    Z_{FCB} = \sigma (W_{FCB}^T X^{L} + b_{FCB} )
    \label{eq:FC}
\end{equation}
where $(W_{FCB},b_{FCB})$ are the weights and the biases of the FFN network.
The Classification Block is made up of a  shallow FFN network followed by a softmax operation, so that the final output is
\begin{equation}
  \hat y= \texttt{Soft}(  Z_{CB}) =\texttt{Soft}\left ( W_{CB}^T Z_{FCB} + b_{CB} \right)
    \label{eq:CB}
\end{equation}
where 
$\texttt{Soft}$ is applied component-wise to the vector $Z_{CB}$ as
$$\texttt{Soft}(z_h) = \frac{e^{z_h}}{\sum_{j=1}^{N} e^{z_j}}.$$

The overall network parameters are
 $$\omega=(W_{CB},b_{CB},W_{FCB},b_{FCB},W_{Conv}^L, \dots, W_{Conv}^1).$$

\dueref{ Mobilenetv2   (\cite{sandler2018mobilenetv2}) and the Resnet50   (\cite{He2016DeepRL}) presnets  differences with respect to the Synthetic Network.
Indeed,
 Mobilnetv2 instead uses depthwise separable convolutions different from \eqref{feature_map_equation}, while 
Resnet50 presents residual connections among layers. Thus, the resulting 
 optimization problem can be different with respect to the one described in this section. }
 
\section{The selected optimization algorithms}
\label{sec:OptAlg}
In this paper, L-BFGS and eight state-of-the-art ML optimization algorithms with multiple hyperparameters setups are compared. Precisely, those are: Adam (\cite{Kingma2015AdamAM}), Adamax (\cite{Kingma2015AdamAM}), Nadam (\cite{Dozat2016IncorporatingNM}), RMSprop\footnote{RMSprop is an adaptive learning rate method devised by Geoff Hinton in one of his Coursera Class (\url{http://www.cs.toronto.edu/~tijmen/csc321/slides/lecture_slides_lec6.pdf}) that is still unpublished}, SGD (\cite{robbins1951stochastic, bottou2018optimization, ruder2016overview, Sutskever2013OnTI}), FTRL (\cite{mcmahan2011follow}), Adagrad (\cite{duchi2011adaptive}), and Adadelta (\cite{Zeiler2012ADADELTAAA}). 
We use the  SciPy\footnote{\url{https://scipy.org/}} version  for  L-BFGS and
 the built-in implementation in TensorFlow library\footnote{\url{https://www.tensorflow.org/probability/api_docs/python/tfp/optimizer/lbfgs_minimize}} for the eight others.

For the sake of completeness, we report the updating rule of each algorithm, assuming it is applied to the problem as in \Cref{eq:loss}, namely 
  $$\min_{\omega\in \mathbf{R}^h} f(\omega)=  \sum_{p=1}^P f_p(\omega).$$
 
 We note that all algorithms require $f$ to be a continuously differentiable function. The use of non-differentiable activation functions (ReLU) in the network layers and the MaxPooling layer, as usually done in CNN, implies that the objective function $f$ does not satisfy this essential property and a finite number of non-differentiable points arise.
When using L-BFGS, this aspect becomes evident as discussed in the  \Cref{sec:opt_default}. {\unoref{We} tried to use L-BFGS with the standard setting in CNNs, but it happened very often that the method \dueref{failed and }ended at a non-stationary point. Indeed, since L-BFGS is a full-batch method using all the samples at each iteration, whenever a point of non-differentiability is reached, the gradient returned by TensorFlow is \texttt{None}, and the method gets stuck. 
The other eight first-order algorithms are instead mini-batch methods, which perform network parameters update using only a small subset of the whole samples. When it happens that the partial gradient is \texttt{None} on a subset of samples, the method continues in the epoch, changing the batch and possibly the new partial gradient can be used to move from the current iteration.
Hence, although convergence of the mini-batch methods requires smoothness, from the computational \dueref{point of view} they can work heuristically without it.
}

% \todo[inline]{Gusto per essere sicura che quanto scritto di seguito sia corretto: SilU e Mean-Pool lo abbiamo usato per TUTTI gli algoritmi ? O solo per L-BFGS ? NOnn si capisce bene \\
% \newline
% Lorenzo: Li abbiamo usati solo per LBFGS, per tutti gli altri, e per la grid search abbiamo usato le configurazioni con ReLU e MaxPooling. Questo perchè tutti gli studi riportati sulla rete semplice, deep, wide, e deepwide sono stati fatti prima di aggiungere LBFGS, e la sua variante con i componenti differenziabili (che prima non conoscevo).}
\dueref{Hence when using L-BFGS,} we need to reduce non-differentiability. To this aim we set SiLU as activation function and we select the \dueref{MeanPooling}, 
which is not a common practice in CNN for image classification. 
We also fully deactivate the Dropout operation \unoref{and set stride s = 2}.
{We also remark that we are  comparing a globally convergent traditional full-batch method with eight different mini-batch methods, that require strong assumptions to prove convergence that do not hold for the problem at hand.
By studying L-BFGS performance against commonly used optimizers we aim  {to assess} whether theoretical convergence really plays an important role in determining the efficiency, the train performances, and, most of all, the generalization capability.}

\subsection{L-BFGS}
Being one of the best-known first-order methods with strong convergence properties (see \cite{liu1989limited}), L-BFGS belongs to the limited memory quasi-Newton methods class. This algorithm is purely deterministic and, at every iteration $k$, exploits an approximated inverse Hessian of the objective function,
and it
performs the following update scheme:
$$ \omega_{k+1} = \omega_k - \eta_k H_k \nabla f(\omega_k)
$$
where $\eta_k$ is a step size obtained via some line search method.

The updating rule for $H_k$ has been formalized by  \cite{nocedal1999numerical}.
Given an initial approximate Hessian $H_0 \sim \left[\nabla^2 f(\omega_0)\right]^{-1}$, the algorithm uses the rule:
$$ H_{k+1} = V^T_k H_k V_k + \frac{s_k s^T_k}{y_k^T s_k} 
$$
where
$$s_k = \omega_{k+1} - \omega_k, \qquad y_k = \nabla f(\omega_{k+1}) - \nabla f(\omega_{k}),\quad V_k = I - \frac{y_k s^T_k}{y_k^T s_k} ,$$ being $I$ is the identity matrix.
% And performs the following update scheme:
% $$ \omega_{k+1} = \omega_k - \eta_k H_k \nabla f(\omega_k)
% $$
% where $\eta_k$ is either the learning rate or, more generally, the step size obtained via some line search method.

L-BFGS is not among the optimizers mostly used in  {machine learning}. However, in force of its strong convergence properties, this algorithm  {has recently gained} increasing interest in the research community. Some multi-batch versions of L-BFGS have been proposed in the past years in (\cite{berahas2016multi, bollapragada2018progressive, berahas2020robust}), in particular for image processing tasks in medicine in (\cite{yun2018small,wang2019accelerating}).

Since L-BFGS is not directly available in TensorFlow, we have used the SciPy version implemented with an open-source wrapper available online\footnote{\href{https://gist.github.com/piyueh/tf\_keras\_tfp\_ L-BFGS.py}{https://gist.github.com/piyueh/712ec7d4540489aad2dcfb80f9a54993}}.

\subsection{SGD}
The Stochastic Gradient  {Descent} (SGD) is the basic algorithm to perform the minimization of the objective function using a direction $d^k$ which is {random} estimate of its gradient (see for details, the comprehensive survey \cite{bottou2018optimization}).
In the TensorFlow implementation, the following mini-batch approximation is used:

$$ g_k (\omega_k)=\displaystyle\frac 1 {|{\cal B}_k|} \sum_{i\in {\cal B}_k}\nabla  f_{i}(\omega_k)\text{ for some }{\cal B}_k\subset\{1,\dots,P\}
$$

{
The updating rule  is given by
$$ \omega_{k+1} = \omega_k - \eta_k g_k(\omega_k) 
$$
The update step can be modified by adding a momentum term \dueref{(which depends on a parameter $\beta$)} or a Nesterov acceleration step  which are an extrapolation steps along the difference between the two past iterations.
TensorFlow  allows the use of a boolean parameter, called \texttt{Nesterov}, which enables the Nesterov acceleration step (see \cite{sutskever2013importance}). 
When \texttt{Nesterov}=False, only a momentum is applied and the basic SGD iteration is modified by  adding
$$\beta_k (\omega_k - \omega_{k-1})$$  with $\beta_k \ge 0$ momentum parameter.
 When \texttt{Nesterov}=True, first 
\begin{equation}\label{nesterov}
z_k = \omega_k + \beta_k (\omega_k - \omega_{k-1})\end{equation}
is computed and the updating rule becomes $$\omega_{k+1} = z_k - \eta_k  g_k(z^k).$$
} 
 In both cases, the value of $\beta$ is a hyperparameter to be tuned.

{In \cite{bertsekas2000gradient} the pure SGD has been proved to converge to a stationary point under strong assumptions. We further refer the reader to \cite{drori2020complexity} for a thorough analysis of SGD complexity and convergence rate.}

\subsection{Adam}
Adam is one of the first SGD extensions, where the gradient estimate is enhanced with the use of an exponential moving average according to two coefficients: $\beta_1$ and $\beta_2$, \unoref{ranging in $(0,1)$}.
The index $i \in \{1,2\}$ is referred to as the moment of the stochastic gradient, i.e., the first moment (expected value) and the second moment (non-centred variance). 
Being $g_k$ the same mini-batch approximation used in SGD,% i.e., $g_k=\displaystyle\frac 1 {|{\cal B}_k|} \sum_{i\in {\cal B}_k}\nabla  f_{i}(\omega_k)$, 
we define the following first and second-moment estimators at iterate $k$:
\begin{equation}\label{eq:mk}
 \mathbb{E}\left[ \nabla f(\omega_k) \right]  \sim  m_k = (1-\beta_1) \sum_{i=1}^k \beta_1^{k-i} g_i
\end{equation}
\begin{equation}\label{eq:vk}
 \mathbb{E}^2\left[ \nabla f(\omega_k) \right] \sim v_k = (1-\beta_2) \sum_{i=1}^k \beta_2^{k-i} {(g_i\odot g_i)}
\end{equation}
{where $\odot$ is the Hadamard component-wise product among vectors.}

{Given the following matrix:
$$ \widetilde{V}_k(\varepsilon) =\frac{1}{1-\beta_2^k} \left[I \varepsilon + \mbox{diag}({v}_k)\right]^{\frac 1 2}$$
where $\mbox{diag}(v)$ denotes the diagonal $h \times h$ matrix with elements
$v_i$ on the diagonal, and $\epsilon >0$, the updating rule is the following:
$$ \omega_{k+1} = \omega_{k} - \eta_k \widetilde{V}_k(\varepsilon)^{-1} \frac{1}{1-\beta_1^k} m_k
$$
where $m_k$ is given in \eqref{eq:mk}.}
{It has been recently proved in \cite{defossez2020simple} that Adam can converge under smoothness assumption and gradients boundness in $L_{\infty}$ norm with convergence rate $O(\frac{h\log(N)}{N})$, being $h$ the number of variables and $N$ the numbers of iterations. 
For a more detailed discussion of Adam complexity (as well as for the other adaptive gradient methods Adamax and Nadam){, we refer} the reader to \cite{zhou2018convergence}.
}

\subsection{Adamax}
Adamax performs mainly the same operations described in Adam, but it does not make use of the parameter $\epsilon$, and the algorithm exploits the infinite norm to average the gradient. 
{Let $$u_k = \max_{i=1,\dots,k} \beta_2^{k-i}|g_i|\qquad U_k=\mbox{diag}(u_k)$$ and $m_k$ as in \eqref{eq:mk}. Thus, the updating rule is:
$$ \omega_{k+1} = \omega_{k} - \eta_k U_k^{-1}\frac{1}{(1 - \beta_1^k)}m_k.
$$
}

\subsection{Nadam}
Nadam, also known as Nesterov-Adam, performs the same updating rule as Adam but employs the Nesterov acceleration step \Cref{nesterov}. Nadam is expected to be more efficient, but the Nesterov trick involves only the order  {in} which operations are carried out and not the updating formula.

\subsection{Adagrad}
Adagrad is the first Adam extension that makes use of adaptive learning rates to discriminate more informative and rare features. The general update rule of $\omega_k \in \mathbf{R}^h$ involves complex matrix operations, for which we need to introduce some other notation.
{At iteration $k$ we introduce the cumulative vector $$G_k = \displaystyle\sum_{t=0}^{k-1} \left( g_t \odot g_t \right)$$ where
 $g_t=\sum_{i\in {\cal B}_t} \nabla f_i(\omega_t)  $. 
Given $\varepsilon > 0$ and the identity matrix $I$, we define the following matrix:
$$ H_k(\varepsilon) = \left[ I \varepsilon + \mbox{diag}(G_k) \right]^{\frac 1 2}
$$
where $\mbox{diag}(v)$, where $v\in  \mathbf{R}^h$, denotes the diagonal ${h\times h}$ matrix with elements $v_i$ on the diagonal. }
Thus, the updating rule resulting after the minimization of a specific proximal function (see \cite{duchi2011adaptive}) is the following:
$$\omega_{k+1} = \omega_k - \eta  H_k (\varepsilon)^{-1} g_k $$

{Convergence results concerning Adagrad, as well as its modification Adadelta, have been reported in (\cite{li2019convergence,defossez2020simple,chen2018convergence}) and they still require strong assumptions.}

\subsection{RMSProp}
Proposed by Hinton et al. in the unpublished lecture \cite{hinton2012neural}, RMSProp (Root Mean Square Propagation) performs a similar operations as Adagrad, but the update rule is modified to slow down the learning rate decrease. 

\unoref{In particular, f}ollowing the notation introduced in the last subsections,  
at iteration $k$
{the following matrix \unoref{is used}:
$$ {V}_k(\varepsilon) = \left[I \varepsilon + \mbox{diag}({v}_k)\right]^{\frac 1 2}
$$
where   $v_k$ be given by \eqref{eq:vk}.
The update rule is:
$$\omega_{k+1} = \omega_k - \eta {V_k (\varepsilon)^{-1}}    g_k  $$}

{Despite being still unpublished, RMSProp has already been studied in the field of convergence theory in (\cite{defossez2020simple,de2018convergence}), 
with an eye to the non-convex ML context.}

\subsection{Adadelta}
Adadelta (\cite{Zeiler2012ADADELTAAA}) can be considered as an extension of Adagrad, which allows for a less rapid decrease in learning rate. 
{Let us consider the same matrix $${V}_k(\varepsilon) = \left[I \varepsilon + \mbox{diag}({v}_k)\right]^{\frac 1 2}
$$ used in RMSProp.
Further, let $\delta_{k}=\omega_{k+1}-\omega_{k}$
$${\widetilde\delta}_{k} = \displaystyle\sum_{t=1}^{k} \rho^{k-t} (1 - \rho) \left( \delta_t \odot \delta_t \right)\qquad  \widetilde{\Delta}_k(\varepsilon) = \left[I \varepsilon + \mbox{diag}({\widetilde\delta}_k)\right]^{\frac 1 2}
$$}
 {Thus, we can write the Adadelta updating rule as follows:}

{$$ \omega_{k+1} = \omega_k - V_k(\varepsilon)^{-\frac 1 2 }
\widetilde{\Delta}_{k-1}(\varepsilon) 
g_k
$$}

\subsection{FTRL}
FTRL (Follow The Regularized Leader), as implemented in TensorFlow following \cite{McMahan2013AdCP}, is a regularized version of SGD, which uses the L1 norm to perform the update of the variables. Given, at every iteration $k$, $d_k = \sum_{t=1}^k g_t$, and fixed the quantity $\sigma_k$ such that $\sum_{t=0}^k \sigma_{t} = \frac{1}{\eta_k}$, the update rule is the following:
$$ \omega_{k+1} = \arg \min_{\omega} \{d_k^T \omega + \sum_{t=1}^k \sigma_t ||\omega - \omega_t||^2 + \lambda ||\omega||_1\}
$$

{As proved in \cite{McMahan2013AdCP}, the minimization problem in the update rule can be solved in closed form, setting:  
\begin{equation}
\label{eq:FTRLupdate}
    \omega_{k+1,i} = \begin{cases}
    0 \ \ \ \text{if } \vert z_{k,i}\vert \le \lambda \\
    -\eta_k (z_{k,i} - \texttt{sgn} (z_{k,i}) \lambda) \text{   otherwise}
    \end{cases}
\end{equation}
where $\displaystyle z_{k,i} = d_k - \sum_{s=1}^k \sigma_s \omega_s$ and $\texttt{sgn}(\cdot)$ is the Signum function.
}

{FTLR convergence can be proved only in the convex case, as explained  {in detail} in \cite{mcmahan2011follow}.}

\section{The datasets}
\label{sec:dataset}
\dueref{We have  carried out our computational test using three datasets:  UC Merced (\cite{uc_merced}),
 CIFAR10 and CIFAR100 (\cite{cifar10}). UC Merced represents the benchmark dataset used to define the \textsc{Baseline} problem, namely the training of the  \textsc{Baseline} network defined in Section \ref{sec:NetArc}. We use the \textsc{Baseline} problem (defined as the pair \textsc{Baseline} network - UC Merced) to assess
 the performance of the different optimization methods.}
%benchmark dataset used to show the performance of the different optimization methods \dueref{for training the } convolutional architectures 
%is the  (\cite{uc_merced}). 

 UC Merced is a balanced dataset that comprises a total of 2100 land samples divided into 21 classes, i.e. 100 images per class. The dataset images have a resolution of $256\times256$ pixels. The high number of classes and the limited number of samples for each class make the multi-class classification a non-trivial task.
% \dueref{We report in Table 
% \ref{tab:architecture_structure} the number of parameters of the \textsc{Baseline} on the UC Merced.
% \begin{table}[h]
%     \footnotesize
%     \centering
%     \caption{\textsc{Baseline} problem: input-output spatial dimensions, reported in the {[Channels, Height, Width]} format, and number of trainable parameters (N° Param.) of each layer.}
%     \begin{tabular}{ l | c c r | c c r}
%     & \multicolumn{3}{c}{UC Merced}\\
%         \multirow{2}*{Operations sequence} & Input Shape & Output Shape & N° Param. \\
%          & {[C,H,W]} & {[C,H,W]}  & [K] \\
%         \hline
%         Convolutional Downsamplig Block$_1$ & (3, 256, 256) & (16, 123, 123) & 5.9  \\
%         Convolutional Downsamplig Block$_2$ & (16, 123, 123) & (32, 57, 57) & 42.6 \\
%         Convolutional Downsamplig Block$_3$ & (32, 57, 57) & (64, 25, 25 & 100.7 \\
%         Convolutional Downsamplig Block$_4$ & (64, 25, 25) & (128, 10, 10)) & 205.4 \\
%         Convolutional Downsamplig Block$_5$ & (128, 10, 10) & (256, 4, 4) & 296.2 \\
% %        Flatten & (4, 4, 256) & (1, 1, 4096) & 0 \\
%    Fully Connected Block & (256, 4, 4) & (512, 1, 1) & 2097.6\\
%         Classification Block & (512, 1, 1) & (21, 1, 1) & 10.8  \\
%     \end{tabular}
%     \label{tab:architecture_structure}
% \end{table}
% }

In order to assess whether the computational results obatined on the \textsc{Baseline} problem % (i.e., the \textsc{Baseline} architecture on UC Merced) 
 generalize to different datasets, we have also carried out additional tests on two larger datasets: CIFAR10 and CIFAR100 (\cite{cifar10}), respectively, with $N=$10 and $N=$100 classes, both containing 60000 samples at a resolution of $H\times W= 32 \times 32$ pixels.

%\dueref{As regards CIFAR10 and CIFAR 100 the number of parameters of the overall network is of the order of $10^7$, from few millions up to approximately $40M$ for Resnet50.}
% \todo[inline]{CHECK: ho modification testo. verifica sia corretto o modifica in modo opportune. Ho inserito un ulteriore referenza. Forse non serve ? }
\unoref{Furthermore, for minibatch methods, we also apply data augmentation, which is a commonly used technique in   {machine learning} for image classification. It consists of random transformation of the selected mini-batch of samples with the aim of increasing the training dataset diversity and achieving better generalization capabilities. 
For a better understanding of this technique, we  {refer} the reader to \citep{van2001art,connor2019survey}; in our case, data augmentation involves random transformations on selected images, such as rotation, scaling, adding noise, and changing brightness and contrast.}

\section{Implementation details} 
\label{sec:impdet}
We implemented the proposed study using TensorFlow 2\footnote{\href{https://www.tensorflow.org/}{https://www.tensorflow.org/}} deep learning high-level API, \dueref{using its implementation of Categorical Cross Entropy (CCE)\footnote{ \url{https://www.tensorflow.org/api_docs/python/tf/keras/losses/CategoricalCrossentropy}}.} 
We set the environment seed (also for the normal initializer of the convolutional kernels) at a randomly chosen value equal to 1699806 or to a specific list\footnote{$[0, 100, 500, 1000, 1500, 10000, 15000, 100000, 150000, 1000000, 1500000, 1699806]$} of values in the multistart analysis. 
Computational tests have been conducted using a mini-batch size $bs = 128$, except for L-BFGS, which is a batch method, i.e., requires the whole gradient at each iteration. Concerning the eight built-in optimizers, in \Cref{sec:opt_default} the network was trained setting to 100 the number of epochs over the whole dataset. We remark that a single epoch consists of $\frac{P}{bs}$ update steps, being $P$ the number of samples in the dataset. In the experiments with tuned hyperparameters in Section \ref{sec:grid_search}, we halved the number of epochs. \unoref{In other experiments with larger problems, the number of epochs was further reduced to 30.}
\unoref{We underline that reducing the number of epochs is a common heuristic procedure in deep learning (\cite{diaz2017effective,yu2020hyper}), where at an early testing phase the number of epochs is set to an arbitrary value (in our case 100) and, then, it is reduced according to the training loss decrease, such that the network is not trained when the loss has already reached values close to zero and is not further improving.
This prevents any waste of computational time that could result from training the network when the loss is already extremely close to zero.}

Eventually, we underline that the TensorFlow implementation of the eight built-in optimizers, as well as the SciPy version of L-BFGS, uses {back-propagation algorithm to compute the gradients.}

The  training have been run on 12GB NVIDIA GTX TITAN V \dueref{GPU}. The L-BFGS algorithm, \unoref{being a full-batch method}, \unoref{cannot be run on a GPU due to the lack of memory storage,} and takes almost 30 seconds on our reference Intel i9-10900X CPU to execute an entire step, i.e. a batch containing all the training samples.

\section{Computational Results}
\label{sec:Res}
We present in this chapter our computational \unoref{experiments} divided into three blocks.
In \Cref{sec:opt_default}, we explain how we have tested L-BFGS and the eight optimizers briefly described in \Cref{sec:OptAlg} on the {baseline problem, i.e., training
the \textsc{Baseline} architecture
on
UC Merced dataset,} using default setting of the hyperparameters.
We have also carried out  {a} multistart test on the three worst-performing algorithms (Adadelta, Adagrad, and FTLR) to assess whether poor performances were caused only by an unfortunate weights initialization or by the inherent behaviour of these optimizers on the dataset. 
We discuss the correlation between the test accuracy performances and the precision with which the problem in \Cref{eq:loss} is solved, the loss profiles produced by the algorithms, as well as the role of  {the} data augmentation technique.
\unoref{Our further analysis in the following sections is focused} only on five of these optimizers since, as shown in \Cref{sec:opt_default}, they achieved the highest accuracy prediction.
In \Cref{sec:grid_search}, we describe the grid search we have carried out \unoref{on the baseline network on UC Merced dataset} to tune optimizers' hyperparameters.
We show how a careful tuning aimed at finding a nearly-optimal hyperparameters setting can result in significant improvements in terms of  test accuracy.
In \Cref{sec:deep_wide}, we discuss the results obtained after modifying the network architecture with respect to the baseline, investigating in particular hyperparameters robustness to the increase in depth and width. 

Finally, \unoref{in \Cref{subsec:new_datasets}, we carry out tests on the two other image classification datasets,  CIFAR10 and CIFAR100 with $H=W=32$ pixels and $N=10, 100$ respectively.}

\begin{table}[h]
 \renewcommand{\arraystretch}{1.25}
       \centering
    \footnotesize
    \caption{Number of trainable parameters (variables) of the different neural architectures; they \unoref{depend on $H, W, N$}. The values are reported in millions [M].}
    \begin{tabular}{ l | c | c | c}
        Architecture & \multicolumn{3}{c}{\# variables  [M]} \\
  %      & 3x256x256 ($N=21$) & 3x32x32 ($N=10$) & 3x32x32 ($N=100$)\\
         & UCMerced & CIFAR10 & CIFAR100\\
       \hline
        \textsc{Baseline} & 4.84 & 0.72 & 0.74 \\
        \textsc{Wide} & 10.97 & 2.72 & 2.74\\
        \textsc{Deep} & 9.83 & 1.57 & 1.62\\
         \textsc{Deep} \& \textsc{Wide} & 22.52 & 5.99 & 6.05\\
        Resnet50 &  24.77 & 24.77 & 24.79 \\
        Mobilenetv2 & 3.05 & 3.05 & 3.07 \\
        
    \end{tabular}
    \label{tab:arch_param}
\end{table}

\subsection{The \textsc{Baseline} \unoref{problem} with default hyperparameters}
\label{sec:opt_default}
The first tests we have carried out are aimed at studying the optimizers' performances both from an optimization perspective (i.e., the value of the final loss) and a   {machine learning} perspective (i.e., the test accuracy). Hyperparameters have been set to their default values (\Cref{tab:def_values}), taken from the TensorFlow documentation.

\begin{table*}[h!]
    \footnotesize
    \centering
    \caption{Default and tuned values of the TensorFlow built-in optimizer hyperparameters. 
    The \unoref{tuned} values are reported in brackets when different from the default ones, 
%    When the \unoref{tuned} values found after the grid search (Section \ref{sec:grid_search}) differ from the default ones, they are reported in brackets. 
   The grid search has not been performed for the worst-performing algorithms: Adadelta, Adagrad, and FTRL.}
    \label{tab:def_values}
    \resizebox{\textwidth}{0.11\textwidth}{\begin{tabular}{l r r r r r  r  r r r}
    \hline
     Algorithm & SGD & Adam & Adamax & Nadam & RMSProp & Adadelta & Adagrad & FTRL \\ \hline \\
    $\eta$ & $10^{-2} (10^{-1})$ & $10^{-3}$ & $10^{-3}$ & $10^{-3}$ & $10^{-2}\ (10^{-3})$ & $10^{-3}$ & $10^{-3}$ & $10^{-3}$ \\
    $\beta$ & 0 (0.9) & - & - & - & 0 & - & - & 0.1 \\
    $\beta_1$ & - & 0.9 & 0.9 (0.6) & 0.9 (0.99) & - & - & - & 0 \\
    $\beta_2$ & - & 0.999 (0.9999) & 0.999 (0.99) & 0.999 (0.99) & - & - & - & 0 \\
    $\epsilon$ & - & $10^{-7}\ (10^{-8})$ & $10^{-7}\ (10^{-6})$ & $10^{-7}\ (10^{-6})$ & $10^{-7}\ (10^{-6})$ & $10^{-7}$ & $10^{-7}$ & - \\
    Amsgrad & - & False (True) & - & - & - & - & - & - \\
    $\rho$ & - & - & - & - & 0.9 & 0.95 & 0.95 & - \\
    Centered & - & - & - & - & False & - & - & - \\
    Nesterov & False & - & - & - & - & - & - & - \\
    \hline
    \end{tabular}}
\end{table*}
%\todo[inline]{I tuned parameter per SGD non ci sono. Devono per forza essere diversi dai default perché i risultati sono tutti diversi}
\dueref{For each algorithm, we report} in \Cref{fig:loss_opt_def} \dueref{the behaviour of } the training losses without and with data augmentation,a and in \Cref{tab:def_accuracy} the test accuracy.

\begin{figure*}[h!]
    \centering
%    \includegraphics[width=\linewidth]{lables.pdf} 
 %   \hfill
    \subfigure[W/out data augmentation]{\includegraphics[width=0.47\linewidth]{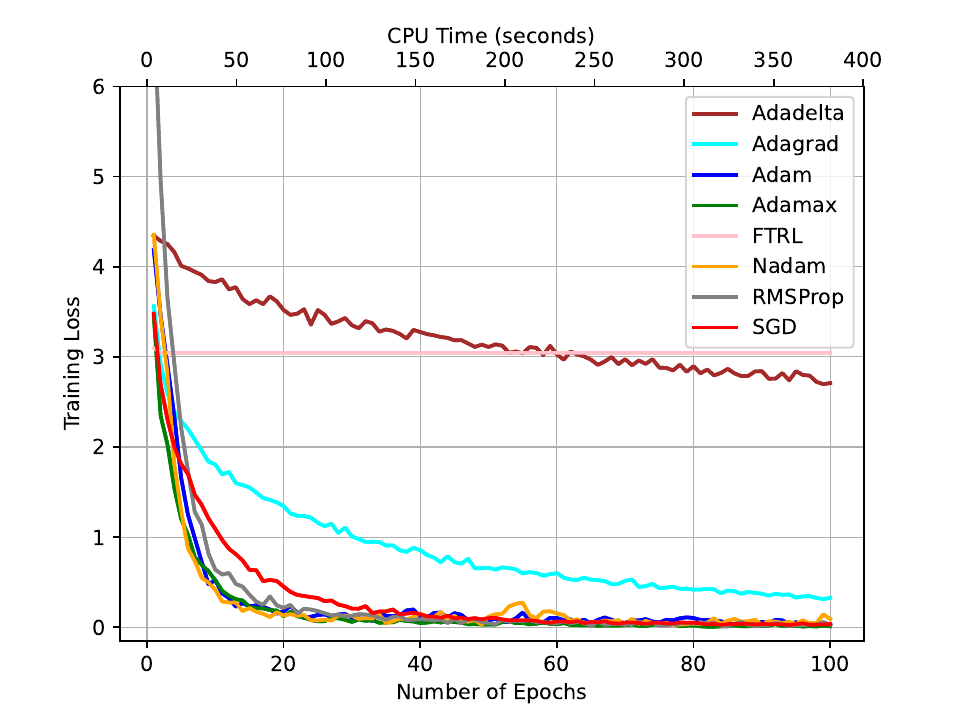}
    \label{img:loss_opt_def_no_aug}
    } \hfill
    \subfigure[with data augmentation]{\includegraphics[width=.47\linewidth]{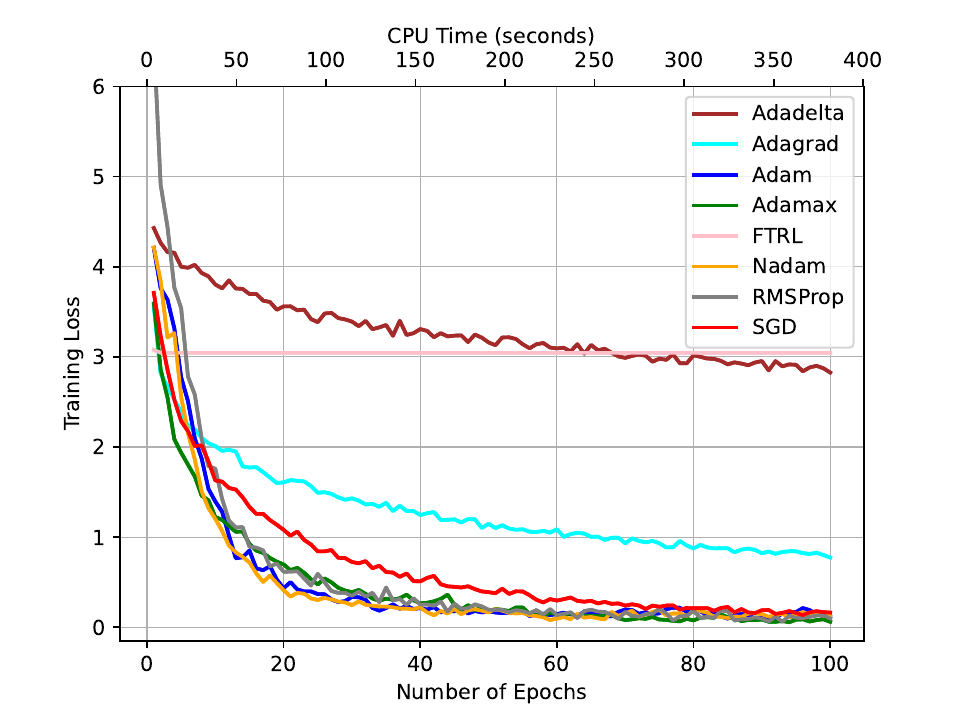}
    \label{img:loss_opt_def_aug}
    }
    \caption{{Training loss} for the \textsc{baseline} problem with default hyperparameters setting: (a) without data augmentation and (b) with data augmentation. A colour is assigned to each algorithm according to the legend. 
}
    \label{fig:loss_opt_def}
\end{figure*}
% \todo[inline]{ATTENZIONE: ho tolto la legenda prché era inutille. PIU' IMPORTANTE: le figure mi sembra fossero invertite a) è senza e b) è con DA. CHECK !! Poi, se vero, check nel testo }

% \todo[inline]{SIETE SICURI CHE QUESTE SIENAO LE FIGURE CON IL DEFAULT ???? il nome sembra quello con "opt" }

We noticed that, while most of the algorithms converge to points in the neighbourhood of the globally optimal solution, i.e., the training loss is close to   zero, Adadelta, Adagrad, and FTRL get stuck in some local minima, as can be seen in \Cref{fig:loss_opt_def}.
This results in quite poor accuracy performances for Adadelta, Adagrad, and FTRL \dueref{as reported in the first column of \Cref{tab:def_accuracy}.} 
{We highlight that FTRL is an extreme case, being the descent so slow that the loss profile looks  {like} an horizontal line. 
This is not surprising because, as explained in \cite{McMahan2013AdCP}, FTRL has been thought to deal with extremely sparse datasets, which is not the case for colored images.
Furthermore, FTLR convergence requires  {a} very strong convexity condition (\cite{mcmahan2011follow}), making  {it} impossible to predict its behaviour in such a non-convex context.}

% \begin{table}[t]
%     \centering
%     \caption{Test accuracy in \% obtained with default values and with tuned values ($\pm$ \% increase) of the optimization algorithms, without and with data augmentation..
%     % with default values ($\pm$ \% increase with tuned values) of the optimization algorithms, without and with data augmentation.
%     }
%     \label{tab:def_accuracy}
%     \renewcommand{\arraystretch}{1.25}
%     \resizebox{0.7\textwidth}{!}{
%     \begin{tabular}{l l l}
%     \hline
%      Algorithm&W/out data Augmentation&With Data Augmentation \\  \hline
%      Adam& 60.0 \% (+2.1 \%) & 72.4 \% (+2.5 \%)\\
%      Adamax& 61.3 \% (+1.7 \%)& 72.5 \% (+0.0 \%)\\
%      Nadam& 61.3\% (+3.0 \%)& 72.1 \% (+2.6 \%)\\
%      RMSProp& 60.2\% (+3.1 \%)& 74 \% (-1.9 \%)\\
%      SGD& 59.4 \% (+1.8\%)& 65.1 \% (+4.4 \%)\\  \hline
%      Adadelta & 17.6 \% &  18.0 \%\\
%      Adagrad& 32.7 \%& 31.0 \%\\
%      FTRL& 4.6 \%& 4.6 \%\\  \hline
%     \end{tabular}}
% \end{table} 

\begin{table}[t]
%\todo[inline]{ho girato la tabella per essere congruente con sopra e recuperare spazio}
    \centering
    \caption{Test accuracy in \% obtained with default values and with tuned values ($\pm$ \% increase) of the optimization algorithms, without data augmentation (W/out DA) and with data augmentation (with DA).
    % with default values ($\pm$ \% increase with tuned values) of the optimization algorithms, without and with data augmentation.
    }
    \label{tab:def_accuracy}
  \renewcommand{\arraystretch}{1.25}
    \resizebox{\textwidth}{!}{
    \begin{tabular}{l | l l l l l l l l}
    \hline
     Algorithm	&	     Adam	&	     Adamax	&	     Nadam	&	     RMSProp	&	     SGD	&	     Adadelta 	&	     Adagrad	&	     FTRL	\\ \hline
W/out DA	&	 60.0 \% (+2.1 \%) 	&	 61.3 \% (+1.7 \%)	&	 61.3\% (+3.0 \%)	&	 60.2\% (+3.1 \%)	&	 59.4 \% (+1.8\%)	&	 17.6 \% 	&	 32.7 \%	&	 4.6 \%	\\
With DA   	&	 72.4 \% (+2.5 \%)	&	 72.5 \% (+0.0 \%)	&	 72.1 \% (+2.6 \%)	&	 74 \% (-1.9 \%)	&	 65.1 \% (+4.4 \%)  	&	  18.0 \%	&	 31.0 \%	&	 4.6 \%  	\\
\hline
    \end{tabular}}
\end{table}

The observed behaviour seems to confirm what has been already pointed out in (\cite{swirszcz2016local,yun2018small}): neural networks can be affected by the local minima issue, which has a direct influence on the performance metrics.
Getting stuck in bad local minima often implies also an accuracy level that makes the entire network useless for the classification task. In the case of FTRL, the test accuracy is so low that the network selects randomly the predicted class.

Furthermore, data augmentation has no substantial effect in modifying the convergence endpoint. Indeed, Adadelta, Adagrad, and FTRL are somehow stable in returning the a bad point, as well as SGD, Adam, Adamax, Nadam, and RMSProp always converge to good solutions, leading to similar values of accuracy, as we can see again in \Cref{tab:def_accuracy}. 
Nonetheless, data augmentation have a boosting effect on test accuracy for all five working algorithms. 
This improvement is obtained because data augmentation {artificially increases the diversity and the quantity of the training data} and, thus, enhances the network generalization capability. \unoref{Nonetheless, data augmentation also makes the task harder and thus the the training loss decrease is slightly slower, i.e., the network needs more time to learn.}

In order to investigate the behaviour of Adadelta, Adagrad, and FTRL and to assess the stability of their bad performance, we have carried out another test employing a multistart procedure. \dueref{To this aim we have initialized the weights using }
%We have used 
two different %weights initialization 
distributions Glorot Uniform (GU) \citep{glorot2010understanding} and Lecun Normal (LN) \citep{lecun1989generalization} and %ran the algorithms with 
16 different seed values, i.e., starting from 16 different initial points \unoref{for each initialization, that is from 32 different points in total}. 
% \todo[inline]{non abbiamo una citazione per le due dostribizioni ? Sono standard ?}
% \todo[inline]{sono 16 per ogni distribuzione, qunidi  32 ? CHECK e cambia in caso. INVECE PER L-BFGS abbaimo usato solo 5 seeds ?}
The three algorithms always get stuck in a point, \dueref{with value of the training loss} quite far from zero with respect to the others.
This behaviour is very stable and does not change with the initialization seeds.

Best accuracy values, not reported in a table for the sake of brevity, are always $4.6\%$ for FTRL, $18.2\%$ for Adadelta, and $32.7\%$ for Adagrad.
% \todo[inline]{IN tabella delle accuracy,  2 algortimi (eccetto Adagrad) sono di norma sotto il 20\%. Con multistart erano migliorati così tanto ?}
These results suggest that %Adadelta, Adagrad, and FTRL do not fit with our task and that their 
the
bad behaviour of Adadelta, Adagrad, and FTRL is not just caused by an unfortunate initial point. These algorithms seem to converge to points which are not  {good}  for our classification task.
Hence, we have discarded Adadelta, Adagrad, and FTRL from the testing phases reported in the next sections.

%\todo[inline]{aggiunto che siamo senza DA}
Concerning L-BFGS, \dueref{we use the original dataset without data augmentation (which is specific for mini-batch methods). We have first trained  the baseline problem using the ReLU as activation function, as well as the MaxPooling \Cref{eq:CNN} and the Dropout operations}. 
Since the final points returned by L-BFGS are influenced by the starting point (\cite{liu1989limited}), we ran the algorithm with different initialization seeds. \dueref{In particular, we have used again the Glorot Uniform and Lecun Normal distributions and, due to the heavy computational effort, only 5 different seed values for each initialization}.
The training loss profile \dueref{of this first set of experiments} are reported in  \Cref{img:loss_ L-BFGS}.
%\todo[inline]{perchè "often fails" ? a me sembra che faòllisca sempre}
We observe that the algorithm \unoref{always} fails before achieving convergence: \unoref{the} lines in \Cref{img:loss_ L-BFGS} stop because the returned loss was infinite at a given iteration. 
We argue that this is caused by a non-differentiability issue. \dueref{Indeed, as we already discussed,} 
L-BFGS convergence is guaranteed exclusively when the objective function is continuously differentiable (\cite{liu1989limited}) and the ReLU, as well as the MaxPooling operation \Cref{eq:CNN}, cause the occurrence of non-differentiable points, i.e., points where the gradient is not defined.
Although, this could in principle happen with any other algorithm, since L-BFGS is a full-batch method{,} once a non-differentiable point is reached the algorithm gets stuck.
%and computing every time the whole gradient to perform a single update step can make  {it} easier to end up in a non-differentiable point.

Hence, we have also trained the baseline network \dueref{in a more differentiable setting, namely} using the  SiLU activation function, the MeanPooling \dueref{and disabling the Dropout operation}.
As we show in \Cref{img:loss_ L-BFGS_NoMPRD}, results significantly improved with almost all the initialization seeds. However, the loss does not always tend to zero and L-BFGS generally converges to points with a far worse loss function value than Adam, Adamax, Nadam, RMSProp, and SGD. 

\dueref{This difference can be seen also in terms of test accuracy. Indeed, when using}
%If we use
the MaxPooling layers and the ReLU activation function, L-BFGS performs quite poorly in terms of test accuracy reaching the maximum value of $15.4 \%$.

When \dueref{using the "differentiable" setting}, we obtain better results reported in 
 \Cref{tab:loss_L-BFGS}. We report the \unoref{average over the 5 runs of the} final \dueref{training} loss values and \dueref{of} the test accuracies.
 %obtained with L-BFGS in this last case. 
% \todo[inline]{sono risulttai on average sui 5 runs per le due differetmti distribuzioni ?}
We observe that even the most unlucky initialization, which is GU 10, results in a test accuracy of 15,6\% which is better than the highest one obtained with ReLU and MaxPooling.
However, L-BFGS is much more sensible to the initialization seed with respect to the other built-in methods, confirming claim ii). 
The final training loss, as well as the test accuracy, are not stable and may vary in a wide range of values. 
Despite this computational result could question the practical effectiveness of traditional batch methods in deep learning, it also confirms our claim i): the quality of local minima matters.
Indeed, looking at \Cref{tab:loss_L-BFGS}, we observe a \unoref{relation} between the final loss value and the test accuracy: 
lower final loss value usually corresponds to higher test accuracy. 
\dueref{In general, the accuracy performances achieved are not satisfactory when compared to mini-batch methods as well as the training loss decrease in unstable and highly influenced by the starting point.}
Finally, we also remark that L-BFGS \unoref{is} significantly less efficient with respect to the other built-in algorithms. \unoref{Indeed, it is practically impossible to run it on a standard GPU, because one needs enough memory storage to access the entire dataset in one single step, which is possible only on CPU and this results in slower training.
%\todo[inline]{verifica se ho riscritto bene la storia delle GPUs}
}
% \todo[inline]{come si vede che è meno efficiente ? E poi perché "sembra" ? DEVE essere meno effciente....}
\begin{figure*}[h!]
    \centering
    \includegraphics[width=0.6\linewidth]{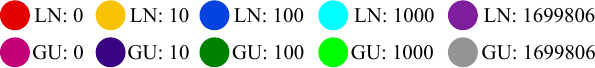}
    \hfill
    \subfigure[Using ReLU activation function, MaxPooling and \dueref{Dropout operations}]{\includegraphics[width=0.8\linewidth]{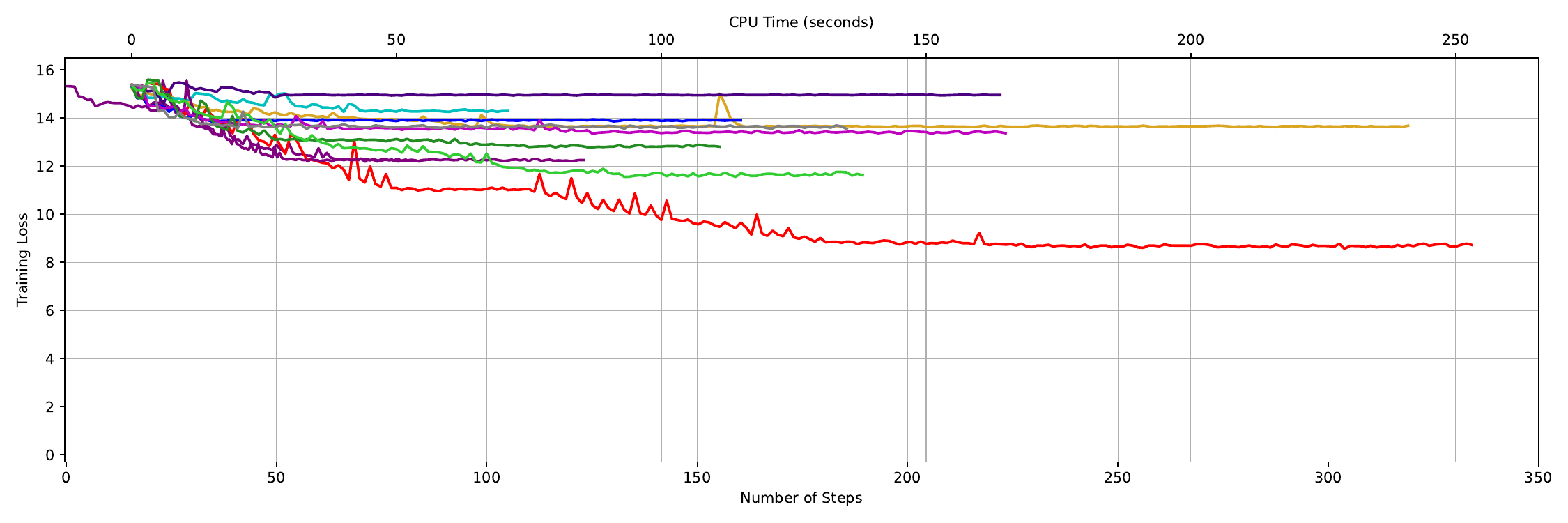}
    \label{img:loss_ L-BFGS}
    }
    \subfigure[Using SiLU activation function, MeanPooling and without Dropout and using Conv with 2 strides]{\includegraphics[width=0.8\linewidth]{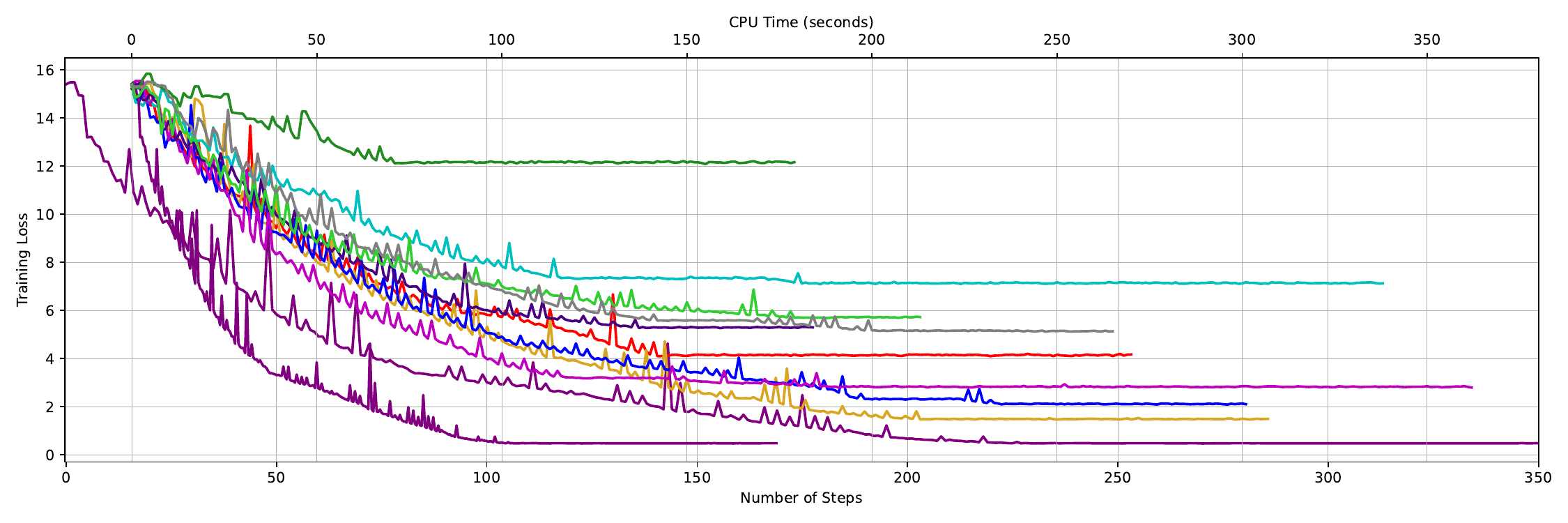}
    \label{img:loss_ L-BFGS_NoMPRD}
    }
    \caption{L-BFGS training loss with 5 different initialization seeds as reported in the legend, using the Lecun Normal (LN) and Glorot Uniform (GU) initialization procedures.}
    \label{fig: L-BFGS}
\end{figure*}
 % \todo[inline]{RIBADISCO perchè la eisposta NON mi convince. Nella caption figura b) si parla di "Conv with 2 strides". Da nessuna parte nel testo sta scritto che per LBFGS usiamo stride 2. Non mi pare che abbiamo mai specificato qual è lo stride per gli altri metodi. In 473 si parklla di stride $s\ge 1$ ma non si dice altro  }

% \begin{table}[ht]
%     \footnotesize
%     \centering
%     \caption{Avg training loss and Avg Test accuracy of L-BFGS \dueref{on the baseline problem using the "differentiable"" setting} with 5 different seeds \dueref{and the two different distributions GU and LN}: different training loss results in different test performances.}
%      \resizebox{0.55\textwidth}{!}{
%     \begin{tabular}{c c c c c c c c c c}
%     \hline
%     & Seed & Distribution & Avg training loss & Avg Test accuracy \\
%     \hline
%     & $0$ & LN & 4.17 & 37.6\% \\
%     & & GU & 2.80 & 51.3\% \\  \hline
%     & $10$ & LN & 1.49 & 43.8\% \\
%     & & GU & 5.30 & 15.6\% \\  \hline
%     & $100$ & LN & 2.10 & 27.6\% \\
%     & & GU & 12.16 & 16.5\% \\  \hline
%     & $1000$ & LN & 7.12 & 27.1\% \\
%     & & GU & 5.73 & 23.3\% \\  \hline
%     & $1699806$ & LN & 0.48 & 49.8\% \\
%     & & GU & 5.13 & 35.2\% \\
%     \hline
%     \end{tabular}}
%     \label{tab:loss_L-BFGS}
% \end{table}

%\todo[inline]{ho girato la tabella}

\begin{table}[ht]
%\todo[inline]{i run sono 10 giusto ? 5 per ogni dustribuzione}
    \footnotesize
    \centering
    \caption{Average value of the training loss (Avg Train value) and Average Test accuracy in \% (Avg Test acc) obtained  on the baseline problem with L-BFGS
        \dueref{in the "differentiable" setting} on 10 multistart runs using 5 different seeds \dueref{for the two different distributions GU and LN}: different training loss results in different test performances.}
  \renewcommand{\arraystretch}{1.25}
     \resizebox{ \textwidth}{!}{
    \begin{tabular}{l | c c |c c |c c |c c |c c }
    \hline
 Seed 		&	\multicolumn{2}{c|}{$0$ } 	 	&	 \multicolumn{2}{c|}{$10$}	 	&	 \multicolumn{2}{c|}{$100$}	 	&	 \multicolumn{2}{c|}{$1000$}		 	&	 \multicolumn{2}{c}{$1699806$}	 	\\
 Distribution 		&	 LN 	&	 GU 	&	 LN 	&	 GU 	&	 LN 	&	 GU 	&	 LN 	&	 GU 	&	 LN 	&	 GU 	\\\hline
Avg Train value		&	 4.17 	&	 2.80 	&	 1.49 	&	 5.30 	&	 2.10 	&	 12.16 	&	 7.12 	&	 5.73 	&	 0.48 	&	 5.13 	\\
 Avg Test acc		&	 37.6\% 	&	 51.3\%   	&	 43.8\% 	&	 15.6\%   	&	 27.6\% 	&	 16.5\%   	&	 27.1\% 	&	 23.3\%   	&	 49.8\% 	&	 35.2\% 	\\
    \hline
    \end{tabular}}
    \label{tab:loss_L-BFGS}
\end{table}

%\subsection{The \textsc{Baseline} \unoref{problem} with tuned hyperparameters}
%\subsection{Impact of hyperparameters tuning on the \textsc{Baseline} \unoref{problem}}
\subsection{Impact of  tuning on the \textsc{Baseline} \unoref{problem}}
\label{sec:grid_search}
\dueref{In this section, we perform tuning of hyperparameters of the optimization algorithm  on the baseline problem, namely on the \textsc{baseline} architecture and the UC Merced problem, to assess their role in the computational efficiency and, in turn, on the final test accuracy. As we mentioned in the introduction, default values for hyper-parameters are often obtained by maximizing the aggregated (in most cases the average) performance across a variety of different tasks, balancing a trade-off between efficiency and adaptability to different datasets (\cite{probst2019tunability,yang2020hyperparameter,bischl2023hyperparameter}). We aim here to assess if 
a specific  tuning on the classification task 
 has a influence on algorithms' behavior. As this will be the case, in the next section we analyse  the impact of the tuning obtained on a baseline problem  to other settings (architecture and/or dataset).
}
%%%%%%%%%%%

%%%%%%%%%%%
We discard Adadelta, Adagrad, FTRL, and L-BFGS from  further analysis due to their extremely poor performance on the baseline problem. Thus, we have carried out a grid search to tune the hyperparameters of Adam, Adamax, Nadam, RMSProp, and SGD \dueref{on the \textsc{baseline} problem}. 

The grid search ranges are reported in \Cref{tab:grid_search_values}. 
\dueref{Concerning numerical hyperparameters, we have chosen ranges centered in the default values, resulting in almost 200 possible combinations for each algorithm.
 We did not perform either a $k$-fold cross-validation or a multistart procedure, as is usually the case in computer vision (see e.g., (\cite{gartner2023transformer})) for  computational reasons.
Indeed, each run (including loading the dataset and the network to the GPU and the effective computational time) takes approximately 9 minutes (around 2 seconds per epoch plus the set-up time). Thus, performing the grid search over about 200 hyperparameter settings takes 1.25 days  on a fully dedicated 12GB NVIDIA GTX TITAN V GPU for each of the five algorithms. 
This implies that each training phase with a complete grid search would require nearly 6.25 days. Therefore, 
performing, e.g. a $k$-fold cross-validation would require $6.25 \cdot k$ days on a fully dedicated machine, being a prohibited amount of time for standard values of $k$ (5 or 10). Similar observations holds for a multistart procedure.}
% \todo[inline]{perché ci vogliono 10 giorni se abbiamo solo 5 ottimizzatori ? a me viene 6,25 giorni CHECK}

% \todo[inline]{Aggiungere come abbiamo selezioonato i parametri ottimi. Abbiamo visto le migliori performance sul test (come abbiamo diviso trian test?) ? oppure i tempi computazinali ? Non è  riportato il criterio}

\begin{table*}[h!]
    \footnotesize
    \centering
    \caption{Range of the grid search for the different hyperparameters. {For the sake of brevity, we have omitted the boolean hyperparameters, \texttt{Amsgrad} for Adam, \texttt{Centered} for RMSProp, and \texttt{Nesterov} for SGD.}}
    \label{tab:grid_search_values}
    \resizebox{\textwidth}{0.075\textwidth}{
        \begin{tabular}{l c c c c c c c c c}
            \hline
            Algorithm & Adam & Adamax & Nadam & RMSProp & SGD \\ \hline
            $\eta$ & $10^{-i}, i=2,3,4$ & $10^{-i}, i=2,3,4$ & $10^{-i}, i=2,3,4$ & $10^{-i}, i=2,3,4$ & $10^{-i}, i=1,2,3,4$ \\
            $\beta$ & - & - & - & $\{0,0.5,0.9\}$ & $\{0,0.5,0.9\}$ \\
            $\beta_1$ & $\{0.6,0.9,0.99\}$ & $\{0.6,0.9,0.99\}$ & $\{0.6,0.9,0.99\}$ & - & - \\
            $\beta_2$ & $\{0.99,0.999,0.9999\}$ & $\{0.99,0.999,0.9999\}$ & $\{0.99,0.999,0.9999\}$ & - & - \\
            $\epsilon$ & $[10^{-i}, i=6,7,8]$ & $[10^{-i}, i=6,7,8]$ & $[10^{-i}, i=6,7,8]$ & $[10^{-i}, i=6,7,8]$ & - \\
            $\rho$ & - & - & - & $\{0.6,0.9,0.99\}$ & - \\
            \hline
        \end{tabular}
    }
\end{table*}

\dueref{The tuned values of the hyperparameters are selected considering the best test accuracy obtained and are reported into brackets in Table \ref{tab:def_values}, when different from default ones.}
Once tuned the hyperparameters to new values, we have used them on the \textsc{Baseline} problem halving the number of epochs. 
% \todo[inline]{cosa si intende per validation ? Abbiamo detto come abbiamo splittato il dataset in train /validation ? poi di solito si usa la parola validation quando si fa una k-fold e si fa la media. Oppure noi abbiamo a parte anche un test ?}
\begin{figure*}[h!]
     \centering
    \hfill
    \\
    \vskip -0.5truecm
    \centering
    \subfigure[W/out data augmentation]{\includegraphics[scale=0.34]{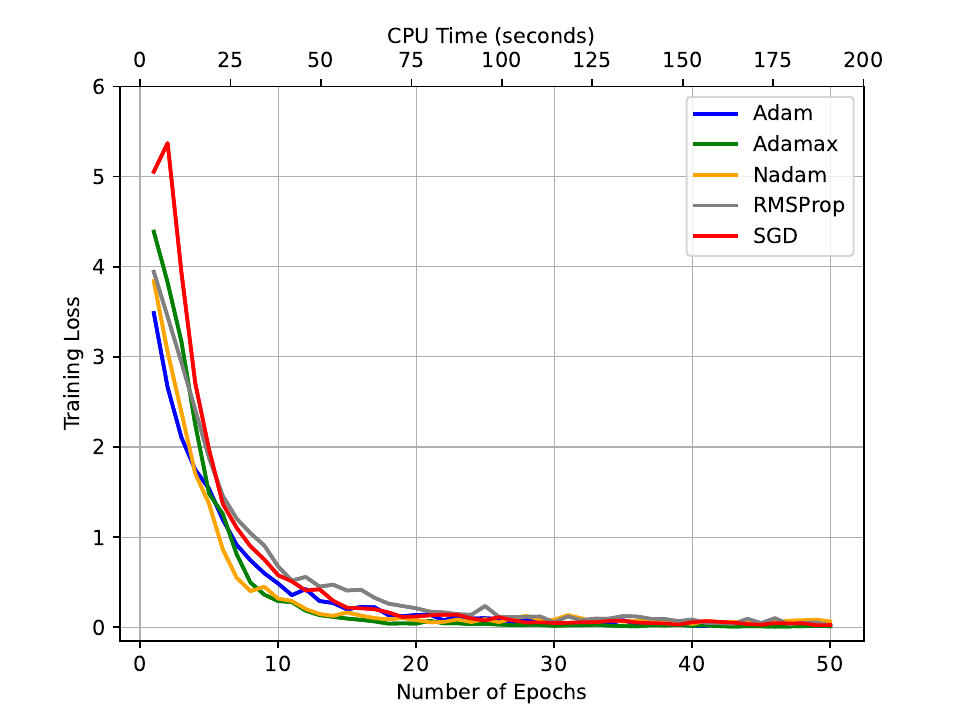}
    \label{img:loss_opt_optimal_no_aug}
    } \hfill
    \subfigure[with data augmentation]{\includegraphics[scale=0.34]{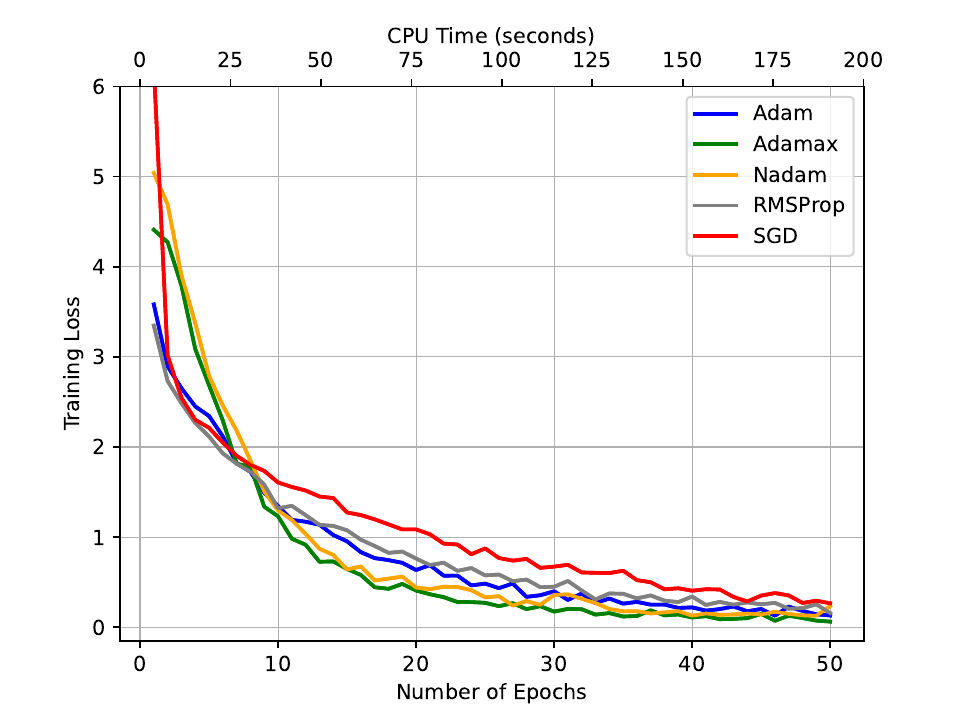}
    \label{img:loss_opt_optimal_aug}
    }
    \caption{Training loss with tuned hyperparameters values for the \textsc{baseline} problem: (a) without data augmentation,(b) with data augmentation.}
    \label{fig:loss_opt_optimal}
\end{figure*}

\dueref{We report in \Cref{img:loss_opt_optimal_no_aug} and \Cref{img:loss_opt_optimal_aug} the training loss profiles  
for the two settings without data augmentation and with data augmentation. Comparing with the corresponding training loss with default values in \Cref{img:loss_opt_def_no_aug}
and \Cref{img:loss_opt_def_aug}, we can state that}
tuning does not directly influence \unoref{the final value of the objective function returned by the algorithms, which was already the global optimal value near to zero.}
However, the loss decrease is faster, and nearly optimal values (near to zero) are reached \unoref{earlier, optaining good results}  despite having halved the number of epochs. 
\unoref{In particular, in \Cref{img:loss_opt_optimal_no_aug} the loss is almost zero already after 15-20 epoch, while in \Cref{img:loss_opt_def_no_aug} after 25-30 epochs.}
\unoref{Considering the case with data augmentation, in \Cref{img:loss_opt_def_aug} the loss is almost zero after 40 epochs in \Cref{img:loss_opt_optimal_aug} the same is true after approximately 60 epochs.}
% \todo[inline]{Questa affermazione non è del tutto corretta vedendo le figure di default, ammesso che siano giuste. INfatti anche quelle arrivavamo quadsi a zero intorno a 4o epoche. Qui però noi ci arriviamo molto prima. bisognerebbe quantiticare quanto prima (in termoin di epoche). Cioè non conta il numero di epoche totali, ma quanto prima potremmo interrompere.
% A me sembra che già ntorno a 20 epoche stiamo a zero}
%We are, in fact, obtaining the same results with  {half the computational time}.
% \todo[inline]{dove si vede che va meglio sul test ?}

However, our computational experience shows that the most relevant benefit of hyperparameters tuning is the gain in terms of test accuracy.
In \Cref{tab:def_accuracy} we report in square brackets the test accuracy \dueref{changes} for the five optimizers, with and without data augmentation.\dueref{The change is always positive except for RMSProp with data augmentation. We also observe that Adam, SGD and Nadam show a larger improvement over Adamax and RMSProp both w/out and with data augmentation.} 

\subsection{\dueref{Impact of  tuning  when changing the architectures}}
\label{sec:deep_wide}

In this  section,
\dueref{we aim to investigate the impact of tuned vs default hyperparameters} when changing the network architecture \dueref{whilst the dataset is UC Merced.}

In particular, we are interested in assessing how optimizers react to the increase in depth and width,
\dueref{using the three synthetic configurations (\textsc{Wide}, \textsc{Deep},   \textsc{Deep}\&\textsc{Wide}) described in \Cref{sec:NetArc},}
as well as in determining the impact of tuning the hyperparameters
\dueref{on state-of-the-art architectures as Resnet50, and Mobilenetv2.}
%Specifically, we consider the three configurations (\textsc{Wide}, \textsc{Deep},   \textsc{Deep}\&\textsc{Wide}) described in \Cref{sec:NetArc}{, Resnet50, and Mobilenetv2}.

% Since in the former experiments data augmentation had lead always to better performances in terms of final accuracy, we have carried out this test only using data augmentation. 
%\todo[inline]{Togliere frase sopra e Aggiugererei una frase del tipo:}
\unoref{In this experiment, we perform a single run for each algorithm starting from the same initial point, fixing the distribution and the random seed to LN 1699806, and considering the use of data augmentation, which gave better results in the former experiments.}

\dueref{The results for the synthetic architectures are reported in \Cref{fig:Wide_Def_Best} (training results) and in \Cref{tab:DW_accuracy} (test accuracy), whereas the results for the state-of-the-art architectures are in  \Cref{fig:stateoftheart_Def_Best} (training) and \Cref{tab:RnMob_accuracy} (test). 
In \Cref{tab:avg_diffs}, we report a cumulative difference in test accuracy (TEST ACC) when using Tuned versus Default hyperparameters setting on average for the Synthetic Networks.}

\dueref{
In terms of training loss, on the synthetic networks, the tuned version seems to reach, on average, slightly smaller values, whereas on the state-of-the.art architectures, there are not noticeable differences in the reached value. Thus, the tuning of the hyperparameters does not improve significantly the decrease rate.} 

% In the \textsc{Wide} configuration, we notice that the \dueref{training} loss, reported in \dueref{(a) and (b) of} \Cref{fig:Wide_Def_Best},  tends very rapidly to zero both using the default  and the \dueref{tuned} hyperparameters setting, while in the \textsc{Deep} (\dueref{(c) and (d) of} \Cref{fig:Wide_Def_Best}) 
% \textsc{Deep}\&\textsc{Wide} (\dueref{(e) and (f) of} \Cref{fig:Wide_Def_Best}%fig:DeepWide_Def_Best}
% ) loss decrease is slower, which is related to the increased complexity of the optimization problem.

%While in terms of loss decrease, we do not observe significant changes after hyperparameters tuning,

\dueref{As regard the test accuracy on synthetic networks, in  \Cref{tab:DW_accuracy}  we report
%we have different impacts as reported in  \Cref{tab:DW_accuracy}  for synthetic networks and  in  \Cref{tab:RnMob_accuracy} for Resnet50 and Mobilenetv2. 
for each algorithm the \% accuracy obtained for UCMerced dataset on the four different architectures and also the average \% over the architectures (column Avg ARCH).}

\dueref{
From  \Cref{tab:DW_accuracy}, it seems that
SGD and Adamax benefit from the tuned setting on the synthetic architectures, significantly improving the average of the \% accuracy (Avg ARCH), whereas Adam and Nadam are slightly worse on average. RMSProp deteriorates significantly, but we remark that this was the only case of worst performance also in the \textsc{Baseline} problem with data augmentation (see \Cref{tab:def_accuracy}).
These results are confirmed by the absolute difference between the tuned and default average test accuracy reported in \Cref{tab:avg_diffs} where SGD and Adam obtain the higher and significant increase.}

\dueref{The accuracy on Resnet50 and Mobilenetv2 are reported in \% in \Cref{tab:RnMob_accuracy} and with absolute variation in 
\Cref{tab:avg_diffs}.}
On these architectures, the tuned configuration does not perform uniformly better. However, we observe an improvement when using SGD  \unoref{on Mobilenetv2, which is more similar in the architecture to the \textsc{Baseline} architecture. Thus, we can conclude that when the architecture is significantly different from the \textsc{Baseline} used for tuning Hyperparmeters, the advantages are limited. } 

% In \cref{tab:avg_diffs}, we report the differences in test accuracy of Tuned versus Default hyperpameters
% setting for Resnet50, Mobilenetv2 and on average for the four Synthetic Networks.
%The results confirm a positive impact of tuning on SGD and Adamax.

% The test accuracy  results for the synthetic networks in \Cref{tab:DW_accuracy} confirm instead the robustness of the tuned configuration.
% We notice that SGD and Adamax perform quite poorly in \textsc{Deep} and  \textsc{Deep}\&\textsc{Wide} configurations with default hyperparameters setting, while test accuracy is still significantly higher with tuned hyperparameters.
% Concerning Adam, Nadam, and RMSProp, the default configuration only slightly outperforms the tuned one.
% Nonetheless, with both configurations we observe that final test accuracy is affected by the depth of the network.

% (see \dueref{(c) and (d) of} \Cref{fig:stateoftheart_Def_Best}).\unoref{Conversely, with Resnet50 in 
% \dueref{(a) and (b) of} \Cref{fig:stateoftheart_Def_Best}
% we can see how, despite reaching similar values of the loss function, the training appears to be slower with tuned configurations.}

\begin{figure*}[h!]   
   \centering
    \hfill
    \\
    \vskip -0.5truecm
     \centering
    \subfigure[\textsc{Wide} - Default HP]{\includegraphics[scale=0.35]{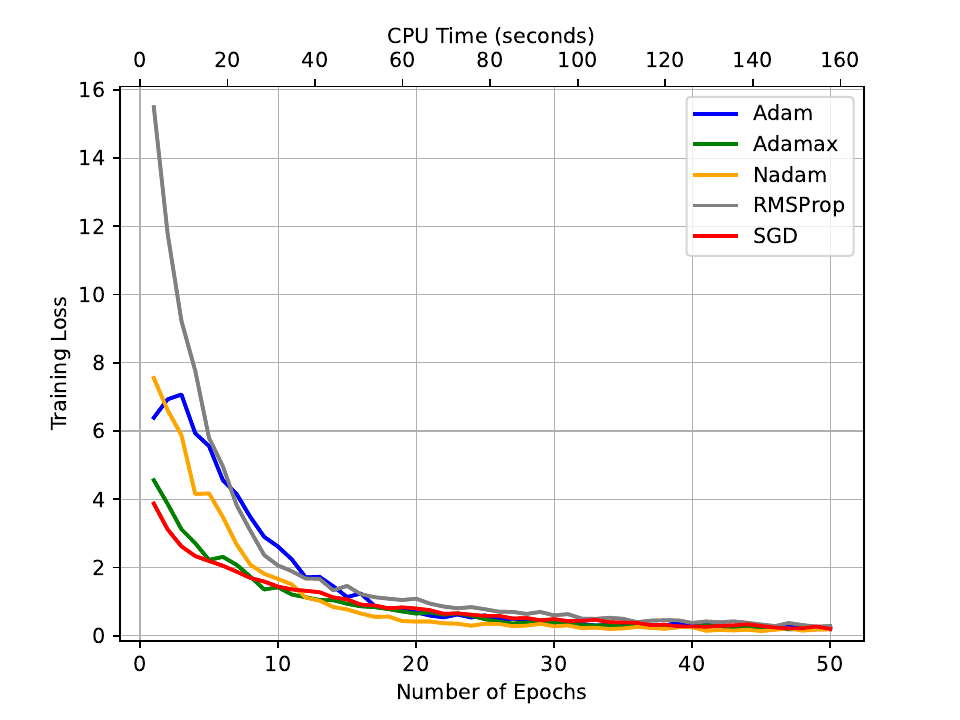}
    } \hfill
    \subfigure[\textsc{Wide} - Tuned HP]{\includegraphics[scale=0.35]{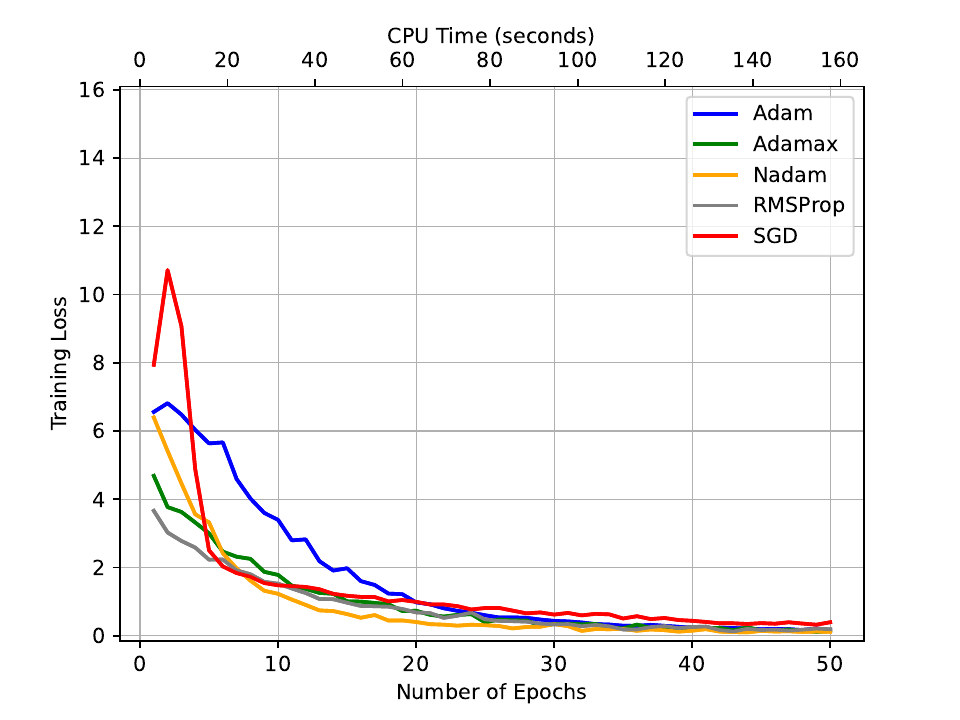}
    } \hfill
    %\caption{\dueref{Training} Loss trends of the 5 algorithms on UC Merced in the \textsc{Wide} configuration.% Colours are assigned to optimizers according to the labels above}
    \hfill
 \\
  \vskip -0.25truecm
     \centering
    \subfigure[\textsc{Deep}  - Default HP]{\includegraphics[scale=0.35]{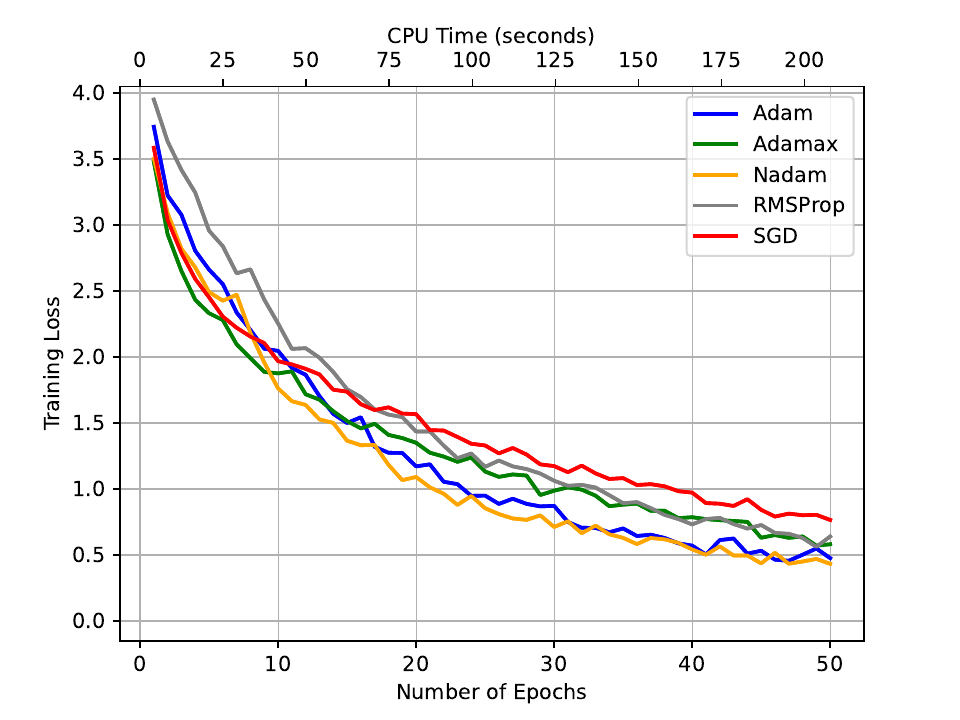}
    } \hfill
    \subfigure[\textsc{Deep}  - Tuned HP]{\includegraphics[scale=0.35]{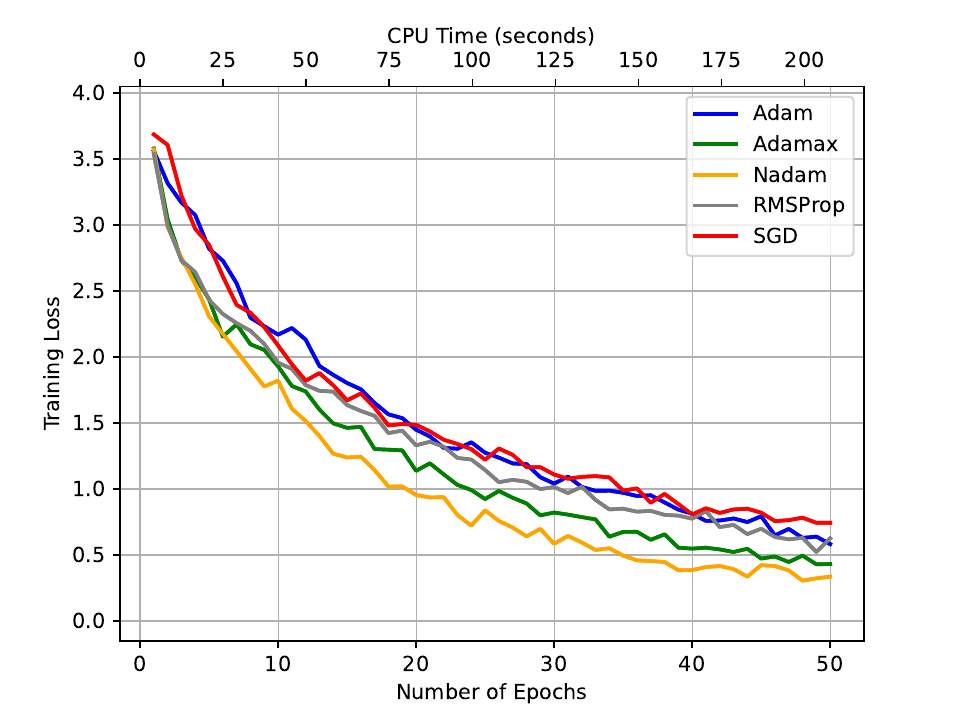}
    } \hfill
%    \caption{\dueref{Training} Loss trends of the 5 algorithms on UC Merced in the \textsc{Deep} configuration.    }
 %   \label{fig:Deep_Def_Best}
    \\
    \vskip -0.25truecm
    \centering
    \subfigure[ \textsc{Deep}\&\textsc{Wide} - Default HP]{\includegraphics[scale=0.35]{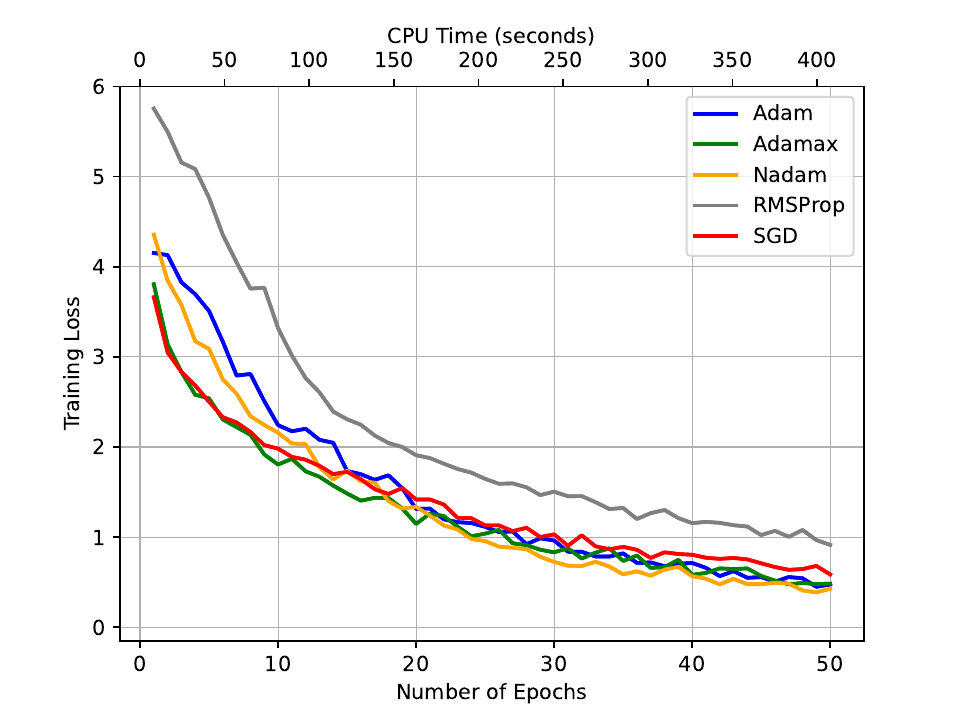}
    } \hfill
    \subfigure[ \textsc{Deep}\&\textsc{Wide} - Tuned HP]{\includegraphics[scale=0.35]{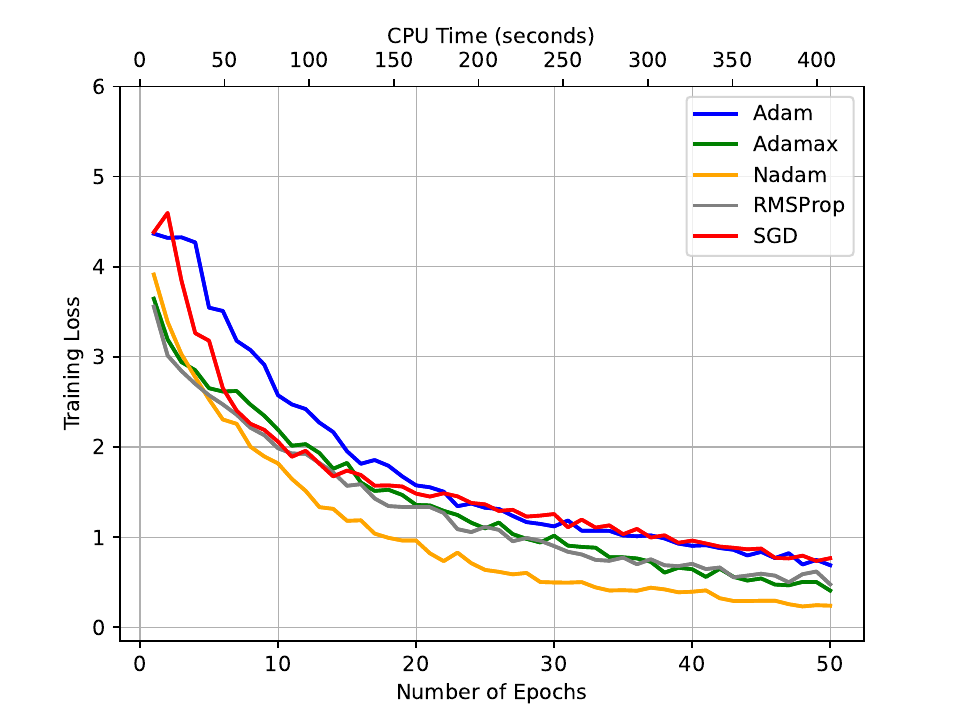}
    } \hfill
    \caption{
    \dueref{Training} Loss trends \dueref{of the 5 algorithms} for UC Merced on \dueref{the synthetic} architectures.}
%    \dueref{Training} Loss trends of the 5 algorithms on UC Merced in the  \textsc{Deep}\&\textsc{Wide} configuration.}
%    \label{fig:DeepWide_Def_Best}
\label{fig:Wide_Def_Best}
\end{figure*}

\begin{figure*}[h!]
    \centering
    \hfill
    \\
    \vskip -0.5truecm
    \centering
    \subfigure[Resnet50 - Default HP]{\includegraphics[scale=0.35]{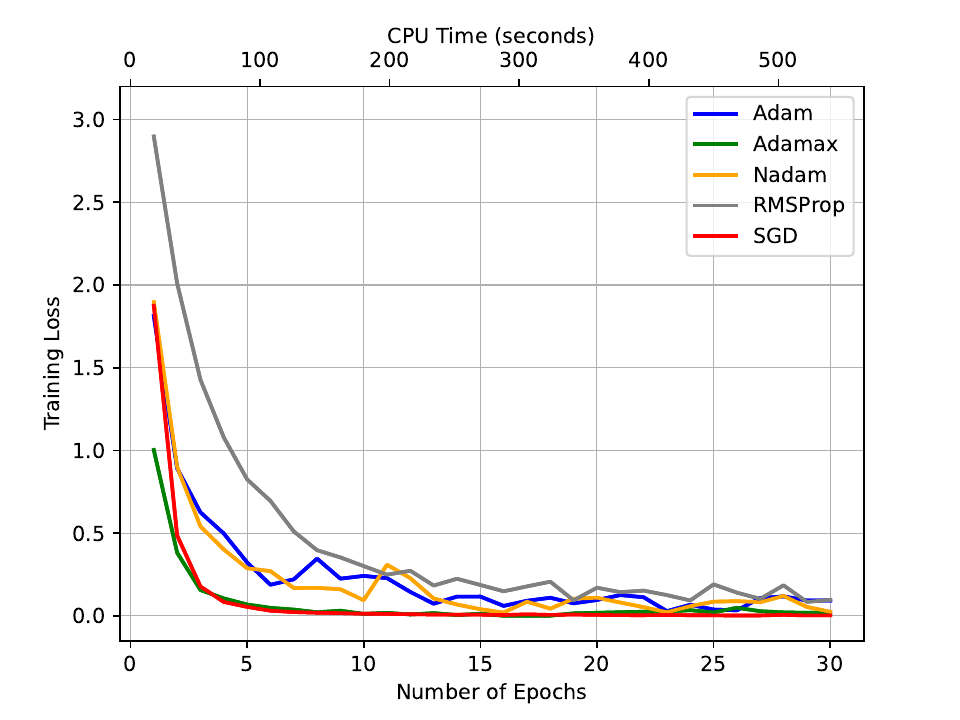}
    } \hfill
    \subfigure[Resnet50 - Tuned HP]{\includegraphics[scale=0.35]{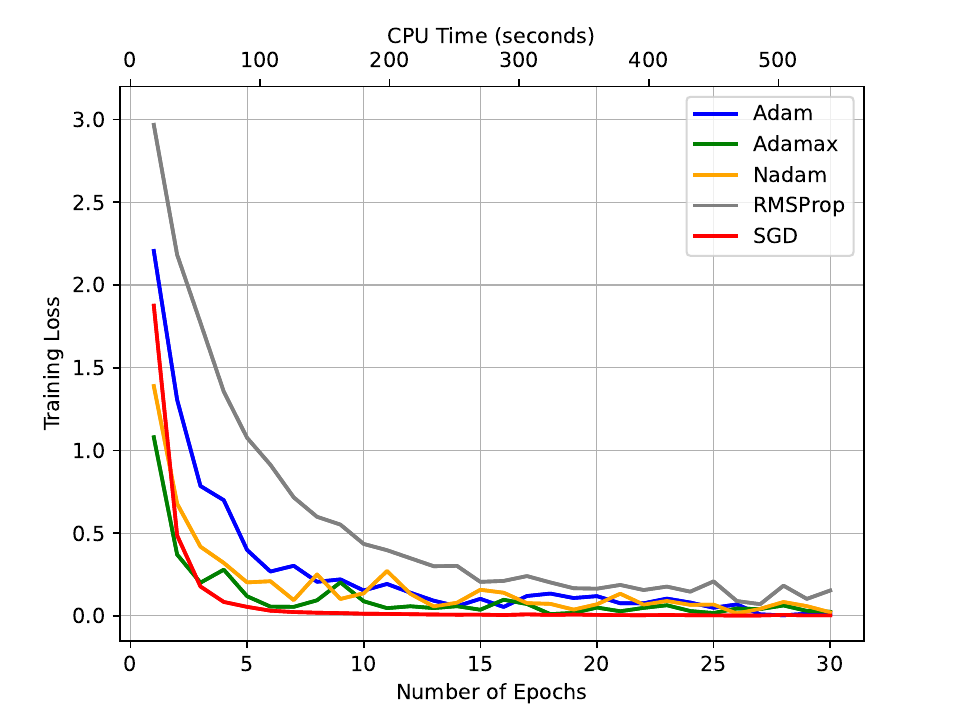}
    } \hfill
%    \caption{\dueref{Training} Loss trends of the 5 algorithms on UC Merced in the Resnet50 configuration.% Colours are assigned to optimizers according to the labels abo   }
    \label{fig:resnet50_Def_Best}
% \end{figure*}
% \begin{figure*}[h!]
%     \centering
%     %\includegraphics[width=0.5\linewidth]{labels5.pdf} 
%     \hfill
%     \\
    \vskip -0.25truecm
    \centering
    \subfigure[Mobilenetv2 - Default HP]{\includegraphics[scale=0.35]{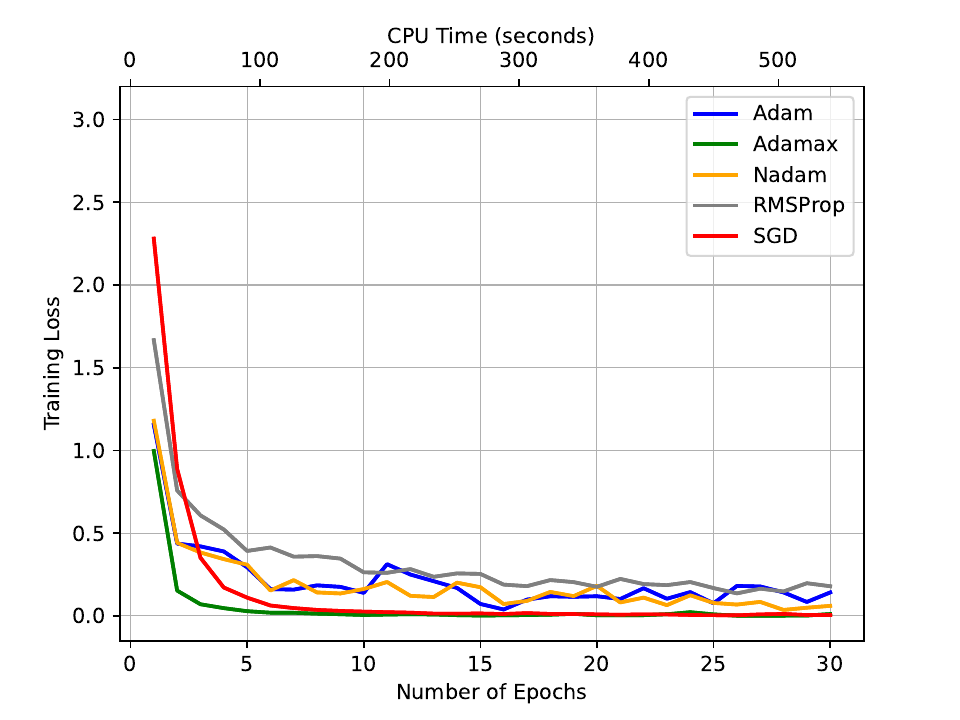}
    } \hfill
    \subfigure[Mobilenetv2 - Tuned HP]{\includegraphics[scale=0.35]{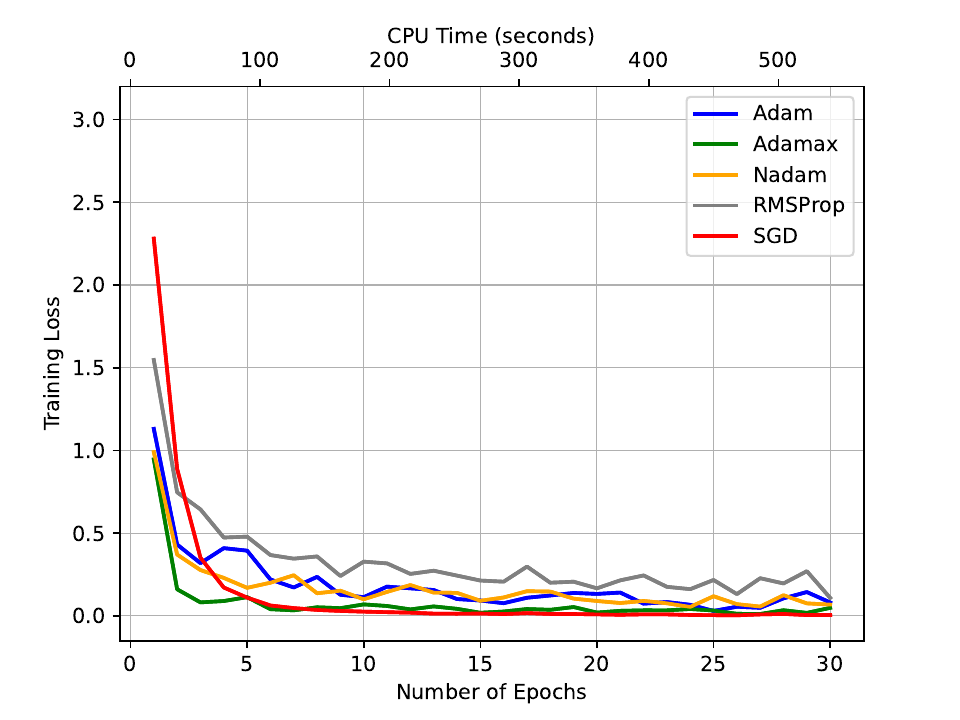}
    } \hfill
    \caption{\dueref{Training} Loss trends of the 5 algorithms for UC Merced on Resnet50 and Mobilenetv2. % Colours are assigned to optimizers according to the labels above.
    }
 %   \label{fig:mobilenetv2_Def_Best}
    \label{fig:stateoftheart_Def_Best}
\end{figure*}

\subsection{\dueref{Impact of tuning when changing the datasets}}
%{Assessing the impact of default versus tuned hyperparameters when changing the datasets}
%\subsection{Performance of all architectures on new datasets}
\label{subsec:new_datasets}

% \todo[inline]{Anche qui nelle figure default e tuned sembra che il punto inziale di partenza sia divesro ! Non so come spiegarmi alrtimenti il valore inioziale molto diverso. CONTROLLARE}
In this last set of experiments, we aim to assess the role of tuning 
when training all the architectures on the two additional datasets CIFAR10 and CIFAR100,  described  in \Cref{sec:dataset}.

% As regards CIFAR10, although it has a similar level of difficulty with respect to UC Merced, it
% is significantly larger.

\dueref{The training loss on CIFAR 10 
are reported in \cref{fig:cifar10_baseline} and \Cref{fig:cifar10_ResMb} for the synthetic architectures and for the state-of-the-art architectures  respectively, whereas the results on CIFAR100 are in
\Cref{fig:CIFAR100_baseline} and \Cref{fig:cifar100_ResMb}.
The test accuracies are reported in \% in \Cref{tab:DW_accuracy} and \Cref{tab:RnMob_accuracy} (test), and as the absolute difference between tuned and default version in \Cref{tab:avg_diffs}.
}

\dueref{Looking at the training losses both on CIFAR10 and CIFAR100, we do not observe significant differences in the profile of the training loss for most of the architectures among the two configurations, default or tuned one.}

%Looking at \Cref{fig:cifar10_baseline}, we do not observe significant differences in the profile of the training loss for most of the architectures.
%However, considering deeper architectures (in particular  \textsc{Deep}\&\textsc{Wide}, Mobilenetv2, and Resnet50), we notice that the loss profile is more stable using default configurations.

%Even stronger considerations can be made by looking at the loss profile on CIFAR100, which is a much harder task, consisting  {of} a 100-class multi-classification. 
%In \Cref{fig:CIFAR100_baseline} we observe that the default configuration leads to a smoother loss profile, especially for Adam and Nadam, as well as for deeper architectures. 
%Conversely, this behavior is much less evident in \Cref{fig:cifar100_ResMb}, where we do not observe particular differences between the default and tuned configurations.

Considering that in almost all the tests, regardless  {of} the architecture and the dataset, the final value of the training loss is not remarkably different between the tuned and the default case, one could conclude that the optimal hyper-parameters setting found in \Cref{sec:grid_search} on the \textsc{Baseline} network is not robust enough.
Nonetheless, if we move to the generalization performance, measured by the final test accuracy, things appear differently.
%We report in \Cref{tab:DW_accuracy} and \Cref{tab:RnMob_accuracy} the test accuracy of the different algorithms with both configurations for all \unoref{dataset-network problems}. 
% Summarizing the results in \Cref{tab:WLT}, we can see that, especially for Adam, Adagrad, and Nadam, our tuned configuration leads to better test accuracy, despite the less stable loss profile.

%\unoref{While, looking at the single problems in \Cref{tab:DW_accuracy} and \Cref{tab:RnMob_accuracy}, it is immediately clear that the tuned configuration outperforms the default one in most of the case, this result could still be misleading.}
%\unoref{In the tables we also report the average test accuracy obtained with the five algorithms for each of the six architectures, averaging both across datasets and across architectures.}
%The table is obtained by averaging, for each of the three datasets, the rows of \Cref{tab:DW_accuracy} and \Cref{tab:RnMob_accuracy}. 
\unoref{Indeed, looking at  \Cref{tab:DW_accuracy} and \Cref{tab:RnMob_accuracy}, we observe that SGD presents a strong advantage of the tuned configuration, while average accuracy values are often very close to each other for the other optimizers.}
\unoref{SGD is the only non-adaptive optimizer, meaning that the learning rate is not adjusted during the training.}
\unoref{We argue that this makes SGD much more sensitive to the hyper-parameters setting than other adaptive algorithms.
Nonetheless, even on Adam and Nadam, the tuned configuration achieves slightly better test accuracy.}
% \unoref{Analogously, in \Cref{tab:avg_acc_nets} we average the test accuracy over the algorithms, obtaining different values per each architecture. 
% Of course, on the \textsc{Baseline} architecture for the UC Merced problem (which is the baseline problem) the tuned version is always the best one. 
% Here we can highlight how the tuned configuration is robust to the increase in depth and width with respect to the baseline architecture. The tuned configuration almost always leads to better values, gaining up to a $+10 \%$ test accuracy.}
\unoref{A remarkable exception to this pattern is Resnet50 in \Cref{tab:RnMob_accuracy}, where the default configuration significantly outperforms the tuned one. 
This result seems to suggest that our hyper-parameter configuration found on the baseline is not robust to more radical architectural changes, like in Resnet50, where residual connections (see \Cref{sec:NetArc}) are added to each layer to prevent the vanishing gradient effect.}
%In \citep{probst2019tunability} has been discussed in detail the fact that the default hyper-parameters configuration  is chosen by maximizing the aggregated (usually the average) performance across a variety of tasks  while ensuring sufficiently stable behavior and balancing a trade-off between efficiency and adaptability to different datasets.
%\unoref{Our computational experiments suggest that, when dealing with different problems belonging to a specific class (i.e., in our case, image classification), it may be useful to tune the hyper-parameters once for all on a representative baseline problem.}
%\unoref{The advantage obtained by tuning can play a crucial role in the generalization capability of the model particularly when using SGD, which is one a widespread used algorithm in this field.}
%\unoref{Conversely, our experiments also show that this holds only when the different problem instances are sufficiently close to each other.}
\unoref{Nonetheless, the overall average effect of using the tuned configuration instead of the default is positive, as it is shown in the last row of \Cref{tab:avg_diffs}. Summing up positive and negative contributions, the average improvement in terms of test accuracy is almost $2\%$ with the only exception of Resnet architecture.}

\begin{table*}[h!]
%\todo[inline]{CHECK: per RMSPROP I valori tuned sono peggiori di quelli default anche nella version BASELINE ! Questo non è possible. Ci deve essere errore. Sono stati invertiti ? erano così anche nella tabella orginaria}
    %\footnotesize
    \centering
    \caption{Test accuracy in \% obtained \dueref{by training} the synthetic network architectures using the Default and Tuned hyperparameter values. \dueref{For each algorithm in the two settings, we report the average over the datasets when the architecture is fixed (Avg DATA) and the average over the architecture when the Dataset is fixed (Avg ARCH).}
    %{Here we report the results on the \textsc{Baseline} architecture and its variants described in \Cref{sec:NetArc}}.
    }
    \label{tab:DW_accuracy}
    \renewcommand{\arraystretch}{1.5}
    \resizebox{\textwidth}{!}{
		\begin{tabular}{| l | l | c c  c c | c || c c c c |c|}		\cline{3-12}			
	\multicolumn{2}{c|}{}		&	\multicolumn{5}{c||}{\textsc{Default}} 	&	\multicolumn{5}{c|}{\textsc{Tuned}} \\\hline
	&	DATASET	&	\textsc{Baseline}	&	\textsc{Wide}	&	\textsc{Deep}	&	\textsc{DeepWide}	&		Avg ARCH	&	\textsc{Baseline}	&	\textsc{Wide}	&	\textsc{Deep}	&	\textsc{DeepWide}		&	Avg ARCH		\\\hline
\multirow{4}{*}{\rotatebox{90}{Adam}}	&	UC Merced	&	72,4	&	65,1	&	56	&	61,9	&	\bf	63,9	&	74,9	&	70,3	&	51,7	&	48,6		&	61,4	\\
	&	CIFAR10	&	77,2	&	82,3	&	77,6	&	85,3	&		80,6	&	77,6	&	83,2	&	79,4	&	86,2		& \bf	81,6	\\
	&	CIFAR100	&	48,8	&	54,3	&	46,4	&	50,2	&		49,9	&	48,6	&	54,7	&	47,4	&	49,7		&\bf	50,1	\\ \cline{2-12}
	&	Avg  DATA 	&	66,1	&	67,2	&\bf	60,0	&\bf	65,8	&			&	\bf 67,0	&\bf	69,4	&	59,5	&	61,5		&		\\\hline
\multirow{4}{*}{\rotatebox{90}{Adamax}}	&	UC Merced	&	72,5	&	68,2	&	18,3	&	30,5	&		47,4	&	72,5	&	67,8	&	58,2	&	64,1		&\bf	65,7	\\
	&	CIFAR10	&	77,7	&	82,5	&	78,7	&	86,4	&		\bf 81,3	&	77,8	&	82,7	&	78,9	&	84,9		&	81,1	\\
	&	CIFAR100	&	46,8	&	55	&	47,5	&	55,5	&		51,2	&	47,3	&	55,4	&	46,2	&	55,9		&	51,2	\\ \cline{2-12}
	&	Avg  DATA 	&	65,7	&	68,6	&	48,2	&	57,5	&			&\bf	65,9	&	68,6	&\bf	61,1	&\bf	68,3		&		\\\hline
 \multirow{4}{*}{\rotatebox{90}{Nadam}}	&	UC Merced	&	72,1	&	70,8	&	51,3	&	60,1	&		\bf 63,6	&	73,7	&	73,8	&	39,7	&	59		&	61,6	\\
	&	CIFAR10	&	78,7	&	82,7	&	77	&	85,4	&		81,0	&	78,8	&	82,6	&	79,5	&	86,6		&\bf	81,9	\\
	&	CIFAR100	&	47,1	&	55,3	&	46,9	&	52	&		50,3	&	48,9	&	53,9	&	47,1	&	53,2		&\bf	50,8	\\ \cline{2-12}
	&	Avg  DATA 	&	66,0	&	69,6	&\bf	58,4	&	65,8	&			&\bf	67,1	&\bf	70,1	&	55,4	&\bf	66,3		&		\\\hline
\multirow{4}{*}{\rotatebox{90}{RMSProp}}	&	UC Merced	&	70.4	&	65,2	&	58,9	&	52,5	&	\bf	61,8	&	72,1	&	64,2	&	28,1	&	44,3		&	52,2	\\
	&	CIFAR10	&	77,2	&	82,8	&	78,6	&	84,2	&	\bf	80,7	&	76,1	&	80,2	&	77,6	&	83,7		&	79,4	\\
	&	CIFAR100	&	44,8	&	51,2	&	44,6	&	55,1	&		48,9	&	47,1	&	51,8	&	45,9	&	55,5		&\bf	50,1	\\ \cline{2-12}
	&	Avg  DATA 	&	64,1	&\bf	66,4	&\bf	60,7	&\bf	63,9	&			&	\bf 65,1	&	65,4	&	50,5	&	61,2		&		\\\hline
\multirow{4}{*}{\rotatebox{90}{SGD}}  	&	UC Merced	&	65,1	&	63,2	&	24,3	&	23,3	&		44,0	&	69,5	&	63,2	&	51,4	&	56,8		&\bf	60,2	\\
	&	CIFAR10	&	75,1	&	81,3	&	75,5	&	77	&		77,2	&	75,1	&	80,7	&	74,9	&	82,9		&\bf	78,4	\\
	&	CIFAR100	&	46,9	&	55,1	&	35,8	&	53,9	&		47,9	&	46,9	&	54,9	&	40,6	&	52,5		&\bf	48,7	\\ \cline{2-12}
	&	Avg over DATA 	&	62,4	&\bf	66,5	&	45,2	&	51,4	&			&\bf	63,8	&	66,3	&\bf	55,6	&\bf	64,1		&		\\\hline
    \end{tabular}}
\end{table*}

\begin{table*}[h!]
    %\footnotesize
    \centering    
    \caption{Test accuracy in \% obtained \dueref{by training} the Resnet50 and Mobilenetv2 using the Default and Tuned hyperparameter values. \dueref{For each algorithm in the two settings, we report the average over the datasets when the architecture is fixed (Avg  DATA) and the average over the architecture when the Dataset is fixed (Avg ARCH).}}
%Test accuracy obtained for Resnet50 and Mobilenetv2 with default and tuned hyperparameters values.}
    \label{tab:RnMob_accuracy}
    \renewcommand{\arraystretch}{1.5}
    \resizebox{0.6\textwidth}{!}{\begin{tabular}{ l |l | c c | c || c c | c| }
    \cline{3-8}
   	\multicolumn{2}{c|}{} 		&	\multicolumn{3}{c||}{\textsc{Default}} 					&	\multicolumn{3}{c|}{\textsc{Tuned}} 				\\\cline{2-8}	
	&	DATASET	&	Resnet50	&	Mobilenetv2	&	Avg ARCH	&	Resnet50	&	Mobilenetv2	&	Avg ARCH	\\\hline
\multirow{4}{*}{\rotatebox{90}{Adam}}	&	UC Merced	&	84,8	&	23,3	&	54,1	&	91,9	&	24,6	&	\bf 58,3	\\
	&	CIFAR10	&	93	&	69	&	\bf 81,0	&	74,8	&	64	&	69,4	\\
	&	CIFAR100	&	42,9	&	39	&	41,0	&	45,3	&	40,3	&	\bf 42,8	\\ \cline{2-8}
	&	Avg  DATA 	&	\bf 73,6	&	43,8	&	\bf 58,7	&	\bf 70,7	&	43,0	&	56,8	\\\hline
\multirow{4}{*}{\rotatebox{90}{Adamax}}	&	UC Merced	&	93	&	69	&	\bf 81,0	&	74,8	&	64	&	69,4	\\
	&	CIFAR10	&	82,4	&	81,7	&	\bf 82,1	&	83,3	&	80,6	&	82,0	\\
	&	CIFAR100	&	51,2	&	45,1	&	48,2	&	49	&	47,9	&	\bf 48,5	\\ \cline{2-8}
	&	Avg  DATA 	&	\bf 75,5	&	65,3	&	\bf 70,4	&	\bf 69,0	&	64,2	&	66,6	\\\hline
 \multirow{4}{*}{\rotatebox{90}{Nadam}}	&	UC Merced	&	90,8	&	18,7	&	\bf 54,8	&	86,2	&	15,5	&	50,9	\\
	&	CIFAR10	&	82,8	&	81,2	&	\bf 82,0	&	83,1	&	80,4	&	81,8	\\
	&	CIFAR100	&	44,2	&	35,4	&	\bf 39,8	&	44,3	&	33,4	&	38,9	\\ \cline{2-8}
	&	Avg  DATA 	&	\bf 72,6	&	45,1	&	\bf 58,9	&	\bf 71,2	&	43,1	&	57,2	\\\hline
\multirow{4}{*}{\rotatebox{90}{RMSProp}}	&	UC Merced	&	70,5	&	18,3	&	\bf 44,4	&	78,7	&	8,7	&	43,7	\\
	&	CIFAR10	&	80,7	&	36,2	&	58,5	&	80,3	&	40,5	&	\bf 60,4	\\
	&	CIFAR100	&	25	&	18,5	&	21,8	&	42,6	&	21,9	&	\bf 32,3	\\ \cline{2-8}
	&	Avg  DATA 	&	\bf 58,7	&	24,3	&	41,5	&	\bf 67,2	&	23,7	&	\bf 45,5	\\\hline
\multirow{4}{*}{\rotatebox{90}{SGD}}  	&	UC Merced	&	97	&	69,2	&	83,1	&	74,8	&	97	&	\bf 85,9	\\
	&	CIFAR10	&	81,7	&	70	&	75,9	&	81,7	&	79,9	&	\bf 80,8	\\
	&	CIFAR100	&	51,2	&	48,6	&	49,9	&	51,2	&	48,6	&	49,9	\\ \cline{2-8}
	&	Avg  DATA 	&	\bf 76,6	&	62,6	&	69,6	&	69,2	&	\bf 75,2	&	\bf 72,2	\\\hline
    \end{tabular}}
\end{table*}

\begin{table}[ht!]
\centering
  \caption{Differences in test accuracy (TEST ACC) of Tuned vs Default hyperpameters setting when training the Resnet50, Mobilenetv2 and on average for the Synthetic Networks. Positive (negative) values indicate better performance of the tuned (default) setting.}
    \label{tab:avg_diffs}
 \renewcommand{\arraystretch}{1.5}
    \resizebox{0.6\textwidth}{!}{
\begin{tabular}{l| l | c c c | }
\renewcommand{\arraystretch}{1.5}
	&		&	\multicolumn{3}{c|}{\textsc{TEST ACC TUNED - TEST  ACC DEFAULT}} 		\\				
	&	DATASET	&	Resnet50	&	Mobilenetv2	&	Avg Synthetic Net.	\\\hline	
\multirow{4}{*}{\rotatebox{90}{Adam}}	&	UC Merced	&	\bf 7,10	&	\bf 1,30	&	-2,48	\\	
	&	CIFAR10	&	-18,20	&	-5,00	&\bf 	1,00	\\	
	&	CIFAR100	&\bf 	2,40	&\bf 	1,30	&\bf 	0,18	\\ \cline{2-5}	
	&	Avg over DATA 	&	-2,90	&	-0,80	&	-0,43	\\\hline	
\multirow{4}{*}{\rotatebox{90}{Adamax}}	&	UC Merced	&	-18,20	&	-5,00	&\bf 	18,28	\\	
	&	CIFAR10	&	\bf 0,90	&	-1,10	&	-0,25	\\	
	&	CIFAR100	&	-2,20	&	\bf 2,80	&	0,00	\\ \cline{2-5}	
	&	Avg over DATA 	&	-6,50	&	-1,10	&\bf 	6,01	\\\hline	
 \multirow{4}{*}{\rotatebox{90}{Nadam}}	&	UC Merced	&	-4,60	&	-3,20	&	-2,03	\\	
	&	CIFAR10	&\bf 	0,30	&	-0,80	&\bf 	0,92	\\	
	&	CIFAR100	&	\bf 0,10	&	-2,00	&\bf 	0,45	\\ \cline{2-5}	
	&	Avg over DATA 	&	-1,40	&	-2,00	&	-0,22	\\\hline	
\multirow{4}{*}{\rotatebox{90}{RMSProp}}	&	UC Merced	&\bf 	8,20	&	-9,60	&	-10,48	\\	
	&	CIFAR10	&	-0,40	&\bf 	4,30	&	-1,30	\\	
	&	CIFAR100	&\bf 	17,60	&\bf 	3,40	&\bf 	1,15	\\ \cline{2-5}	
	&	Avg over DATA 	&\bf 	8,47	&	-0,63	&	-3,54	\\\hline	
\multirow{4}{*}{\rotatebox{90}{SGD}}  	&	UC Merced	&	-22,20	&\bf 	27,80	&\bf 	16,25	\\	
	&	CIFAR10	&	0,00	&\bf 	9,90	&\bf 	1,18	\\	
	&	CIFAR100	&	0,00	&	0,00	&	\bf 0,80	\\ \cline{2-5}	
	&	Avg over DATA 	&	-7,40	&\bf 	12,57	&\bf 	6,08	\\\hline	
 \multicolumn{2}{c|}{\textsc{OVERALL AVERAGE}} 		&	-1,95	&\bf 	1,61	&	\bf 1,58	\\\hline
  \end{tabular}}
\end{table}

\begin{figure*}[h!]
   \centering
    \hfill
    \\
    \vskip -0.25truecm
     \centering
    \subfigure[Default HP - \textsc{Baseline}]{\includegraphics[scale=0.35]{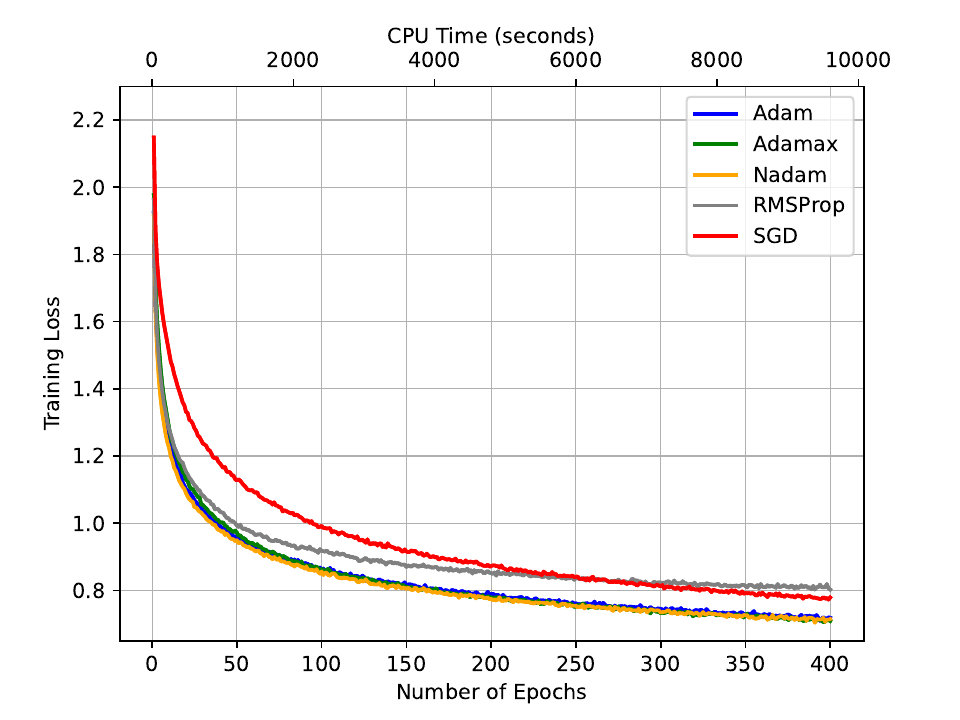}
    } \hfill
    \subfigure[Tuned HP - \textsc{Baseline}]{\includegraphics[scale=0.35]{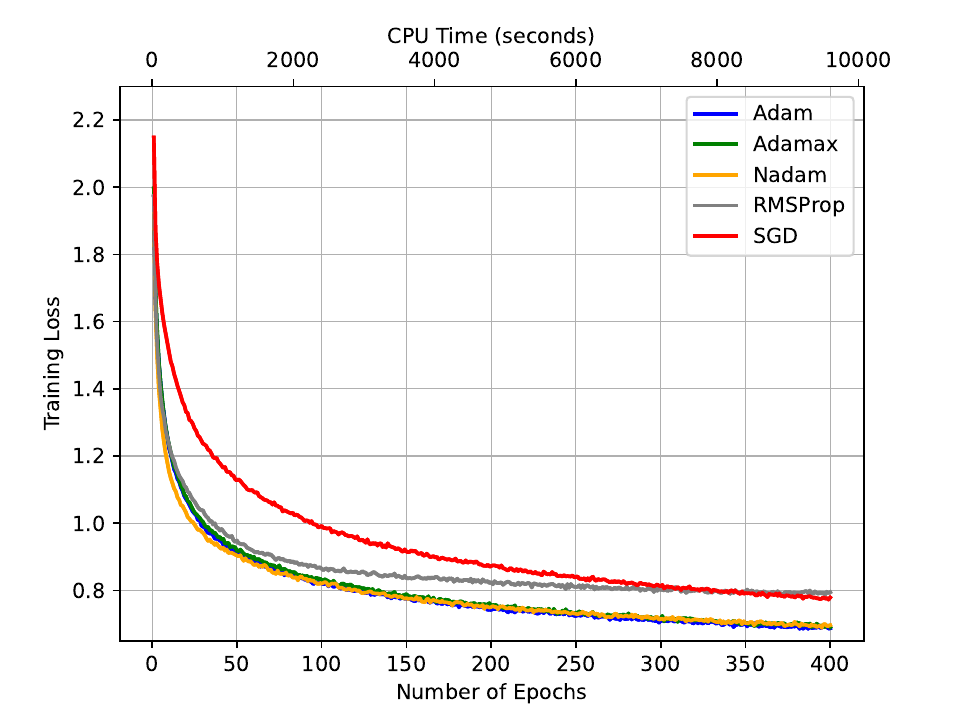}
    } \hfill
 \\
  \vskip -0.25truecm
     \centering
    \subfigure[\textsc{Wide} - Default HP]{\includegraphics[scale=0.35]{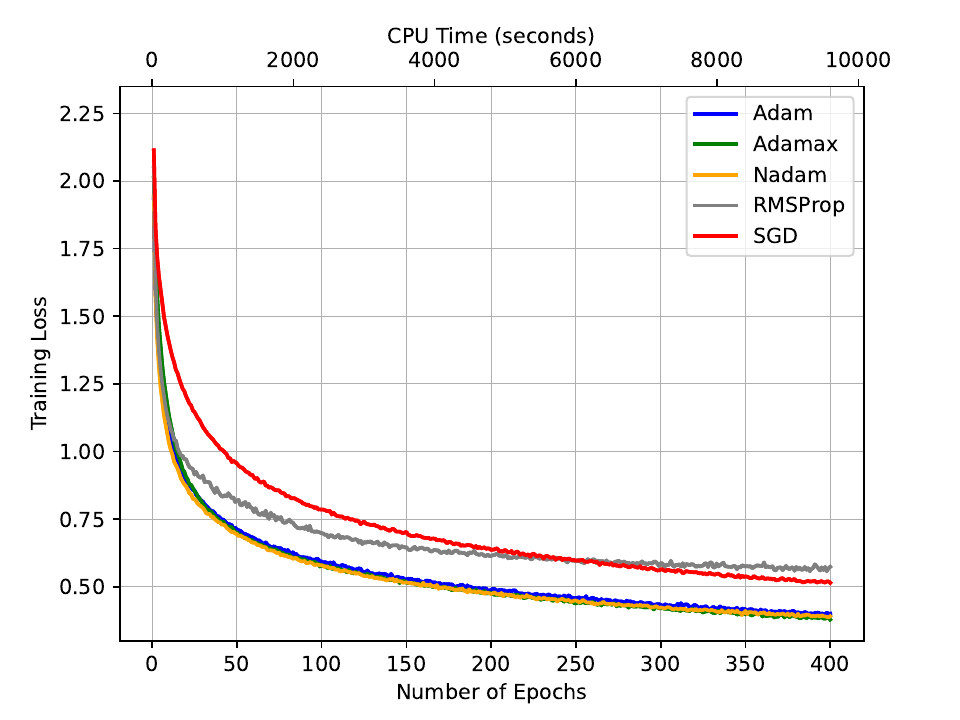}
    } \hfill
    \subfigure[\textsc{Wide} - Tuned HP]{\includegraphics[scale=0.35]{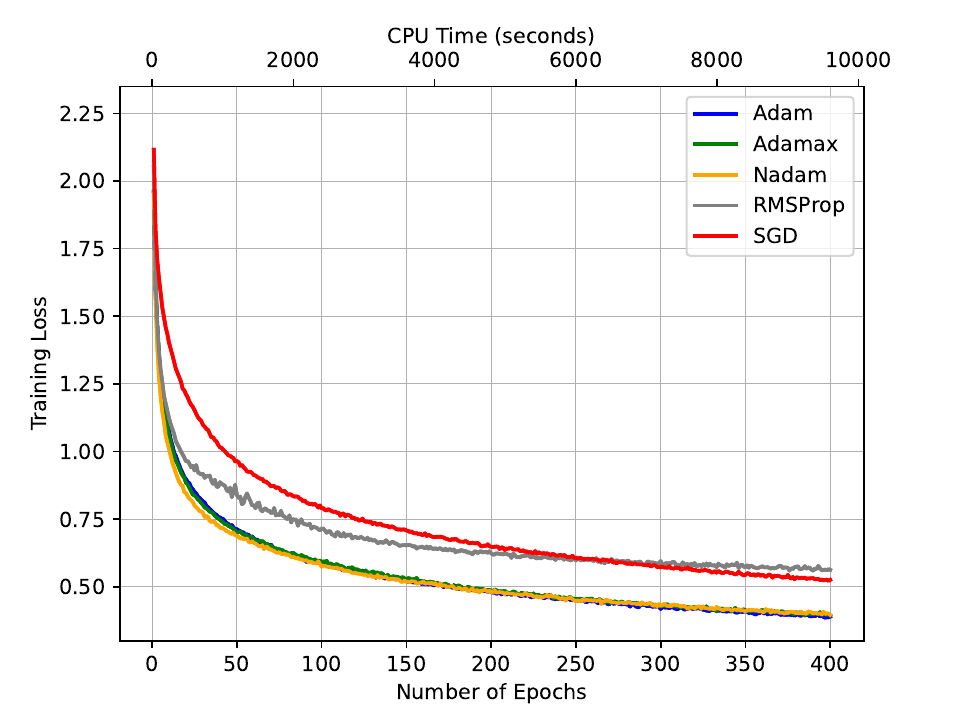}
    } \hfill
 \\
    \vskip -0.25truecm
     \centering
    \subfigure[\textsc{Deep}  -  Default HP]{\includegraphics[scale=0.35]{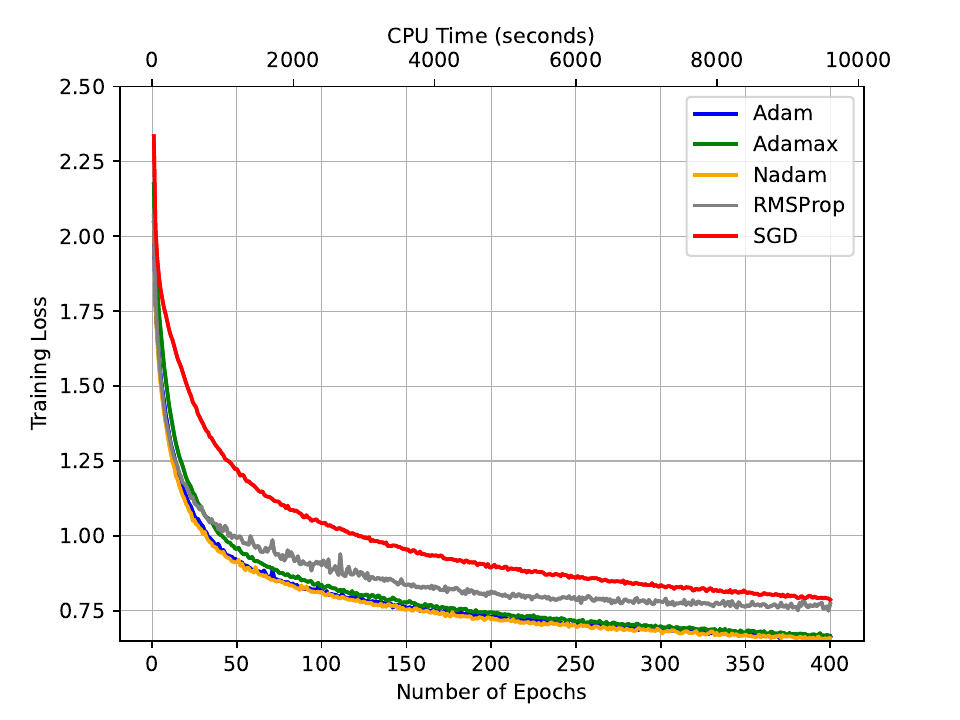}
    } \hfill
    \subfigure[\textsc{Deep}  -  Tuned HP]{\includegraphics[scale=0.35]{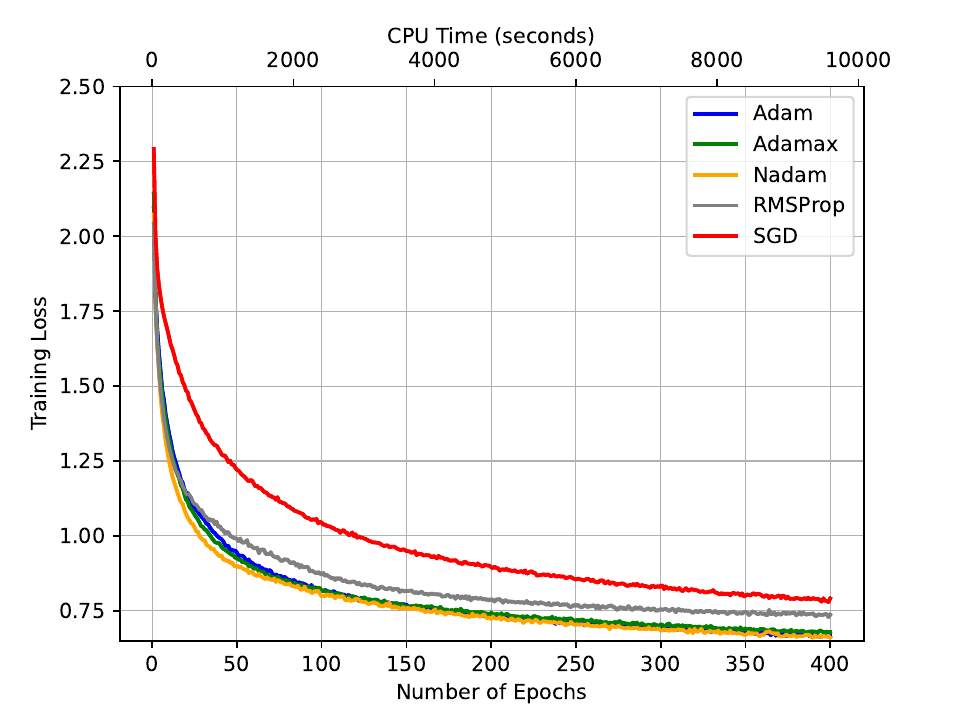}
    } \hfill
\hfill
    \\
    \vskip -0.25truecm
     \centering
    \subfigure[\textsc{Deep}\&\textsc{Wide} - Default HP]{\includegraphics[scale=0.35]{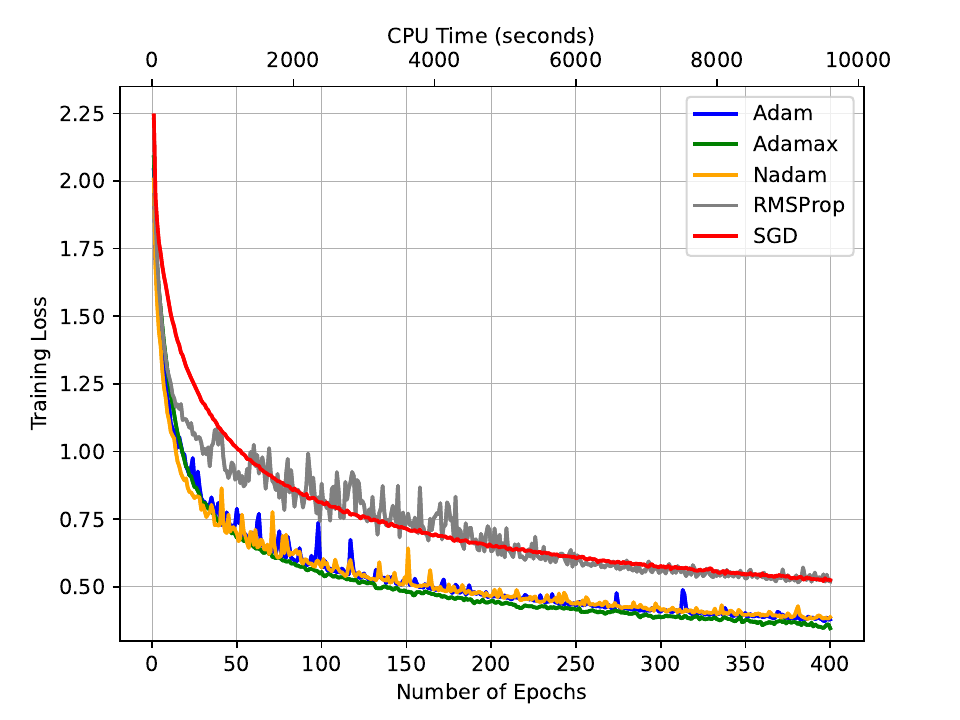}
    } \hfill
    \subfigure[\textsc{Deep}\&\textsc{Wide} - Tuned HP]{\includegraphics[scale=0.35]{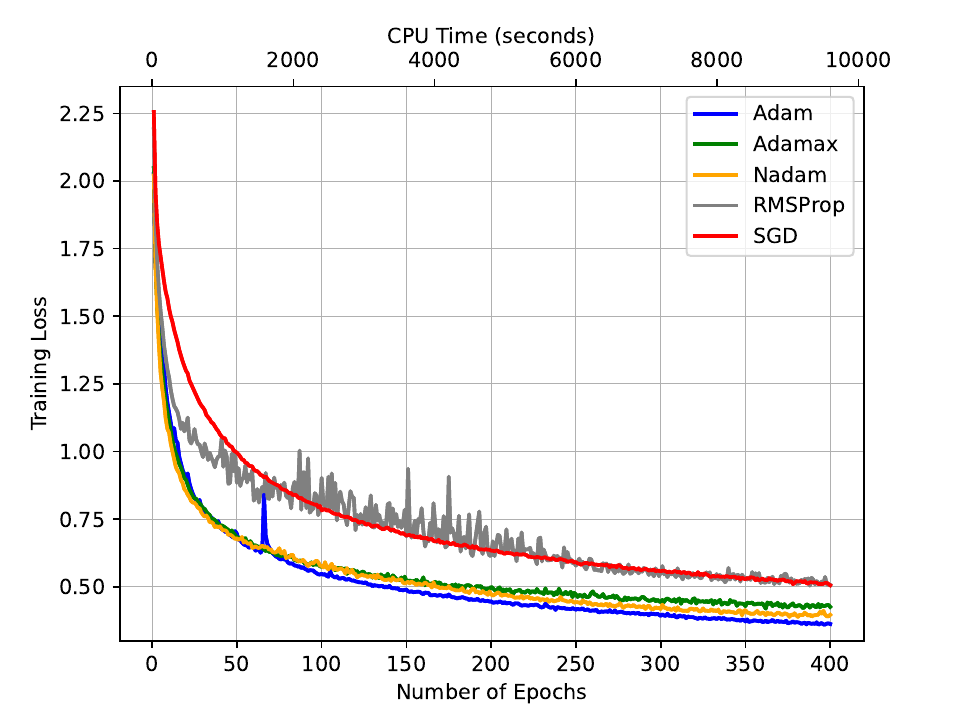}
    } \hfill
    \caption{\dueref{Training} Loss trends \dueref{of the 5 algorithms} for CIFAR10 on \dueref{the synthetic} architectures. 
    %Colours are assigned to optimizers according to the labels above
    }
    \label{fig:cifar10_baseline}
\end{figure*}

\begin{figure*}[h!]
   \centering
    \hfill
    \\
    \vskip -0.25truecm
     \centering
    \subfigure[Default HP - \textsc{Baseline}]{\includegraphics[scale=0.35]{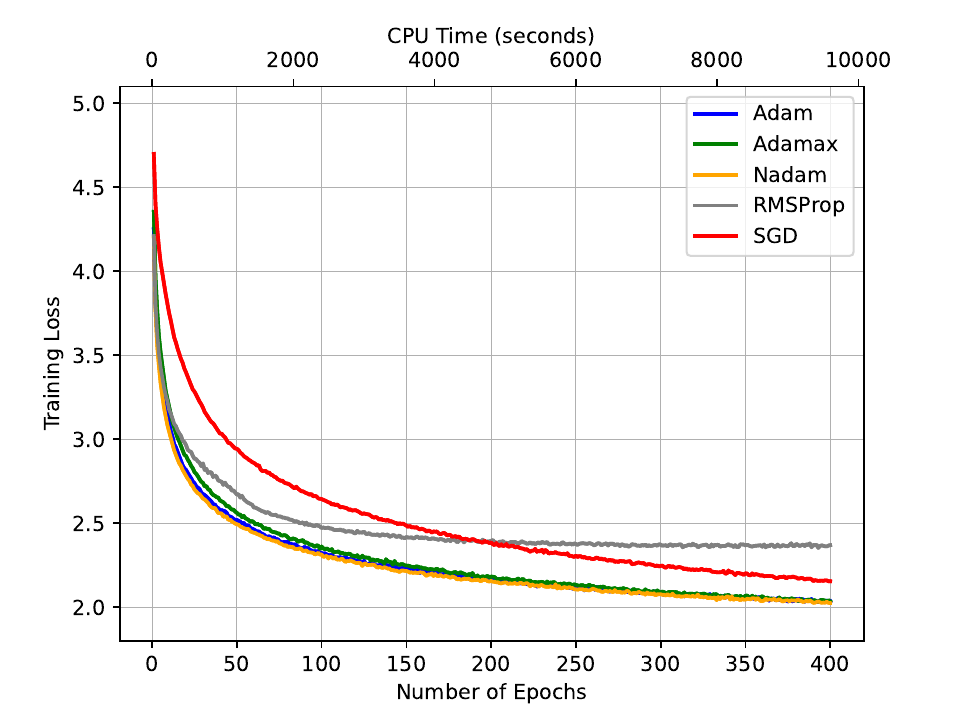}
    } \hfill
    \subfigure[Tuned HP - \textsc{Baseline}]{\includegraphics[scale=0.35]{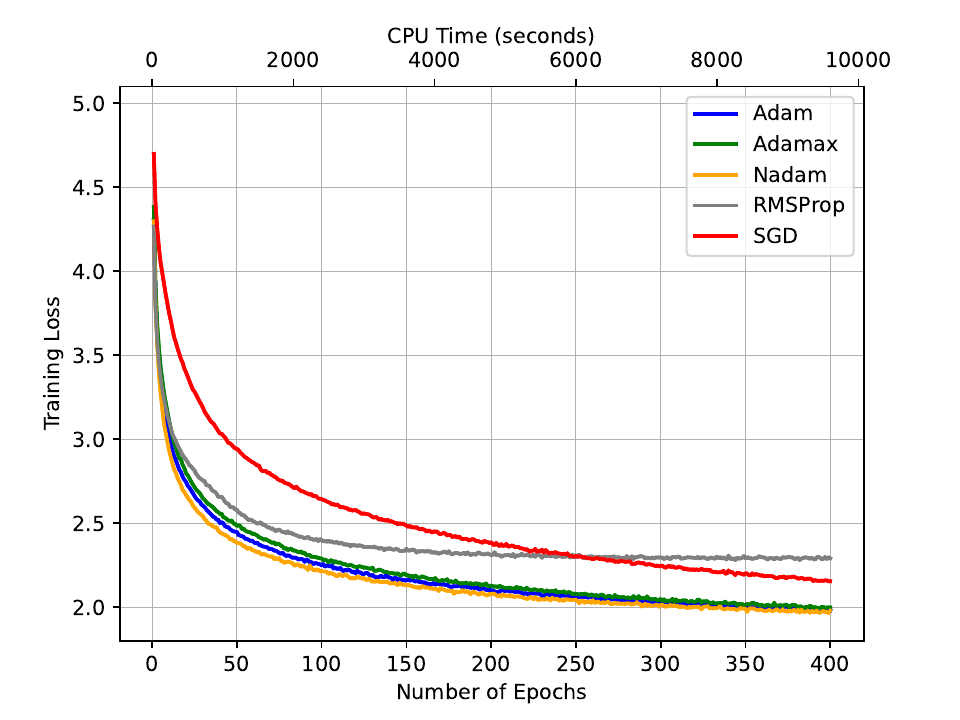}
    } \hfill
 \\
  \vskip -0.25truecm
     \centering
    \subfigure[Default HP - \textsc{Wide}]{\includegraphics[scale=0.35]{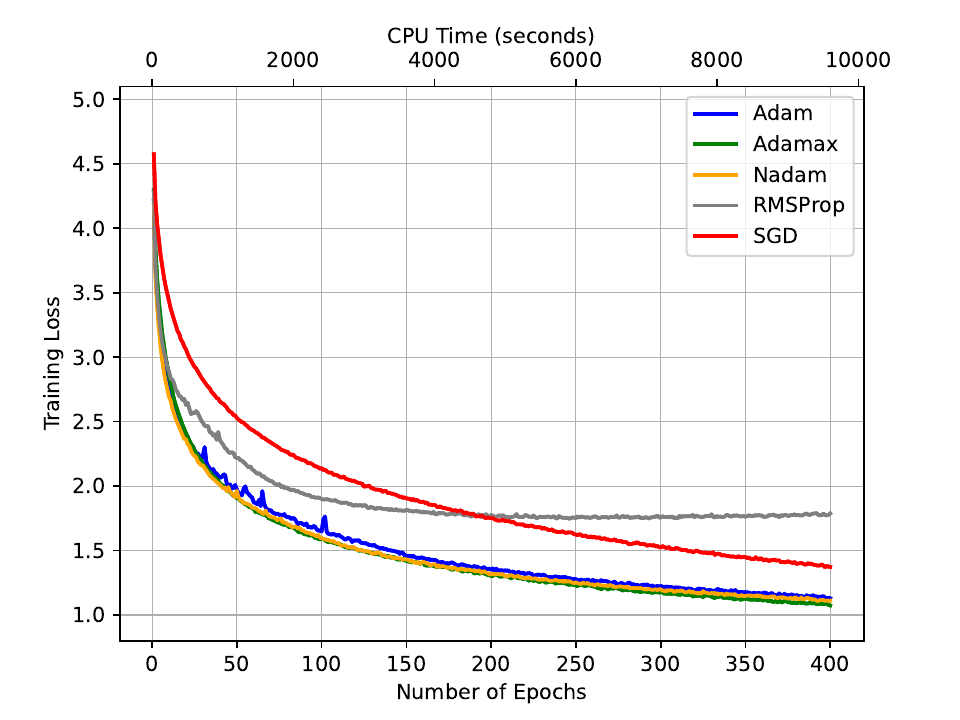}
    } \hfill
    \subfigure[Tuned HP - \textsc{Wide}]{\includegraphics[scale=0.35]{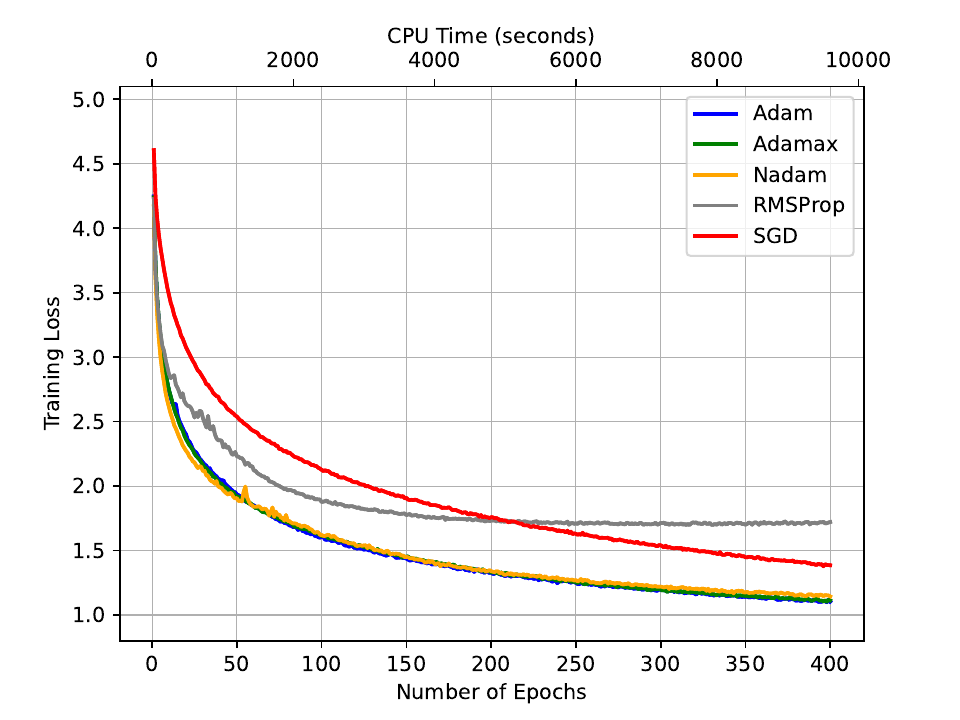}
    } \hfill
 \\
    \vskip -0.25truecm
     \centering
    \subfigure[Default HP - \textsc{Deep}]{\includegraphics[scale=0.35]{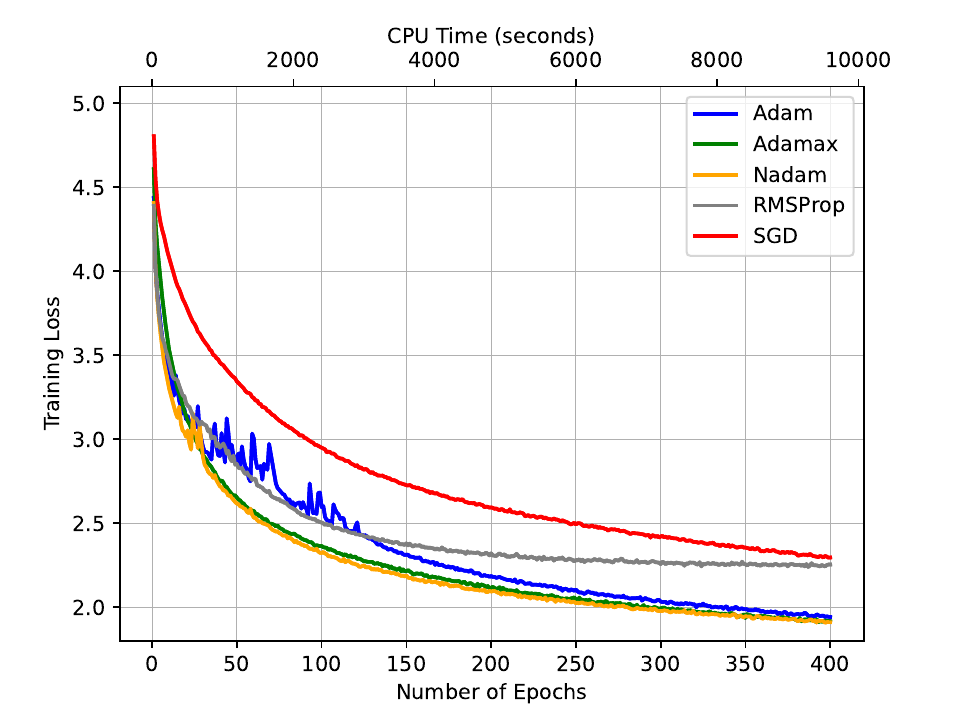}
    } \hfill
    \subfigure[Tuned HP - \textsc{Deep}]{\includegraphics[scale=0.35]{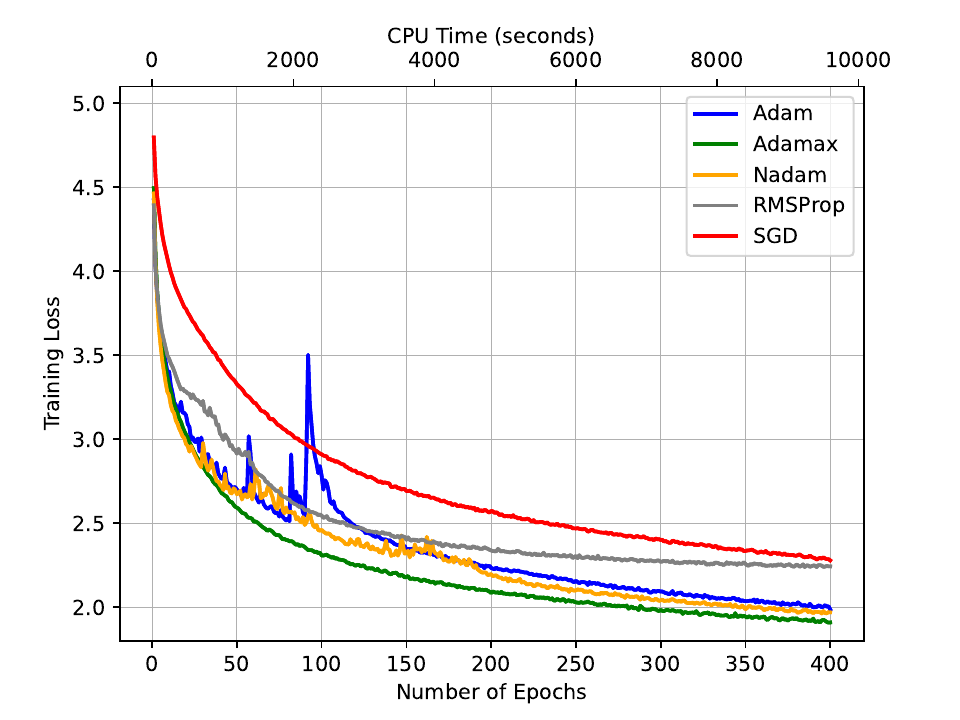}
    } \hfill
\hfill
    \\
    \vskip -0.25truecm
     \centering
    \subfigure[Default HP -  \textsc{Deep}\&\textsc{Wide}]{\includegraphics[scale=0.35]{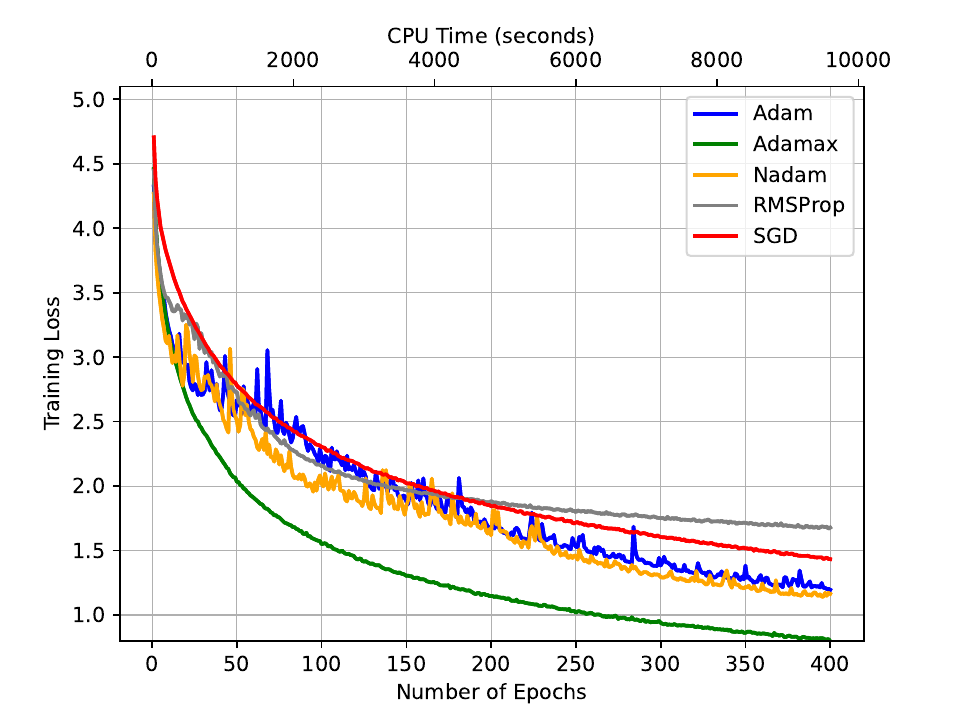}
    } \hfill
    \subfigure[Tuned HP -  \textsc{Deep}\&\textsc{Wide}]{\includegraphics[scale=0.35]{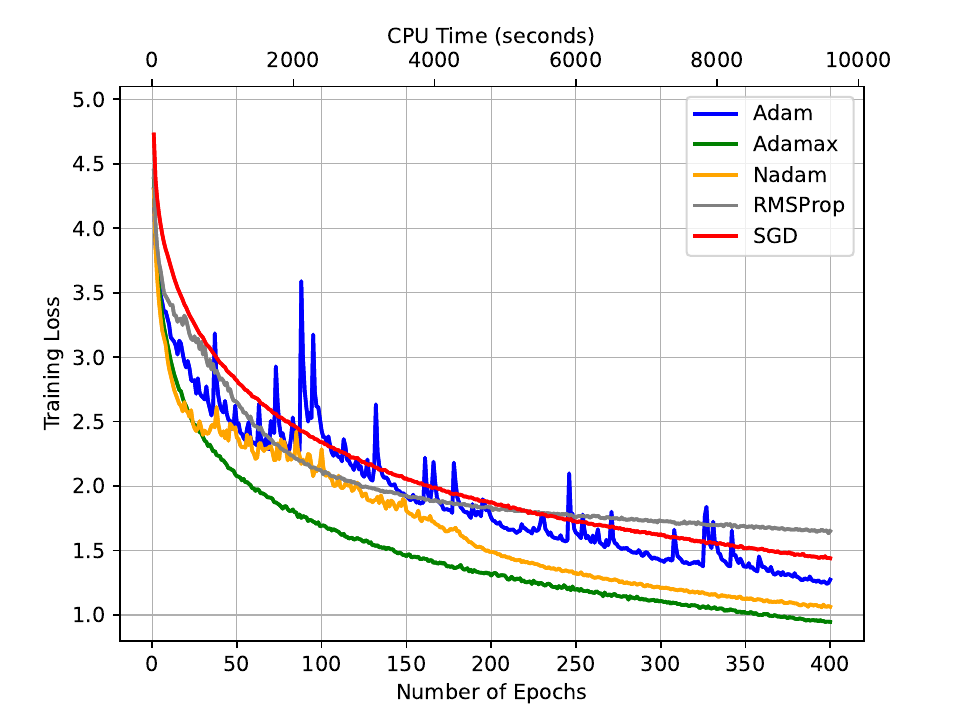}
    } \hfill
    \caption{\dueref{Training} Loss trends \dueref{of the 5 algorithms} for CIFAR100 on \dueref{the synthetic} architectures. }
    \label{fig:CIFAR100_baseline}
\end{figure*}
%%%%%

\begin{figure*}[h!]
   \centering
    \hfill
    \\
    \vskip -0.25truecm
     \centering
    \subfigure[Resnet50 - Default HP]{\includegraphics[scale=0.35]{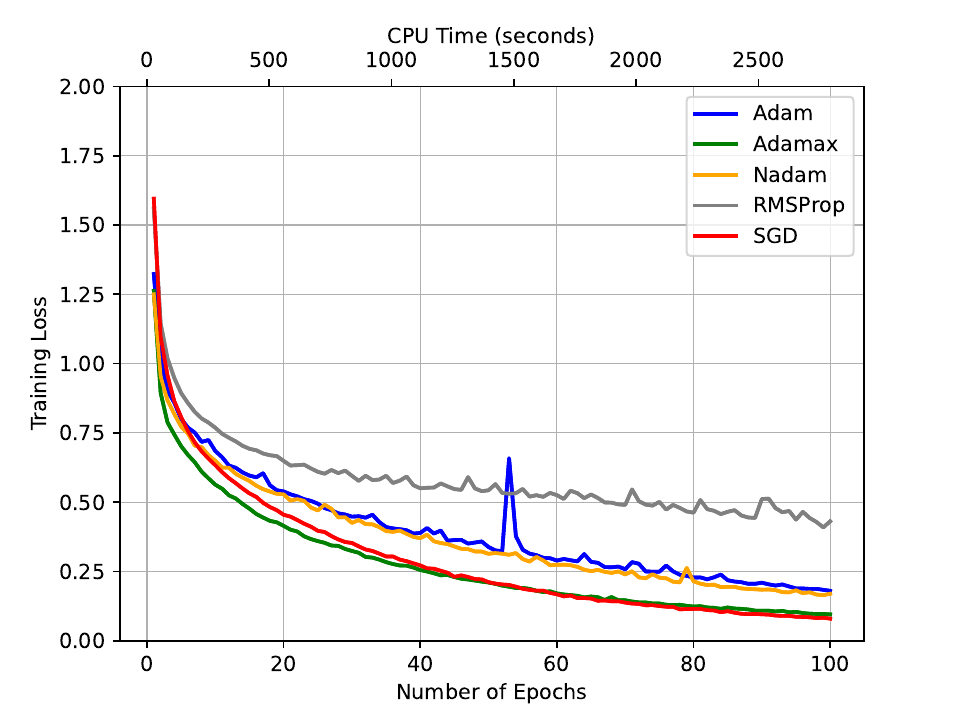}
    } \hfill
    \subfigure[Resnet50 - Tuned HP]{\includegraphics[scale=0.35]{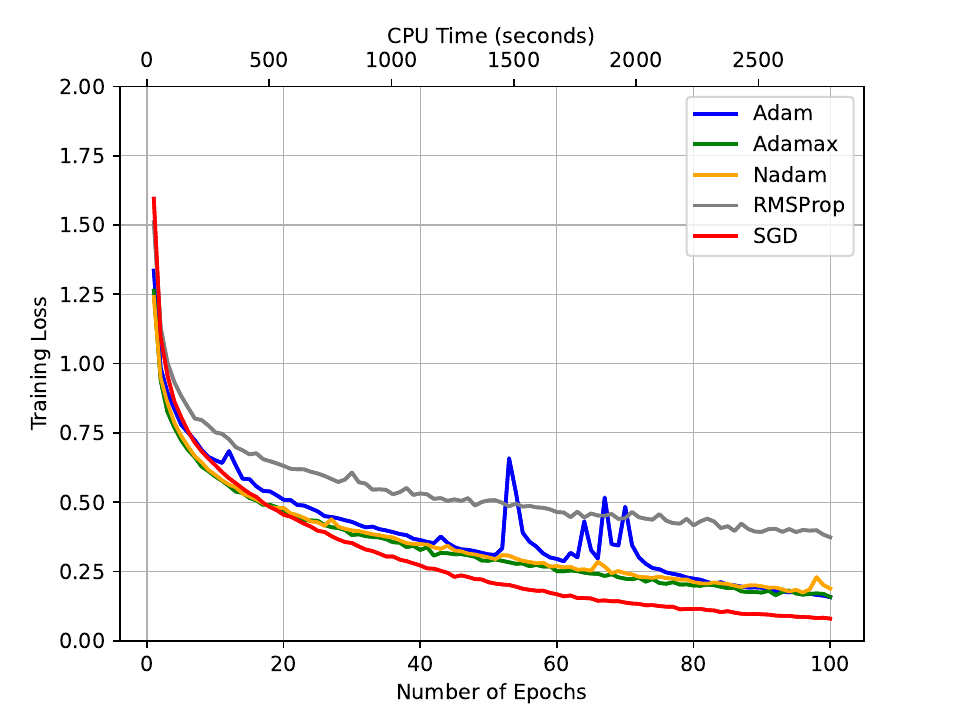}
    } \hfill
 \\
  \vskip -0.25truecm
     \centering
    \subfigure[Mobilenetv2 - Default HP]{\includegraphics[scale=0.35]{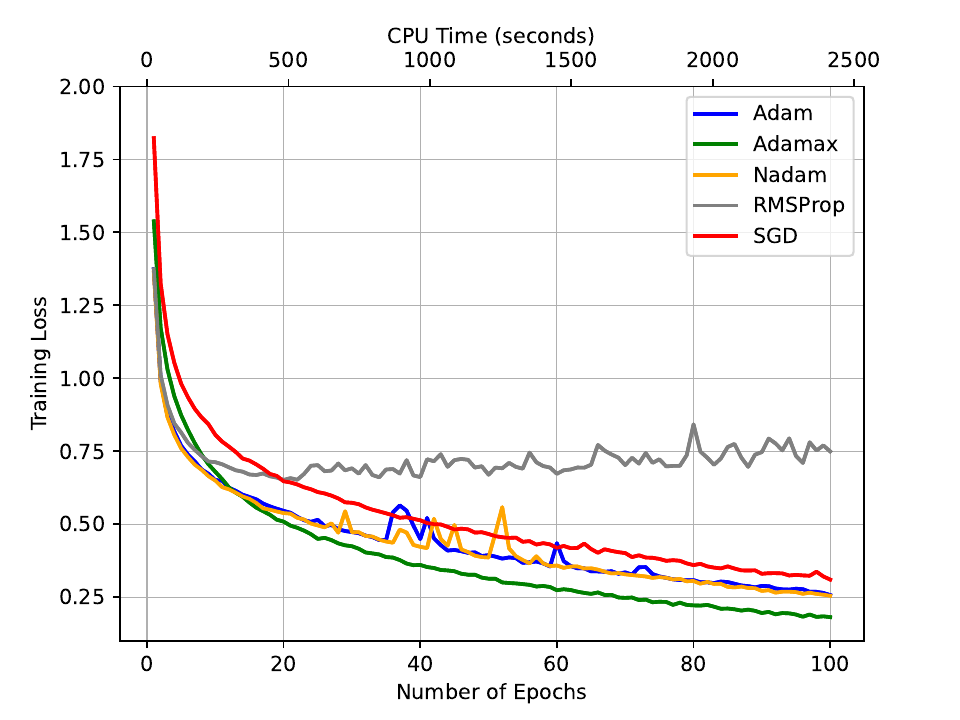}
    } \hfill
    \subfigure[Mobilenetv2 - Tuned HP]{\includegraphics[scale=0.35]{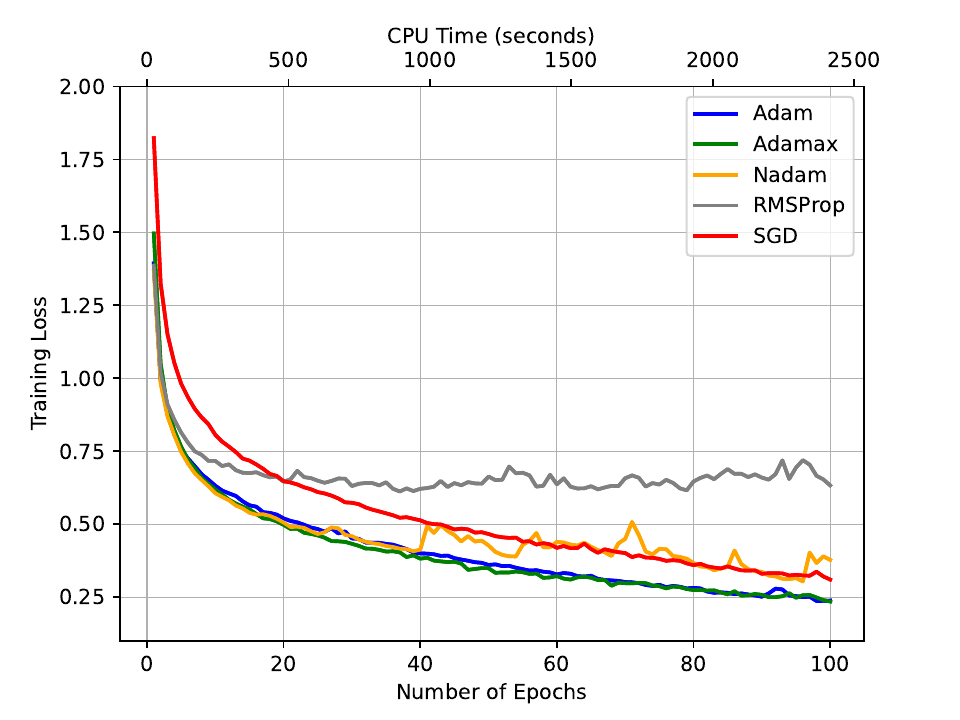}
    } \hfill
    \caption{\dueref{Training} loss trends for CIFAR10 on Resnet50 and MobilenetV2.}
    \label{fig:cifar10_ResMb}
\end{figure*}

\begin{figure*}[h!]
   \centering
    \hfill
    \\
    \vskip -0.25truecm
     \centering
    \subfigure[Resnet50 - Default HP]{\includegraphics[scale=0.35]{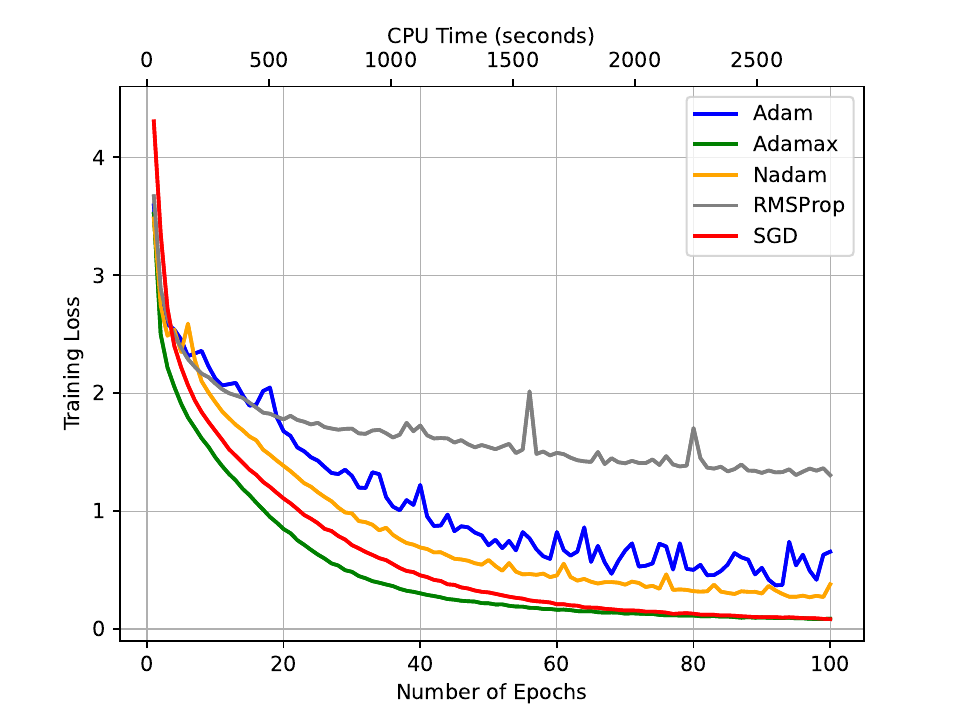}
    } \hfill
    \subfigure[Resnet50 - Tuned HP]{\includegraphics[scale=0.35]{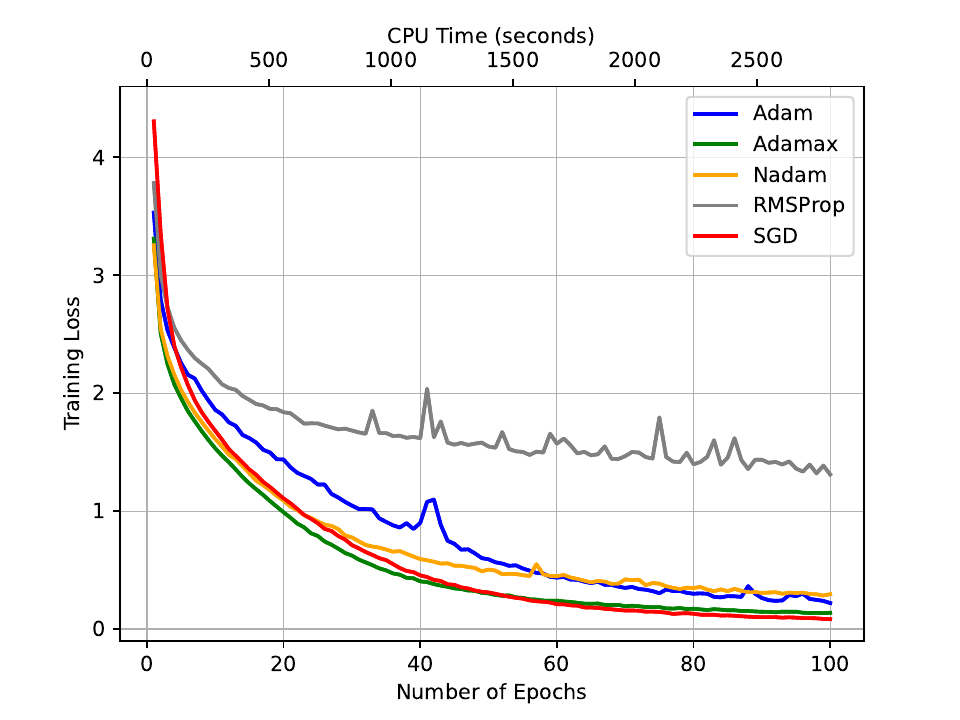}
    } \hfill
 \\
  \vskip -0.25truecm
     \centering
    \subfigure[Mobilenetv2 - Default HP]{\includegraphics[scale=0.35]{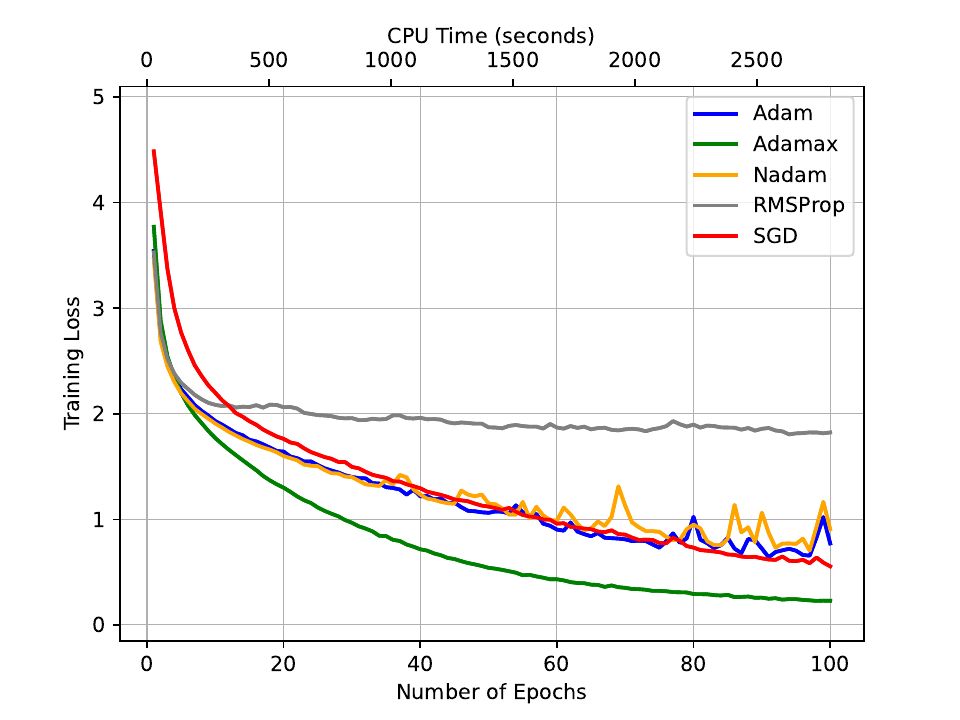}
    } \hfill
    \subfigure[Mobilenetv2 - Tuned HP]{\includegraphics[scale=0.35]{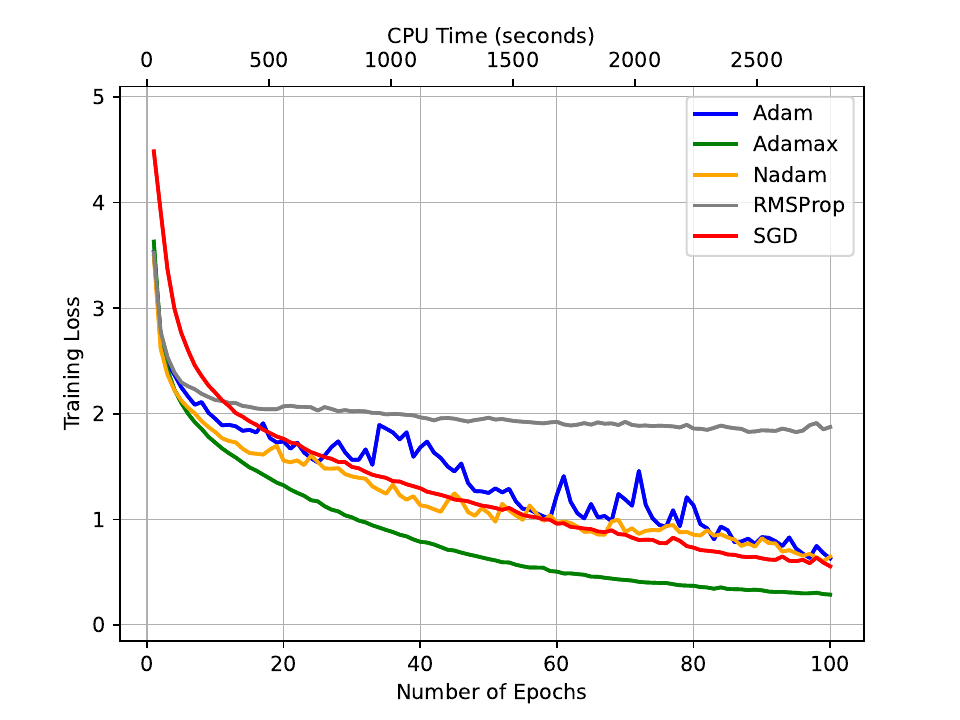}
    } \hfill
    \caption{\dueref{Training} Loss trends \dueref{of the 5 algorithms} for CIFAR100 on Resnet50 and MobilenetV2.}
    \label{fig:cifar100_ResMb}
\end{figure*}

%\todo[inline]{cambia titolo per evidenziare Training results}
\subsection{Collective impact of tuning: performance profiles}
\label{subsec:ppfr}

\dueref{In this section, we consider a collective representation of the training results with the aim of assessing the impact of tuned versus default setting in reaching the global solution.}
Following the underlying idea of the bench-marking method proposed in \citep{dolan2002benchmarking} and \citep{more2009benchmarking}, we consider a variant of the performance profiles as an additional tool of comparison between the five algorithms: Adam, Adamax, Nadam, RMSProp, and SGD.

Following \citep{more2009benchmarking}, given the set of problems $\mathcal{P}$ and the set of solvers $\mathcal{S}$,
\unoref{a problem $p\in \mathcal{P}$ is solved by a solver $s\in \mathcal{S}$ with precision $\tau$ if 
$$f(\omega_{p,s}) \le f^p_L + \tau (f(\omega^0_p) - f^p_L),$$ being
 $\omega^0_p$ the starting point of problem $p \in \mathcal{P}$, which is the same for all solvers,
$f(\omega_{p,s})$ the final value of the objective function in \eqref{eq:loss} after the training process and 
$f^p_L = \min_{s \in \mathcal{S}} f(\omega_{p,s})$.
}
In our case $\mathcal{P}$ {is made up of the 18}  different versions of the problem \eqref{eq:loss} corresponding to all the possible combinations \unoref{of the 6 network architectures} (\textsc{Baseline}, \textsc{Deep}, \textsc{Wide},  \textsc{Deep}\&\textsc{Wide}, Resnet50, Mobilenetv2), with \unoref{the three datasets}. \dueref{Further, in our set of experiments we have fixed the epochs, i.e. the computational time. Hence differently from \citep{more2009benchmarking}, we are interested in checking how many problems are solved to a given accuracy $\tau$. Thus we introduce the}
% \todo[inline]{DOMANDA: lo starting point $\omega^0_p$ cambia con l'architettura perché i parametri da settare sono diversi, ma NON dipende dal dataset. QUindi ad architettura fissata, il punto $\omega^0_p$ è lo stesso. Il valore $f(\omega^0_p)$ cambia invece con il dataset. Quindi fissato architettura e dataset, il valore $f(\omega^0_p)$ dovrebbe essere lo stesso per tutti}
%\remove{we define $f(\omega_{p,s})$ the final value of \Cref{eq:loss} after the training process.  In our case $\mathcal{P}$ includes the different versions of the problem of minimizing \Cref{eq:loss} corresponding to the 24 possible combinations (\textsc{Baseline}, \textsc{Deep}, \textsc{Wide},  \textsc{Deep}\&\textsc{Wide}, Resnet50, Mobilenetv2, with our   without data augmentation, before \change{or}{and} after the grid search). We define $\omega^0_p$ the starting point of problem $p \in \mathcal{P}$, which is the same for all solvers, as the seeds have been set to the same value (see \Cref{sec:impdet}). Given $f^p_L = \min_{s \in \mathcal{S}} f(\omega_{p,s})$, we consider problem $p$ to be solved by solver $s$ with precision $\tau$ if $$f(\omega_{p,s}) \le f^p_L + \tau (f(\omega^0_p) - f^p_L).$$ }
 success rate performance profile $\sigma_s ( \tau)$  for a solver $s$ in \Cref{fig:performanceprofiles} as:
$$ \sigma_s ( \tau) = \frac{1}{\vert \mathcal{P} \vert} \vert \{ p \in P: f(\omega_{p,s}) \le f^p_L + \tau (f(\omega^0_p) - f^p_L)\} \vert .$$ 

\dueref{
In \Cref{fig:performanceprofiles} we plot $ \sigma_s ( \tau)$ with $\tau\in [10^{-4},1]$. The higher the plot on the left, the better.
Looking at the performance profiles in \Cref{fig:performanceprofiles}, we confirm that it is not possible to state the superiority of the tuning versions in the training performance. However, the tuned versions of the algorithms differ less from each other, being more stable. 
}
%can conclude that \unoref{in the default configuration} Adamax significantly outperforms the other algorithms on the task in terms of efficiency, \unoref{while in the tuned configuration, the algorithms seem to have more similar performances}.
%\unoref{On the one hand we} remark that, to the best of our knowledge, at least two computational studies confirm that Adamax, together with Adam and Nadam, seems to be the best-performing optimizer for image classification tasks (\cite{gupta2021adam,bera2020analysis}).
%\unoref{On the other hand, we observe that the optimizers are more stable with the tuned configuration, meaning they are almost interchangeable.}
%\unoref{It is remarkable, for instance, how Adam performs significantly better in the tuned configuration even for high-precision $\tau$.}
%\todo[inline]{NOn capisco questa affernazione. A me NADAM mi sembra sotto ADAM. Come si leggono questi performance profile ? Più in alto a sx e meglio è ?}
\unoref{
\subsection{Data availability}
All the data we have presented in this section are fully reproducible from the source code, which is available on the public Github repository at \url{https://github.com/lorenzopapa5/Computational_Issues_in_Optimization_for_Deep_networks}.}

% \begin{figure*}
%    \centering
%     %\includegraphics[width=0.5\linewidth]{labels5.pdf} 
%     \includegraphics[scale=0.5]{PerformanceProfiles_all.pdf}
%     \caption{Success rate performance profiles \unoref{on the 24 problems arising by considering} all the architectures \unoref{with all the hyper parameters settings} over all the three datasets.}
%     \label{fig:performanceprofiles}
% \end{figure*}

\begin{figure*}[h!]
%\todo[inline]{è possibile togliere la scritta in testa alla figura ? per di più c'è scritto optimal e no tuned, ma è una ripetizione}
    \centering
    \hfill
    \\
    \vskip -0.5truecm
    \centering
    \subfigure[Default HP]{\includegraphics[scale=0.35]{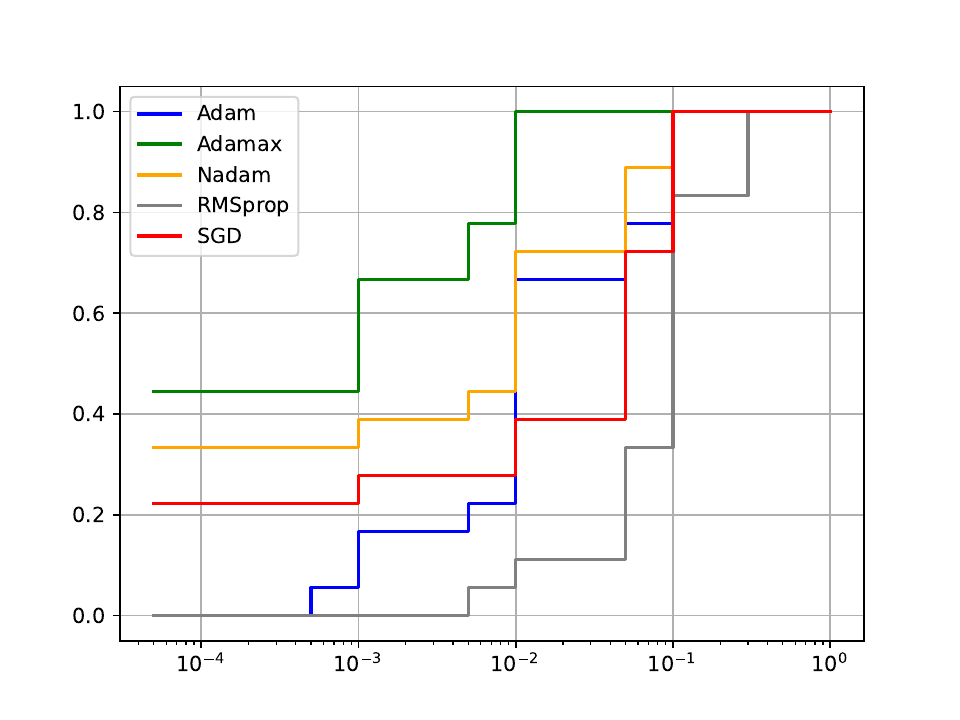}
    } \hfill
    \subfigure[Tuned HP]{\includegraphics[scale=0.35]{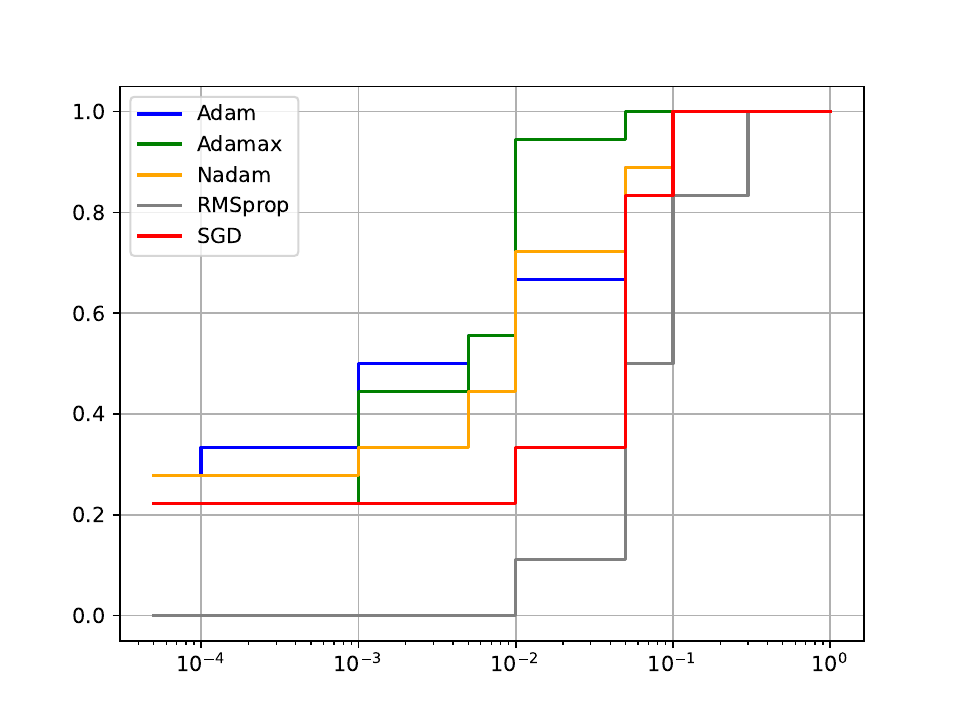}
    } \hfill
    \caption{Success rate performance profiles $\sigma_s ( \tau) $ of the five algorithms $s$ \unoref{on the 18 problems arising by considering} all the 6 architectures over all the three datasets: a) \unoref{with default} and b) {tuned hyper parameters settings} }
    \label{fig:performanceprofiles}
\end{figure*}

\clearpage
\newpage
\section{Conclusions}
\label{sec:Conc}
In this paper, nine optimization open-source algorithms have been extensively tested in training a deep CNN network on a multi-class classification task.
%The first set of tests analyses convergence to a global minimizer by embedding the local algorithm in a multistart framework. 
Computational experience shows that not all the algorithms reach a neighbourhood of a global solution, and some of them get stuck in local minima, independently of the choice of the starting point. Algorithms reaching a local non-global solution have {test performances}{, i.e., accuracy on the test set,} far below the minimal required threshold for such a task. 
This result confirms the initial claim that reaching a neighbourhood of a global optimum is extremely important for generalization performance.

A fine grid search on the optimization hyperparameters leads to hyperparameter choices that give remarkable improvements in test accuracy when the network structure and the dataset do not change. 
Thus, using a default setting might not be the better choice.

Finally, the tests on different architectures {and datasets} suggest that \unoref{when the architectural changes are not too radical, it could be  convenient to use the tuned configuration than the default one.}
\unoref{We believe that this result can have a remarkable impact, especially on ML practitioners asked to train similar models on different datasets belonging to the same problem class, e.g., image classification, which is often the case in real-world applications.}
\unoref{Performing a grid-search on a representative problem of a given class and tuning the hyperparameters on it instead of using the default configuration can be seen as creating a new customized setting, which is reusable for larger instances and achieves better generalization performance.}
% Instead, increasing the depth of the network seems to make the optimization task more difficult both with standard and refined hyperparameters settings.

\section*{Acknowledgements}
 Authors thank the ALCOR laboratory of DIAG Sapienza University of Rome (\url{https://alcorlab.diag.uniroma1.it/}) for making the workstations available for the tests. We also thank Nicolas Zaccaria for the extensive editing performed on the graphs.
Reviewers' comments helped us to improve significantly the paper. 
Laura Palagi acknowledges financial support from Progetto di Ricerca Medio Sapienza Uniroma1 (2022) - n. RM1221816BAE8A79.
Corrado Coppola acknowledges financial support from Progetto Avvio alla Ricerca Sapienza Uniroma1 (2023) - n. AR123188B03A3356.

\clearpage
\newpage
\setlength{\bibsep}{3pt}
\renewcommand{\bibfont}{\small}
\bibliography{Manuscript}

\end{document}